\documentclass[11pt,a4paper]{article}
\usepackage[T1]{fontenc}
\usepackage{lmodern}
\usepackage[utf8]{inputenc}
\usepackage{indentfirst}

\usepackage{amsmath,amsthm,mathtools}
\usepackage{amssymb}
\usepackage[margin=25mm,headheight=14pt]{geometry}
\usepackage[expansion=false]{microtype}
\usepackage{booktabs,array,longtable}
\usepackage{enumitem}
\usepackage{xcolor}
\usepackage{hyperref}
\hypersetup{colorlinks=true,linkcolor=black,citecolor=black,urlcolor=blue!45!black,
 pdftitle={Gaussian limits for Rademacher random multiplicative functions from short to near full intervals},
 pdfauthor={},pdfsubject={Random multiplicative functions and short intervals}}
\usepackage{aliascnt}
\usepackage[nameinlink,noabbrev]{cleveref}
\numberwithin{equation}{section}
\allowdisplaybreaks[2]
\theoremstyle{plain}
\newtheorem{theorem}{Theorem}[section]
\theoremstyle{plain}
\newaliascnt{proposition}{theorem}
\newtheorem{proposition}[proposition]{Proposition}
\aliascntresetthe{proposition}
\crefname{proposition}{Proposition}{Propositions}
\Crefname{proposition}{Proposition}{Propositions}
\newaliascnt{lemma}{theorem}
\newtheorem{lemma}[lemma]{Lemma}
\aliascntresetthe{lemma}
\crefname{lemma}{Lemma}{Lemmas}
\Crefname{lemma}{Lemma}{Lemmas}
\newaliascnt{corollary}{theorem}
\newtheorem{corollary}[corollary]{Corollary}
\aliascntresetthe{corollary}
\crefname{corollary}{Corollary}{Corollaries}
\Crefname{corollary}{Corollary}{Corollaries}
\newaliascnt{extinput}{theorem}

\aliascntresetthe{extinput}
\crefname{extinput}{External input}{External inputs}
\Crefname{extinput}{External input}{External inputs}
\theoremstyle{definition}
\newaliascnt{definition}{theorem}

\aliascntresetthe{definition}
\crefname{definition}{Definition}{Definitions}
\Crefname{definition}{Definition}{Definitions}
\newtheoremstyle{boldremark}
  {.5\topsep}{.5\topsep}
  {\normalfont}{}{\bfseries}{.}{.5em}{}
\theoremstyle{boldremark}
\newaliascnt{remark}{theorem}
\newtheorem{remark}[remark]{Remark}
\aliascntresetthe{remark}
\crefname{remark}{Remark}{Remarks}
\Crefname{remark}{Remark}{Remarks}
\newcommand{\EE}{\mathbb E}
\newcommand{\PP}{\mathbb P}
\newcommand{\RR}{\mathbb R}
\newcommand{\NN}{\mathbb N}
\newcommand{\ind}{\mathbf 1}
\newcommand{\cF}{\mathcal F}

\newcommand{\cH}{\mathcal H}
\newcommand{\cQ}{\mathcal Q}
\newcommand{\cE}{\mathcal E}

\newcommand{\eps}{\varepsilon}
\newcommand{\dd}{\,d}
\DeclareMathOperator{\Var}{Var}
\DeclareMathOperator{\Cov}{Cov}

\newcommand{\vR}{v_R(x,y)}
\newcommand{\Vf}{V_f(x,y)}

\title{Distribution of Rademacher random multiplicative functions in short intervals, with proper normalization}
\author{%
  Renjie Zhu\textsuperscript{a} and
  Tianping Zhang{%
    \thanks{%
      Corresponding author: Tianping Zhang.\newline
      E-mail addresses: \texttt{zhurenjie@snnu.edu.cn} (Renjie Zhu),
      \texttt{tpzhang@snnu.edu.cn} (Tianping Zhang).%
    }%
  }\\[0.8em]
  {\scriptsize
   \textsuperscript{a}School of Mathematics and Statistics,
   Shaanxi Normal University, Xi'an, Shaanxi 710119, P.\ R.\ China}\\[0.25em]
  {\scriptsize
   \textsuperscript{b}Research Center for Number Theory and Its Applications,
   Northwest University, Xi'an, Shaanxi 710127, P.\ R.\ China}%
}
\date{}
\renewenvironment{abstract}
  {\par\smallskip\noindent\textbf{Abstract.}\enspace\ignorespaces}
  {\par\smallskip}
\begin{document}
\maketitle
\vspace{-16pt}
\begin{abstract}
We determine the limiting distribution of partial sums $\sum_{x<n\le x+y}f(n)$ of a Rademacher random multiplicative function over short intervals. Combined with the theorem of Soundararajan and Xu, our result gives a standard real Gaussian limit, with a suitable deterministic normalization, throughout $x^{1/5}\log x\ll y=o(x)$. This resolves the Rademacher extension raised by Harper, Soundararajan and Xu. The normalization is smaller than the standard deviation by a factor tending to zero when $\log(x/y)=o(\sqrt{\log\log x})$.
As an application of our method, we also obtain a complex Gaussian limit for short interval sums twisted by a fixed primitive nonreal Dirichlet character.
\end{abstract}
\smallskip


\section{Introduction}\label{sec:intro}

The limiting distributions of partial sums of random multiplicative functions have attracted much attention in recent years. A Steinhaus random multiplicative function $f:\NN\to\mathbb C$ is obtained by choosing the variables $f(p)$ independently and uniformly on the complex unit circle and defining $f(n)=\prod_{p^a\parallel n}f(p)^a$. A Rademacher random multiplicative function is constructed from independent random variables $(\varepsilon_p)_p$, each taking the values $1$ and $-1$ with probability $1/2$, by setting $f(n)=\mu^2(n)\prod_{p\mid n}\varepsilon_p$.

A striking result of Harper \cite{Harper}, resolving a conjecture of Helson \cite{Helson}, is that in both models
\[
 \EE\left|\sum_{n\le x}f(n)\right|
 \asymp\frac{\sqrt{x}}{(\log\log x)^{1/4}}.
\]
Thus the long sum exhibits better than square root cancellation, and division by $\sqrt{x}$ gives convergence to zero in probability. A natural question is whether short interval sums have a Gaussian limit when $y\to\infty$ and $y=o(x)$. 

We write $$S_f(x,y):=\sum_{x<n\le x+y}f(n).$$
Chatterjee and Soundararajan \cite{CS} proved a Gaussian limit in the Rademacher model when $x^{1/5}\log x\ll y=o(x/\log x)$, after normalization by the standard deviation. Later, Soundararajan and Xu \cite{SX} extended this range and established the corresponding Steinhaus result. For Steinhaus functions, their Corollary~1.2 yields the standard complex Gaussian limit
\begin{equation*}
 \frac{S_f(x,y)}{\sqrt y}\xrightarrow{d}\mathcal N_{\mathbb C}(0,1),
\end{equation*}
where $y\to\infty$ and $y\le x/(\log x)^{2\log2-1+\varepsilon}$. For Rademacher functions, their Corollary~9.2 gives
\begin{equation*}
 \frac{S_f(x,y)}{\sqrt{(6/\pi^2)y}}\xrightarrow{d}N(0,1),
\end{equation*}
where $x^{1/5}\log x\ll y\le x/(\log x)^{2\log2-1-\varepsilon}$. Here $\varepsilon>0$ is fixed and sufficiently small. The lower bound ensures that the interval contains roughly the expected number of square-free integers.

A related problem fixes a realization of $f$ and varies the starting point of the interval. Pandey, Wang and Xu \cite[Theorem~1.2]{PWX} proved that, for a set of Steinhaus realizations of measure $1-o(1)$,
\[
 \frac{1}{\sqrt H}\sum_{U_X<n\le U_X+H}f(n)
 \xrightarrow{d}\mathcal N_{\mathbb C}(0,1),
\]
where $U_X$ is uniformly chosen from $\{1,\ldots,X\}$ and the length $H\ll X/(\log X)^{W(X)}$ tends to infinity with $W(X)\to\infty$ arbitrarily slowly. In \cite[Section~4]{PWX}, they proposed such a problem: whether a different normalization can yield a Gaussian limit for fixed intervals when the square root normalization fails.

Caich \cite[Theorems~1.1 and~1.2]{Caich} showed that such a change of normalization is necessary when $y$ is sufficiently close to $x$. For both the Steinhaus and Rademacher models, he proved that
\[
 \frac{S_f(x,y)}{\sqrt y}\xrightarrow{\mathbb P}0,
\]
where $\log(x/y)=o(\sqrt{\log\log x})$. Thus a nondegenerate Gaussian limit in this range requires a different normalization.

Recently, Harper, Soundararajan and Xu \cite{HSX} surprisingly answered this question for Steinhaus functions. They proved that
\[
 \frac{S_f(x,y)}{\sqrt{V(x,y)}}
 \xrightarrow{d}\mathcal N_{\mathbb C}(0,1),
\]
where $y\to\infty$ and $y=o(x)$. The deterministic scale satisfies $V(x,y)\sim yp_x$, with
\[
 p_x=\frac1{\sqrt{2\pi}}\int_{-\kappa_x}^{\kappa_x}e^{-u^2/2}\dd u,\qquad
 \kappa_x=\frac{\log(x/y)}{\sqrt{2\log\log x}}.
\]

The results above give a Gaussian limit throughout the full short interval range for the Steinhaus model, after a suitable normalization. The corresponding result for the Rademacher model was left open. In \cite[Section~1.4]{HSX}, Harper, Soundararajan and Xu expected most of their argument to extend to this model, but noted that some nontrivial work was still needed.

Their proof uses translation invariance in law for Steinhaus model to reduce a varying frequency $t$ to $0$ and a pair $(t,s)$ to $(0,s-t)$ when they deal with the one and two point weighted expectations in the concentration of the conditional variance. These reductions fail for Rademacher signs. We combine Berry-Esseen estimates with suitable prime blocks and increment grids. This allows us to compare rounded barrier probabilities under the one point and two point tilted measures with their Gaussian counterparts. The relative errors are summable along every allowed path, even when the barrier probabilities tend to zero. For nearby frequencies, we condition on the correlated initial blocks and restart the comparison at a later scale. The second moment estimate retains the product of the survival probabilities of the two Gaussian random walks over the remaining blocks.

The conditional fourth moment requires a separate arithmetic argument. After separating a greatest common factor, we use a six edge parametrization of quadruples whose product is a square. Two complementary bounds for the resulting harmonic sum reduce the problem to a near diagonal two point integral. We estimate this integral using endpoint barrier bounds that account for the changing tilted mean and variance, uniformly down to zero frequency. These estimates lead to the following result.

\begin{theorem}\label{thm:main}
Let $f$ be a Rademacher random multiplicative function. Suppose that $x\to\infty$ and $\frac{x}{\sqrt{\log x}}\le y=o(x)$. With $p_x$ defined above, put $\vR=\frac6{\pi^2}yp_x$. Then
\begin{equation}\label{eq:targetclt}
 \frac{S_f(x,y)}{\sqrt{\vR}}\xrightarrow{d}N(0,1).
\end{equation}
\end{theorem}
The range in Theorem~\ref{thm:main} overlaps the Rademacher range of Soundararajan and Xu \cite[Corollary~9.2]{SX}. Combining the two results, we obtain the following Gaussian limit for a wider range of short intervals.
\begin{corollary}\label{cor:fullpicture}
Suppose $x^{1/5}\log x\ll y=o(x)$ and define
\begin{equation}\label{eq:unifiedvariance}
 \mathcal V_R(x,y)=
 \begin{cases}
 (6/\pi^2)y,&y\le x/\sqrt{\log x},\\
 v_R(x,y),&x/\sqrt{\log x}<y=o(x).
 \end{cases}
\end{equation}
Then $S_f(x,y)/\sqrt{\mathcal V_R(x,y)}\xrightarrow{d}N(0,1)$.
\end{corollary}
\begin{proof}
Use \cite[Corollary~9.2]{SX} below the joining scale and \cref{thm:main} above it. At $y=x/\sqrt{\log x}$ one has $\kappa_x=\sqrt{\log\log x}/(2\sqrt2)\to\infty$.
\end{proof}

\begin{remark}\label{rem:transition}
The normalization in Theorem~\ref{thm:main} agrees with the classical
one when $\log(x/y)/\sqrt{\log\log x}\to\infty$, since
$v_R(x,y)\sim(6/\pi^2)y$. Thus division by $\sqrt{(6/\pi^2)y}$
gives a standard Gaussian limit in this range.

When $y$ is sufficiently close to $x$ that
$\log(x/y)=o(\sqrt{\log\log x})$, the normalization instead satisfies
\begin{equation}\label{eq:smallkappa}
 \sqrt{v_R(x,y)}
 \sim \sqrt{\frac6{\pi^2}y}\,
 \frac{\sqrt{\log(x/y)}}{(\pi\log\log x)^{1/4}}.
\end{equation}
The additional factor tends to zero, so a smaller normalization is needed to obtain a standard Gaussian limit. This means the Gaussian fluctuations occur on a smaller scale, while the variance still satisfies $\operatorname{Var}(S_f(x,y))=\mathbb E|S_f(x,y)|^2\sim(6/\pi^2)y$ by orthogonality.
\end{remark}

As an application of our method, we also obtain a Gaussian limit for a fixed arithmetic twist. Let $\chi$ be a fixed primitive nonreal Dirichlet character modulo $q$, and define
\[
 S_{f,\chi}(x,y)=\sum_{x<n\le x+y}f(n)\chi(n).
\]
Put $c_q=\frac6{\pi^2}\prod_{p\mid q}\frac p{p+1}$ and $v_\chi(x,y)=c_qyp_x$. Let $\mathcal N_{\mathbb C}(0,1)$ denote the circular complex Gaussian random variable with independent real and imaginary parts distributed as real Gaussian random variable $N(0,1/2)$. 
\begin{theorem}\label{thm:twisted}
Let $f$ be a Rademacher random multiplicative function supported on square-free integers. Fix a positive integer $q$ and a primitive nonreal Dirichlet character $\chi$ modulo $q$. Suppose that $x\to\infty$ and $x^{1/5}\log x\ll y=o(x)$. Then
\begin{equation}\label{eq:twisted-main-intro}
 \frac{S_{f,\chi}(x,y)}{\sqrt{v_\chi(x,y)}}
 \xrightarrow{d}\mathcal N_{\mathbb C}(0,1).
\end{equation}
\end{theorem}
\begin{remark}
In the classical short interval range $x^{1/5}\log x\ll y\le x/(\log x)^{2\log 2-1+\varepsilon}$, the corresponding Gaussian limit with normalization $\sqrt{c_qy}$ can be deduced from the arguments of Soundararajan and Xu \cite[Section~9]{SX}. Theorem~\ref{thm:twisted} extends this conclusion throughout $x^{1/5}\log x\ll y=o(x)$, with normalization $\sqrt{c_qyp_x}$. In particular, it covers the regime $\log(x/y)=o(\sqrt{\log\log x})$, where this normalization is smaller than the standard deviation by a factor tending to zero.
\end{remark}
\subsection{Outline of Proof strategy}
We use the martingale central limit theorem framework of Harper, Soundararajan and Xu \cite{HSX}. Put $z=x^{1/\log\log\log x}$ and condition on $\cF_z=\sigma(\varepsilon_p:p\le z)$. Separating the smooth part and grouping the other terms by their largest prime factor, we write
\begin{equation}\label{eq:intro-martingale-chain}
 S_f(x,y)=S_{\le z}+\sum_{p>z}\Delta_p,
\end{equation}
where $\Delta_p=\varepsilon_pA_p$ and
\[
A_p:=\sum_{\substack{x/p<n\le(x+y)/p\\P(n)<p}}f(n).
\]
Here $P(n)$ is the largest prime factor of $n$, and $S_{\le z}$ is the part with $P(n)\le z$. It can be verified that the primes in increasing order makes $(\Delta_p)_p$ a martingale difference sequence. Put $\cQ_z=\sum_{p>z}A_p^2$ and $\Vf=\EE[\cQ_z\mid\cF_z]$. Apart from the smooth part and increment moment bounds proved in Section~\ref{sec:elementary}, the martingale theorem requires two estimates. The first is concentration of the conditional variance:
\begin{equation}\label{eq:intro-two-targets}
 \frac{\Vf}{\vR}\xrightarrow{\PP}1,
\end{equation}
which means the conditional variance is asymptotic in probability to the deterministic scale. The second controls the conditional fluctuation of $\cQ_z$:
\begin{equation}\label{eq:intro-fourth-target}
 \frac{\EE[(\cQ_z-\Vf)^2\mid\cF_z]}{\vR^2}
 \xrightarrow{\PP}0.
\end{equation}
Thus the fluctuations of $\mathcal Q_z$ around its conditional mean $V_f(x,y)$ are negligible on the scale $v_R(x,y)$.

For the first estimate, we adapt Iwaniec's linear sieve and Mellin-Parseval argument of Harper, Soundararajan and Xu to obtain
\[
 \Vf
 =\frac{xe^{1/T}}{2\pi T^2}
   \prod_{p\le z}(1-p^{-1})
   \int_{\mathbb R}|F_z(\tfrac12+it)|^2K_T(t)^2\,dt
   +o_{\mathbb P}(\vR),
\]
where $F_z(s)=\prod_{p\le z}(1+\varepsilon_pp^{-s})$ and $T=2/\log(1+y/x)$. The kernel and error are specified in Section~\ref{sec:CVtransfer}. We must therefore prove concentration of this random integral. After estimating the central frequencies and the tails separately, the identity $|F_z(\tfrac12-it)|=|F_z(\tfrac12+it)|$ allows us to restrict to positive frequencies. In particular, the two possible correlations $t\approx s$ and $t\approx-s$ can be reduced to $|t|\approx|s|$.

After restricting the Euler product integral to the relevant frequency range, a natural first attempt is to use Chebyshev's inequality and the second moment method to prove concentration around its mean. However, rare large values of the Euler product make the untruncated moments too large to describe its typical size. Therefore, Harper, Soundararajan and Xu introduce weak and strong barriers events. The weak barrier event preserves the integral with high probability, while the strong barrier event controls nearby frequencies in the second moment. The loss between the two weighted expectations is negligible, which allows concentration around the strongly restricted mean to transfer to the original integral. They compute the weighted expectations by splitting the logarithm into independent prime block increments. An approximate Girsanov change of measure and Gaussian comparison reduce the calculation to Gaussian barrier probabilities. Brownian motion estimates then give the factor $p_x$ in the restricted mean and hence in the normalization.

Thus, for the strongly restricted integral
\[
 J_{\mathcal H}
 =\int |F_z(\tfrac12+it)|^2K_T(t)^2
       \ind_{\mathcal H(t)}\,dt,
\]
over the truncated frequency range, the main task is to prove
\[
 \EE[J_{\mathcal H}^2]
 =(1+o(1))(\EE J_{\mathcal H})^2.
\]
To estimate the first and second moments here, we need to estimate one point and two point weighted expectations respectively.

In both calculations, Harper, Soundararajan and Xu use translation invariance in law: a single frequency reduces to $0$, and a frequency pair $(t,s)$ reduces to $(0,s-t)$. However, Rademacher signs do not have this invariance. Our key idea here is to compare the logarithmic increments with Gaussian variables under the one point and two point tilted measures directly. For the rounded logarithmic walk $\mathbf S$ and its Gaussian reference walk $\mathbf G$, an absolute error $\varepsilon$ is insufficient if the Gaussian barrier probability is smaller than $\varepsilon$. What we need, and prove for the barriers used here, is
\[
 \mathbb P_{\mathrm{tilt}}(\mathbf S\in\mathcal B)
 =(1+o(1))\mathbb P(\mathbf G\in\mathcal B).
\]
We obtain this relative estimate from Berry-Esseen bounds by choosing the prime blocks and increment grids so that the relative errors on admissible grid cells are summable. 

Therefore, the first moment estimate follows from the one point relative Gaussian comparison. For the second moment, we distinguish between separated and nearby frequencies. For separated positive frequencies, we compare the joint barrier
probability for the large prime part with that of two independent Gaussian reference walks, with relative error tending to zero. Nearby frequencies require a further change. Write $d=|t-s|$ and consider the Euler product over primes $v<p\le z$.
Its two point moment is controlled once $d\log v\gg1$. The relative path comparison, however, requires the accumulated mean and covariance errors to tend to zero. We condition on the initial blocks and restart that comparison only above a later cutoff satisfying $\log v\gg d^{-2}$. This keeps both Gaussian barrier probabilities over the remaining blocks, which are needed to compare the second moment with the square of the mean. The remaining small prime integral is controlled by an orthogonal
expansion in the prime signs.

For the conditional fourth moment, Rademacher orthogonality gives
\[
 \EE\prod_{j=1}^4f(r_j)=\ind_{r_1r_2r_3r_4=\square}
\]
for square-free $r_j$, where $\square$ denotes a perfect square. The additional quadruples in this relation require a separate arithmetic argument. After separating their greatest common factor, we parametrize them by six edge variables. Following the approach of Harper, Soundararajan and Xu, smoothing and Mellin inversion then lead to two complementary bounds for an Euler product kernel. After folding to positive frequencies, one bound controls separation within either of two fixed pairs, and the other controls separation between the pairs. The remaining four frequencies are all close together. Integrating out two leaves a narrow two point integral.

This last integral needs barrier estimates that remain valid near zero. For $\ell_p(t)=\log|1+\varepsilon_pp^{-1/2-it}|$, the one point tilt gives $\EE_t\ell_p(t)=p^{-1}+\cos(2t\log p)/(2p)+O(p^{-2})$ and $\Var_t\ell_p(t)=(1+\cos(2t\log p))/(2p)+O(p^{-2})$. The extra cosine terms are absent in the Steinhaus case and cannot be discarded near zero. We choose a different lower cutoff for the prime partition and adjust the weak and endpoint barriers. At each frequency, we compare the tilted walk with Brownian motion whose drift and time scale follow these moments. The resulting endpoint strip estimate is uniform down to $t=0$. It controls the weighted contribution where the weak barrier holds but the endpoint constraint fails; on the remaining part, that constraint bounds the Euler product directly. These bounds combined together prove \eqref{eq:intro-fourth-target} on the deterministic scale $v_R^2$.

For the fixed character twist, cancellation from the nonprincipal character $\chi^2$ controls the extra tilted terms and the additional complex quadratic term. Applying the real martingale central limit theorem to every real projection gives the circular complex Gaussian limit in the near full range. For shorter intervals, we use a weighted real central limit theorem
for the same projections.
\subsection{Organization and notation}
Section~\ref{sec:concentration} gives the martingale reduction and proves
the conditional variance estimate \eqref{eq:intro-two-targets}.
Section~\ref{sec:conditionalCLT} proves the conditional fourth moment bound
\eqref{eq:intro-fourth-target} and completes the proof of Theorem~\ref{thm:main}.
Section~\ref{sec:twists} treats the fixed character application.

We use $P(n)$ and $P^-(n)$ for the largest and smallest prime factors, with
$P(1)=1$ and $P^-(1)=\infty$. Implied constants may depend on fixed barrier
constants and, in Section~\ref{sec:twists}, on the fixed character.

\section{Concentration of the conditional variance}\label{sec:conditional}\label{sec:concentration}\label{sec:assembly}
Following Harper, Soundararajan and Xu \cite{HSX}, we begin by separating the smooth part of the short interval sum. Put $z=x^{1/\log\log\log x}$ and write
\[
 S_f(x,y)=S_{\le z}+S_{>z},
\]
where the smooth part is
\[
 S_{\le z}=\sum_{\substack{x<n\le x+y\\P(n)\le z}}f(n).
\]
The remaining terms have a prime factor exceeding $z$. Grouping them by their largest prime factor gives
\[
 S_{>z}=\sum_{z<p\le x+y}\Delta_p,
\]
where $\Delta_p=\eps_pA_p$ and
\[
 A_p=\sum_{\substack{x/p<n\le(x+y)/p\\P(n)<p}}f(n).
\]
It can be verified that $(\Delta_p)_p$ is a real martingale difference sequence as $p$ varies. We condition on $\cF_z=\sigma(\eps_p:p\le z)$ and consider its conditional variance
\[
 \Vf=\EE[S_{>z}^2\mid\cF_z].
\]
Our main aim is to compare $\Vf$ with $v_R(x,y)=(6/\pi^2)yp_x$. For $X_p=\Delta_p/\sqrt{\vR}$, the real martingale estimate \cite[Proposition~2.2]{SX} reduces the proof to four estimates. We first remove the smooth part in $L^2$:
\begin{equation}\label{eq:smooth-target}
 \frac{\EE S_{\le z}^2}{\vR}\to0.
\end{equation}
The maximum increment error $\EE\max_p|X_p|$ is at most $(\sum_{p>z}\EE\Delta_p^4/\vR^2)^{1/4}$ by H\"older's inequality, so we need to prove
\begin{equation}\label{eq:increments-target}
 \frac{\sum_{p>z}\EE\Delta_p^4}{\vR^2}\to0.
\end{equation}
For the quadratic variation error $\EE\min\{1,|\sum_pX_p^2-1|\}$, martingale orthogonality gives conditional mean $\Vf/\vR$. We first prove
\begin{equation}\label{eq:variance-target}
 \frac{\Vf}{\vR}\xrightarrow{\PP}1.
\end{equation}
The conditional Cauchy-Schwarz inequality then reduces the fluctuation about this mean to
\begin{equation}\label{eq:fluctuation-target}
 \frac{\EE\left[\left(\sum_{p>z}\Delta_p^2-\Vf\right)^2\mid\cF_z\right]}{\vR^2}
 \xrightarrow{\PP}0.
\end{equation}
Proposition~\ref{thm:assembly} records this deduction in quantitative form.

Section~\ref{sec:elementary} proves the first two estimates. For \eqref{eq:variance-target}, Section~\ref{sec:variance_tools} gives the Euler product representation; Sections~\ref{sec:barriers}-\ref{sec:barrierintegration} introduce the barriers, compare their means and prove concentration. Section~\ref{sec:conditionalCLT} proves \eqref{eq:fluctuation-target}.

\subsection{Elementary estimates}\label{sec:elementary}
We first prove \eqref{eq:smooth-target}. Write $M=\log\log z$, so that $M\sim\log\log x$.

\begin{lemma}\label{lem:smooth}
For every fixed $B>0$, we have
\[
 \EE S_{\le z}^2\ll_B yM^{-B}.
\]
\end{lemma}
\begin{proof}
Expanding the square and using orthogonality gives
\[
 \begin{aligned}
 \EE S_{\le z}^2
 &=\sum_{\substack{x<m,n\le x+y\\P(m),P(n)\le z}}
      \EE[f(m)f(n)]\\
 &=\sum_{\substack{x<n\le x+y\\P(n)\le z}}\mu^2(n)\\
 &\le\#\{x<n\le x+y:P(n)\le z\}.
 \end{aligned}
\]
For $z=x^{1/\log\log\log x}$ and $x/\sqrt{\log x}\le y=o(x)$, Hildebrand's short interval estimate \cite{HildebrandSmooth}, combined with the smooth number bound in \cite[(1.12)]{GranvilleSmooth}, gives the following bound; see also \cite[proof of Proposition~2.1]{HSX}:
\[
 \#\{x<n\le x+y:P(n)\le z\}
 \ll_B y(\log\log x)^{-B}
\]
for every fixed $B>0$. Since $M\sim\log\log x$, this proves the lemma.
\end{proof}

The definition of $p_x$ implies $v_R(x,y)\gg y/\sqrt M$. Therefore taking a fixed $B>1/2$ in Lemma~\ref{lem:smooth}, we obtain $\EE S_{\le z}^2/v_R(x,y)\ll_B M^{1/2-B}=o(1)$, which proves \eqref{eq:smooth-target}.

Next we focus on proving \eqref{eq:increments-target}. Since $\Delta_p=\eps_pA_p$, its fourth moment is the fourth moment of a restricted Rademacher sum. We first establish a bound that applies to each $A_p$.

\begin{lemma}\label{lem:fourth}
Let $W\ge1$ and $\mathcal I\subseteq[1,W]\cap\NN$. For fixed complex coefficients satisfying $|c_n|\le1$, we have
\[
 \EE\left|\sum_{n\in\mathcal I}c_nf(n)\right|^4
 \ll W|\mathcal I|(\log(2W))^8.
\]
\end{lemma}
\begin{proof}
We give the one sign argument and then iterate it. For complex numbers $a,b$,
\[
 \begin{aligned}
 \frac{|a+b|^4+|a-b|^4}{2}
 &=(|a|^2+|b|^2)^2+4\{\Re(a\overline b)\}^2\\
 &\le(|a|^2+3|b|^2)^2.
 \end{aligned}
\]
Let $\|\cdot\|_r$ denote the $L^r$ norm over the prime signs. A polynomial in finitely many signs can be written as $F=F_0+\eps_pF_1$, where $F_0,F_1$ depend only on the other signs. Averaging first over $\eps_p$ and then using the triangle inequality in $L^2$, we get
\[
 \begin{aligned}
 \|F\|_4^2
 &=\left(\EE\left[\frac{|F_0+F_1|^4+|F_0-F_1|^4}{2}\right]\right)^{1/2}\\
 &\le\left(\EE(|F_0|^2+3|F_1|^2)^2\right)^{1/2}\\
 &\le\|F_0\|_4^2+3\|F_1\|_4^2.
 \end{aligned}
\]
Apply this inequality successively to the prime signs occurring in $\sum_{n\in\mathcal I}c_nf(n)$. Each prime dividing a square-free index contributes a factor of $3$, while terms with a square factor vanish. Thus
\[
 \left\|\sum_{n\in\mathcal I}c_nf(n)\right\|_4^2
 \le\sum_{n\in\mathcal I}3^{\omega(n)}\mu^2(n)|c_n|^2.
\]
Squaring and applying the Cauchy-Schwarz inequality, we obtain
\[
 \begin{aligned}
 \EE\left|\sum_{n\in\mathcal I}c_nf(n)\right|^4
 &\le\left(\sum_{n\in\mathcal I}3^{\omega(n)}\mu^2(n)|c_n|^2\right)^2\\
 &\le\left(\sum_{n\in\mathcal I}\mu^2(n)|c_n|^2\right)
       \left(\sum_{n\in\mathcal I}9^{\omega(n)}\mu^2(n)|c_n|^2\right)\\
 &\le|\mathcal I|\sum_{n\le W}d_9(n)\\
 &\le W|\mathcal I|\left(\sum_{m\le W}\frac1m\right)^8\\
 &\ll W|\mathcal I|(\log(2W))^8.
 \end{aligned}
\]
Here $d_9$ is the 9-fold divisor function. We used $9^{\omega(n)}=d_9(n)$ for square-free $n$ and $|c_n|\le1$. The summatory bound follows by writing $d_9(n)$ as a count of factorizations $n=n_1\cdots n_9$ and summing over the last factor first.
\end{proof}

We apply Lemma~\ref{lem:fourth} to the sums $A_p$. Write $\mathcal L_x=\sum_{p>z}\EE\Delta_p^4/v_R(x,y)^2$.

\begin{corollary}\label{cor:increments}
In the range of Theorem~\ref{thm:main}, we have
\[
 \sum_{p>z}\EE\Delta_p^4
 \ll\left(\frac{xy}{z}+x\log(2x)\right)(\log(2x))^8.
\]
Consequently,
\[
 \frac{\sum_{p>z}\EE\Delta_p^4}{v_R(x,y)^2}\to0.
\]
\end{corollary}
\begin{proof}
For $z<p\le x+y$, put $W_p=(x+y)/p$ and
\[
 \mathcal I_p=\{n\in\NN:x/p<n\le(x+y)/p,\ P(n)<p\}.
\]
Then $A_p=\sum_{n\in\mathcal I_p}f(n)$, $W_p\ge1$, $W_p\ll x/p$ and $|\mathcal I_p|\le y/p+1$. Since $|\Delta_p|=|A_p|$, Lemma~\ref{lem:fourth} yields
\[
 \begin{aligned}
 \sum_{p>z}\EE\Delta_p^4
 &\ll\sum_{z<p\le x+y}
        \frac xp\left(\frac yp+1\right)(\log(2x))^8\\
 &\ll(\log(2x))^8
        \left(xy\sum_{p>z}\frac1{p^2}
             +x\sum_{p\le2x}\frac1p\right)\\
 &\ll\left(\frac{xy}{z}+x\log(2x)\right)(\log(2x))^8.
 \end{aligned}
\]
In the last step we used $\sum_{p>z}p^{-2}\ll z^{-1}$ and $\sum_{p\le2x}p^{-1}\ll\log(2x)$.
Since $v_R(x,y)^2\gg y^2/M$ and $y\ge x/\sqrt{\log x}$, dividing by $v_R(x,y)^2$ gives
\[
 \begin{aligned}
 \mathcal L_x
 &\ll\frac{M}{y^2}
       \left(\frac{xy}{z}+x\log(2x)\right)(\log(2x))^8\\
 &=M\left(\frac{x(\log(2x))^8}{yz}
          +\frac{x(\log(2x))^9}{y^2}\right)\\
 &\ll M\left(\frac{(\log(2x))^{17/2}}z
             +\frac{(\log(2x))^{10}}x\right)=o(1).
 \end{aligned}
\]
The last expression tends to zero for our choice of $z$.
\end{proof}

This proves \eqref{eq:increments-target} and completes the first two estimates.

\subsection{Euler product representation and tilted moments}\label{sec:variance_tools}
We now turn to \eqref{eq:variance-target}, using an Euler product representation and estimates for the tilted moments. Write $\Lambda=\log z$.

\subsubsection{Mellin-Parseval representation}\label{sec:CVtransfer}
Following \cite[Proposition~2.1]{HSX}, we first average over the rough factor by a sieve and then apply Mellin-Parseval to the smooth sum. In the Rademacher model, we must also account for the square-free restriction on the rough factor.

Let $\mathcal A_z$ be the square-free integers $r>1$ all of whose prime divisors exceed $z$. Factoring each integer into its smooth and rough parts gives
\begin{equation}\label{eq:roughdecomp}
 S_{>z}=\sum_{r\in\mathcal A_z}f(r)a(r),
\end{equation}
where, for real $r>0$, we put
\[
 a(r)=\sum_{\substack{x/r<n\le(x+y)/r\\P(n)\le z}}f(n).
\]
The coefficients $a(r)$ are $\cF_z$-measurable, and the rough factors are orthogonal after conditioning on $\cF_z$. Consequently
\[
 \Vf=\sum_{r\in\mathcal A_z}a(r)^2.
\]
Only $r\le x+y$ contribute. To apply the sieve, it is convenient to sum over all rough integers. Define
\[
 V_f^{\mathrm{all}}(x,y)
 =\sum_{\substack{r>1\\P^-(r)>z}}a(r)^2.
\]
The added rough factors have a factor $p^2$ with $p>z$. We first bound their mean contribution.

\begin{lemma}\label{lem:roughsquare}
We have
\[
 0\le\EE\big(V_f^{\mathrm{all}}(x,y)-\Vf\big)\ll \frac{y}{z}+\sqrt x.
\]
\end{lemma}
\begin{proof}
Expanding $a(r)^2$ and using orthogonality counts pairs $(r,n)$ with weight $\mu^2(n)\le1$. The factorization $N=rn$ into smooth and rough parts is unique, and a rough factor that is not square-free forces $p^2\mid N$ for some $p>z$. Hence
\[
 \begin{aligned}
 \EE\big(V_f^{\mathrm{all}}(x,y)-\Vf\big)
 &\le\sum_{z<p\le\sqrt{2x}}
       \#\{x<N\le x+y:p^2\mid N\}\\
 &\le\sum_{z<p\le\sqrt{2x}}(y/p^2+1)
 \ll y/z+\sqrt x.\qedhere
 \end{aligned}
\]
\end{proof}

Applying Markov's inequality to the nonnegative difference, we obtain
\[
 \PP\big(V_f^{\mathrm{all}}(x,y)-\Vf>yM^{-55}\big)
 \ll M^{55}\big(z^{-1}+\sqrt x/y\big)=o(M^{-10}).
\]
Here we used $y\ge x/\sqrt{\log x}$ and the choice of $z$. Thus the added rough factors contribute $O(yM^{-55})$ outside this exceptional event. We restore the square-free restriction at the end.

The smooth sums are encoded by the finite Euler product
\[
 F_z(s)=\sum_{\substack{n\ge1\\P(n)\le z}}\frac{f(n)}{n^s}
       =\prod_{p\le z}\left(1+\frac{f(p)}{p^s}\right).
\]
Write $Z_z=\prod_{p\le z}(1+p^{-1})$ and $d_z=\prod_{p\le z}(1-p^{-1})$. Orthogonality gives, for every real $t$,
\begin{equation}\label{eq:Fzsecond}
 \EE\left|F_z(\tfrac12+it)\right|^2=Z_z.
\end{equation}
Put $T=2/\log(1+y/x)$, so that $T\sim2x/y\to\infty$, and define
\[
 K_T(t)=T\left|
       \frac{e^{(1/2+it)/T}-e^{-(1/2+it)/T}}{1/2+it}
       \right|.
\]
The following representation retains $d_z$ to identify the exact normalization.

\begin{proposition}\label{prop:transfer}
The following asymptotic formula holds with probability $1-O(M^{-10})$,
\begin{equation}\label{eq:CVtransfer}
 \Vf
 =d_z\frac{xe^{1/T}}{2\pi T^2}
 \int_{\mathbb R}\left|F_z(\tfrac12+it)\right|^2K_T(t)^2\,dt
 +O\left(\frac{y}{M^{50}}\right).
\end{equation}
\end{proposition}
\begin{proof}
We use the sieve argument in \cite[proof of Proposition~2.1, (2.4)-(2.7)]{HSX}, keeping the density $d_z$ exact. Put $K=\lfloor x^{1/10}\rfloor$ and $R=(x+y)/(K+1)$. For $1<r\le R$, orthogonality and the smooth number estimate used there give $\EE a(r)^2\ll_D y/(rM^D)$ for every fixed $D>0$. Thus the expected contribution from this range is $\ll_D yM^{-D}(\log x)/\Lambda$.

For the remaining rough factors, partition $(R,x+y]$ into $I_k=((x+y)/(k+1),(x+y)/k]$, $1\le k\le K$, and put $B_k=a((x+y)/k)$. On each $I_k$, the upper integer endpoint of $a(r)$ stays fixed and the lower endpoint moves by less than one, so $|a(r)-B_k|\le1$. Since $|I_k|\gg x^{4/5}$, the linear sieve estimate of Iwaniec \cite{Iwaniec}, in the form recorded in \cite[Lemma~2.2]{HSX}, gives
\[
 \#\{r\in I_k:P^-(r)>z\}
 =d_z|I_k|\big(1+O_D(M^{-D})\big).
\]
The sieve parameter is comparable to $\log\log\log x$, so $D$ can be any fixed constant. Using $|a(r)^2-B_k^2|\le2|B_k|+1$, $\EE B_k^2\le yk/x+2$ and the Cauchy-Schwarz inequality, the endpoint and sieve errors satisfy
\[
 \begin{aligned}
 \EE\left|V_f^{\mathrm{all}}(x,y)-d_z\int_R^{x+y}a(w)^2\,dw\right|
 &\ll_D\frac{y\log x+x}{\Lambda M^D}
       +\frac{x+\sqrt{xy}}{\Lambda}\\
 &\ll yM^{-80}.
 \end{aligned}
\]
Here the first term also includes the initial rough range, and the last bound uses $y\ge x/\sqrt{\log x}$ and sufficiently large $D$. We apply Markov's inequality to the absolute difference on the left. It is $O(yM^{-55})$ outside an event of probability $O(M^{-25})$. The argument uses only orthogonality and the upper bound by the number of smooth integers, so it applies to the present coefficients.

Write $\delta_0=y/x$ and put
\[
 B(u)=\sum_{\substack{u<n\le(1+\delta_0)u\\P(n)\le z}}f(n).
\]
The substitution $u=x/w$ gives an integral up to $u_1=x/R\asymp x^{1/10}$. As in the passage from (2.6) to (2.7) of \cite{HSX}, orthogonality and Rankin's inequality give, for every fixed $A>0$,
\[
 \begin{aligned}
 d_zx\,\EE\int_{u_1}^{\infty}|B(u)|^2\frac{du}{u^2}
 &\le d_zyu_1^{-A/\Lambda}
       \prod_{p\le z}(1-p^{-1+A/\Lambda})^{-1}\\
 &\ll_A yM^{-A/11}.
 \end{aligned}
\]
We choose $A$ large enough that the last expectation is $O(yM^{-80})$. By Markov's inequality, the omitted tail is $O(yM^{-55})$ outside an event of probability $O(M^{-25})$. Combining this with the sieve approximation, we obtain
\[
 V_f^{\mathrm{all}}(x,y)
 =d_zx\int_{1/(1+\delta_0)}^{\infty}|B(u)|^2\frac{du}{u^2}
  +O(yM^{-55}),
\]
with probability $1-O(M^{-25})$.

To evaluate this integral, set $u=e^{v-1/T}$ and $H(v)=e^{-v/2}B(e^{v-1/T})$. Since $1+\delta_0=e^{2/T}$, termwise integration gives the Fourier transform
\[
 \widehat H(t):=\int_{\mathbb R}H(v)e^{-itv}\,dv
 =F_z(\tfrac12+it)
       \frac{e^{(1/2+it)/T}-e^{-(1/2+it)/T}}{1/2+it}.
\]
The change of variables and Parseval's identity now give
\[
 \begin{aligned}
 x\int_{1/(1+\delta_0)}^{\infty}|B(u)|^2\frac{du}{u^2}
 &=xe^{1/T}\int_{\mathbb R}|H(v)|^2\,dv\\
 &=\frac{xe^{1/T}}{2\pi T^2}
       \int_{\mathbb R}|F_z(\tfrac12+it)|^2K_T(t)^2\,dt.
 \end{aligned}
\]
Substituting into the approximation for $V_f^{\mathrm{all}}$ and restoring the square-free restriction proves \eqref{eq:CVtransfer}, with error $O(yM^{-55})$ and total exceptional probability $O(M^{-10})$.
\end{proof}

Since $v_R(x,y)\gg y/\sqrt M$, the error in \eqref{eq:CVtransfer} is negligible on the required scale. To use the representation, we record a pointwise bound, an exact integral identity and a variation estimate for its kernel. We write $\operatorname{TV}(g)$ for the total variation of $g$.

\begin{lemma}\label{lem:kernel}
For real $t$, we have
\[
 K_T(t)^2\ll\min\{1,T^2/(1+t^2)\}.
\]
The kernel has the exact normalization
\[
 \frac{xe^{1/T}}{2\pi T^2}
       \int_{\mathbb R}K_T(t)^2\,dt=y.
\]
Moreover, if $U=T(\log T)^{100}$ and $g(t)=K_T(t)^2\mathbf1_{[1,U]}(t)$, then $\int|g|\asymp T$ and $\|g\|_\infty+\operatorname{TV}(g)\ll1$.
\end{lemma}
\begin{proof}
The pointwise estimate follows directly from the formula for $K_T$. The derivative calculation in \cite[proof of Lemma~3.2, (3.2)]{HSX}, for the same kernel, gives $\int_{\mathbb R}|(K_T^2)'(t)|\,dt\ll1$. Hence
\[
 \operatorname{TV}(g)
 \le K_T(1)^2+K_T(U)^2+\int_1^U|(K_T^2)'(t)|\,dt
 \ll1,
\]
which also proves $\|g\|_\infty+\operatorname{TV}(g)\ll1$.

For the exact integral, apply Parseval's identity to $q(v)=e^{v/2}\mathbf1_{[-1/T,1/T]}(v)$, whose Fourier transform has modulus $K_T(t)/T$:
\[
 \int_{\mathbb R}K_T(t)^2\,dt
 =2\pi T^2\int_{-1/T}^{1/T}e^v\,dv
 =4\pi T^2\sinh(1/T).
\]
Since $y=2xe^{1/T}\sinh(1/T)$, this is the stated normalization. Finally, the pointwise bound shows that removing $|t|>U$ loses $O(T^2/U)=o(T)$, while removing $[-1,1]$ loses $O(1)$. By evenness,
\[
 \int_1^U K_T(t)^2\,dt\sim2\pi T,
\]
which also proves $\int|g|\asymp T$.
\end{proof}

The mean factor from the Euler product combines with the rough density as
\[
 d_zZ_z=\prod_{p\le z}(1-p^{-2})\to\frac6{\pi^2}.
\]
Together with the exact kernel normalization, this identifies the square-free density in $v_R(x,y)$. The kernel gives a weighted average over a frequency range of length $T\asymp x/y$. We will use barriers to control large values and estimate the dependence between separated frequencies to prove concentration.

Before introducing the barriers, we explain the restriction on the frequency range. With $U=T(\log T)^{100}$ as in Lemma~\ref{lem:kernel}, the reduction we need is
\begin{equation}\label{eq:frequency-reduction}
 \begin{aligned}
 \int_{\mathbb R}|F_z(\tfrac12+it)|^2K_T(t)^2\,dt
 &=2\int_1^U|F_z(\tfrac12+it)|^2K_T(t)^2\,dt\\
 &\quad+o_{\PP}(TZ_zp_x).
 \end{aligned}
\end{equation}
The central interval requires a low moment bound, since its ordinary mean need not be negligible on the scale $TZ_zp_x$.

\begin{lemma}\label{lem:central}
For some fixed $\nu\in(2/3,1)$,
\begin{equation}\label{eq:centralmoment}
 \EE\left(\int_{-1}^1|F_z(\tfrac12+it)|^2\,dt\right)^\nu
 \ll_\nu\left(\frac{Z_z}{\sqrt M}\right)^\nu.
\end{equation}
Consequently
\begin{equation}\label{eq:centralnegligible}
 \int_{-1}^1|F_z(\tfrac12+it)|^2K_T(t)^2\,dt
 =O_{\PP}(Z_z/\sqrt M)
 =o_{\PP}(TZ_zp_x).
\end{equation}
\end{lemma}
\begin{proof}
Take the large parameter in \cite[Section~4.4]{Harper} to be $X=z^e$ and take $k=0$, so that its Euler product $F_0$ has prime cutoff $z$. That section bounds the $\nu$th moment of the integral over $[N-1/2,N+1/2]$, for each fixed integer $N$, by $O_\nu((\log z/\sqrt{\log\log z})^\nu)$. Apply this with $N=-1,0,1$ and use $(a+b+c)^\nu\le a^\nu+b^\nu+c^\nu$. Since $Z_z\asymp\log z$ by Mertens's theorem (see \cite[Section~2.2]{MV}) and $M=\log\log z$, this proves \eqref{eq:centralmoment}.

Using the pointwise bound for $K_T$ and applying Markov's inequality to the $\nu$th moment in \eqref{eq:centralmoment}, we obtain the first estimate in \eqref{eq:centralnegligible}. For the second, $T\sim2x/y$ and $M\sim\log\log x$ imply $p_x\asymp\min\{1,\log T/\sqrt M\}$. Hence
\[
 \frac{Z_z/\sqrt M}{TZ_zp_x}
 \asymp\frac1{T\min\{\sqrt M,\log T\}}\to0.\qedhere
\]
\end{proof}

The remaining estimate is
\begin{equation}\label{eq:tail-negligible}
 \int_{|t|>U}|F_z(\tfrac12+it)|^2K_T(t)^2\,dt
 =o_{\PP}(TZ_zp_x).
\end{equation}
Its proof uses the weak barrier and the one point Gaussian comparison, so we defer it to Section~\ref{sec:CVcompletion}. Since $|F_z(\tfrac12-it)|=|F_z(\tfrac12+it)|$ and $K_T(t)^2$ is even, \eqref{eq:centralnegligible} and \eqref{eq:tail-negligible} together imply \eqref{eq:frequency-reduction}. The barrier argument will therefore concern only $[1,U]$. We first record the tilted moments needed for that argument.

\subsubsection{Prime block tilts}\label{sec:local}
We record the estimates needed for the change of measure of Girsanov type and the relative Gaussian comparison in Section~\ref{sec:relative}. The first and second moments of the Euler product integral lead to one point and two point tilted measures, respectively. We compute the means and covariances of the logarithmic increments under these measures, together with the third moment bounds needed for Gaussian approximation. Since these moments depend on the frequencies in the Rademacher
model, we begin with an estimate for the oscillatory prime sums that arise.

\begin{lemma}\label{lem:primeosc}
For $v$ sufficiently large and real $\xi\ne0$, we have
\[
 \sum_{v<p\le v^e}\frac{\cos(\xi\log p)}p
 \ll \frac1{|\xi|\log v}+(1+|\xi|)e^{-c\sqrt{\log v}},
\]
where $c>0$ is absolute.
\end{lemma}
\begin{proof}
Partial summation with the classical prime number theorem, followed by integration by parts, gives
\[
 \begin{aligned}
 \sum_{v<p\le v^e}\frac{\cos(\xi\log p)}p
 =\left[\frac{\sin(\xi u)}{\xi u}\right]_{\log v}^{e\log v}
   +\frac1\xi\int_{\log v}^{e\log v}\frac{\sin(\xi u)}{u^2}\,du+O\big((1+|\xi|)e^{-c\sqrt{\log v}}\big).
 \end{aligned}
\]
The boundary and integral terms are $O((|\xi|\log v)^{-1})$. The factor $1+|\xi|$ in the error term is needed when the frequency grows with $x$; compare \cite[Sections~3-4]{HSX}.
\end{proof}

Put $\ell_p(t)=\log|1+\eps_pp^{-1/2-it}|$. The normalized square of an Euler factor is
\begin{equation}\label{eq:localdensity}
 \frac{|1+\eps_pp^{-1/2-it}|^2}{1+1/p}=1+\eps_pc_p(t),
\end{equation}
where $c_p(t)=2\sqrt p\cos(t\log p)/(p+1)$. For a fixed prime $p$ and any function $H:\{-1,1\}\to\mathbb C$, define the one point tilted expectation by
\[
 \EE_t H(\eps_p)
 =\EE\big[H(\eps_p)(1+\eps_pc_p(t))\big].
\]
The weight has expectation one because $\EE\eps_p=0$. For two frequencies $t,s$, define
\[
 \EE_{t,s}H(\eps_p)
 =\frac{\EE\big[H(\eps_p)
          (1+\eps_pc_p(t))(1+\eps_pc_p(s))\big]}
        {1+c_p(t)c_p(s)}.
\]
Here the denominator is the expectation of the product of the two weights. Since $|c_p(u)|<1$ for every real $u$, these normalized weights are positive and define probability laws. For functions of finitely many prime signs, we use the product of the corresponding tilted laws at the primes. In particular, the prime signs remain independent under either tilt.

\begin{lemma}\label{lem:tilt}
Uniformly for real $t,s$, the one point tilt satisfies
\[
 \EE_t\ell_p(t)=\frac{1}{p}+\frac{\cos(2t\log p)}{2p}+O(p^{-2}),
\]
and
\[
 \Var_t\ell_p(t)=\frac{1}{2p}+\frac{\cos(2t\log p)}{2p}+O(p^{-2}).
\]
Under the two point tilt, we have
\[
 \EE_{t,s}\ell_p(t)
 =\frac{1}{p}+\frac{\cos(2t\log p)}{2p}
 +\frac{\cos((t-s)\log p)+\cos((t+s)\log p)}p+O(p^{-2}),
\]
and
\[
 \Cov_{t,s}(\ell_p(t),\ell_p(s))
 =\frac{\cos((t-s)\log p)+\cos((t+s)\log p)}{2p}+O(p^{-2}).
\]
\end{lemma}
\begin{proof}
Write $\theta=t\log p$, $\phi=s\log p$ and $\ell_p(t)=A_p(t)+\eps_pB_p(t)$. The even and odd Taylor terms give $A_p(t)=-\cos(2\theta)/(2p)+O(p^{-2})$ and $B_p(t)=\cos\theta/\sqrt p+O(p^{-3/2})$.

The one point tilted mean of $\eps_p$ is $c_p(t)$. For two points, the normalizing factor is $1+c_p(t)c_p(s)$, and hence
\[
 m_{t,s}:=\EE_{t,s}\eps_p
 =\frac{c_p(t)+c_p(s)}{1+c_p(t)c_p(s)}
 =\frac{2(\cos\theta+\cos\phi)}{\sqrt p}+O(p^{-3/2}).
\]
Substituting this mean into the logarithmic increment gives
\[
 \begin{aligned}
 \EE_{t,s}\ell_p(t)
 &=A_p(t)+B_p(t)m_{t,s}\\
 &=-\frac{\cos(2\theta)}{2p}
   +\frac{2\cos\theta(\cos\theta+\cos\phi)}p+O(p^{-2}).
 \end{aligned}
\]
The two point mean follows from the product formula for cosines, with $m_{t,s}$ replaced by $c_p(t)$ for the one point case. For the covariance we use $\Cov_{t,s}(\ell_p(t),\ell_p(s))=B_p(t)B_p(s)(1-m_{t,s}^2)$. Compare \cite[Section~3.3]{Harper}.
\end{proof}

The same calculation also shows that, uniformly for real $t,s$,
\[
 \Var_{t,s}\ell_p(u)
 =\frac{1+\cos(2u\log p)}{2p}+O(p^{-2}),
 \qquad u\in\{t,s\}.
\]
Since $|B_p(u)|\ll p^{-1/2}$ and the centered sign is bounded by $2$, we have
\[
 \EE_t\big|\ell_p(t)-\EE_t\ell_p(t)\big|^3\ll p^{-3/2}.
\]
For $L_p=(\ell_p(t),\ell_p(s))$, the corresponding vector bound is
\[
 \EE_{t,s}\big\|L_p-\EE_{t,s}L_p\big\|_2^3\ll p^{-3/2}.
\]
Prime contributions in disjoint blocks remain independent under these tilted measures. These third moment bounds allow us to apply the following form of the Berry-Esseen inequality when the block covariance is nondegenerate.

\begin{lemma}\label{lem:BE}
Let $X_i$ be independent centered vectors in $\RR^d$, where $d=1$ or $2$, and put $\Sigma=\sum_i\Cov(X_i)$. Suppose the eigenvalues of $\Sigma$ lie in a fixed compact subset of $(0,\infty)$. If $G$ is a centered Gaussian vector with covariance $\Sigma$, then for every convex Borel set $C$,
\[
 |\PP(\sum_iX_i\in C)-\PP(G\in C)|
 \ll\sum_i\EE\|X_i\|_2^3,
\]
with the constant depending only on $d$ and the fixed covariance bounds.
\end{lemma}
This follows from \cite[Theorem~1.1]{Raic} after standardizing the covariance. Rectangles and their nonsingular linear images are convex. We check the covariance condition for the prime blocks in the relative comparison in Section~\ref{sec:relative}.

\subsection{Weak and strong barriers}\label{sec:barriers}
We now study the integral over $[1,U]$ in \eqref{eq:frequency-reduction}. Rare large values prevent a direct second moment argument. Following \cite[Section~3]{HSX}, we introduce two barrier events. The weak barrier will hold throughout this interval with high probability. The strong barrier gives the additional saving needed for the second moment.

Let $m=\lceil M-(\log T)^{1/200}\rceil$ and put $b=M-m$. For $0\le j\le m$, set $w_j=\exp(e^{b+j})$ and write $w=w_0$. Thus $b=(\log T)^{1/200}+O(1)$, $\log\log w_j=b+j$ and $w_m=z$. The partial Euler products are
\[
 R_j(t)=\prod_{w<p\le w_j}(1+\eps_pp^{-1/2-it}).
\]
These products are cumulative and share prime factors. The logarithmic increments
\[
 \log|R_j(t)|-\log|R_{j-1}(t)|
 =\sum_{w_{j-1}<p\le w_j}
       \log|1+\eps_pp^{-1/2-it}|
\]
involve disjoint prime blocks and are independent. We choose the starting point $w$ to balance two requirements. The relation $\log w=e^b\gg(\log T)^C$ for every fixed $C>0$ makes the Gaussian approximation errors small enough for the relative comparison. On the other hand, $e^{5b}/T\to0$, which will make the small prime variance bound in Lemma~\ref{lem:smallvariance} negligible compared with its mean squared.

\paragraph{The weak barrier $\mathcal G_T(t)$.}
Fix a sufficiently large constant $B_0>2$. For each real $t$, let $\mathcal G_T(t)$ be the event that
\begin{equation}\label{eq:weakbarrier}
 (\sqrt T\,e^j)^{-B_0}\le |R_j(t)|
 \le \sqrt T\,e^j\exp\big((\log T)^{1/100}\big)\log\log w_j
 \qquad(1\le j\le m).
\end{equation}
The upper bound excludes large values of the partial Euler products. The lower bound is a technical restriction that will keep the logarithmic increments bounded in the relative comparison.

\paragraph{The strong barrier $\mathcal H_T(t)$.}
We retain the lower bound and impose a smaller upper bound. Let $\mathcal H_T(t)$ be the event that
\begin{equation}\label{eq:strongbarrier-product}
 (\sqrt T\,e^j)^{-B_0}\le |R_j(t)|
 \le \sqrt T\,e^j
       \frac{\exp\big(- (\log T)^{1/100}\big)}{(\log\log w_j)^5}
 \qquad(1\le j\le m).
\end{equation}
In particular, $\mathcal H_T(t)\subseteq\mathcal G_T(t)$. The factor $(\log\log w_j)^{-5}$ makes the later sum over separation scales convergent; see \eqref{eq:nearseries}. The exponential factor absorbs the contribution from primes at most $w$ and leaves a vanishing error for nearby frequencies. More precisely, after normalization that error contains $\exp(4b-2(\log T)^{1/100})=o(1)$. A rounded version of this event will be introduced in Section~\ref{sec:barrierprobabilities}, where the Gaussian comparison is applied.

We first prove that the weak barrier holds throughout $[1,U]$.

\begin{proposition}\label{lem:weakbarrier}
We have
\begin{equation}\label{eq:weakgood}
 \PP\left(\mathcal G_T(t)\text{ fails for some }t\in[1,U]\right)=o(1).
\end{equation}
\end{proposition}
\begin{proof}
Fix $j$ and write $D_j=\log w_j=e^{b+j}$. We need bounds throughout $[1,U]$, for both $R_j$ and its reciprocal. Cover this interval by $O(UD_j)$ intervals of length $h_j=D_j^{-1}$. For $Q=R_j$ or $R_j^{-1}$, we use the Sobolev inequality on each such interval $I$:
\[
 \sup_{t\in I}|Q(t)|^2
 \ll h_j^{-1}\int_I|Q(t)|^2\,dt
       +h_j\int_I|Q'(t)|^2\,dt.
\]
We therefore estimate the second moments of these functions and their derivatives. The length $h_j$ will balance the two terms above.

Orthogonality gives, uniformly in real $t$,
\[
 \EE|R_j(t)|^2\ll e^j,\qquad
 \EE|R_j'(t)|^2
 =\sum_{n\in\mathcal S_j}\frac{(\log n)^2}{n}
 \ll e^jD_j^2,
\]
where $\mathcal S_j$ is the set of square-free integers supported on $(w,w_j]$. For the derivative bound, expand $(\log n)^2$ and use
$\sum_{p\le w_j}(\log p)/p\ll D_j$ and $\sum_{p\le w_j}(\log p)^2/p\ll D_j^2$.

For the lower boundary we apply the same argument to $R_j^{-1}$. Every Euler factor is nonzero, and averaging one sign, with $\theta=t\log p$, gives
\[
 \EE|1+\eps_pp^{-1/2-it}|^{-2}
 =1+\frac{1+2\cos(2\theta)}p+O(p^{-2})
 \le1+\frac3p+O(p^{-2}).
\]
Multiplication over primes yields $\EE|R_j(t)|^{-2}\ll e^{3j}$. To estimate the derivative, let $\widetilde{\EE}$ denote expectation under the probability measure with density $|R_j(t)|^{-2}/\EE|R_j(t)|^{-2}$. The signs remain independent. Put $a_p=p^{-1/2-it}$ and $L_p=i(\log p)\eps_pa_p/(1+\eps_pa_p)$, so that $(R_j^{-1})'=R_j^{-1}\sum_{w<p\le w_j}L_p$. Under this measure, $\widetilde{\EE}\eps_p=-2\Re a_p/(1+|a_p|^2)$, and hence
\[
\begin{gathered}
\widetilde{\EE}L_p
=i(\log p)\big(a_p\widetilde{\EE}\eps_p+O(p^{-1})\big)
=O((\log p)/p),\\
\widetilde{\EE}|L_p|^2
\le\frac{(\log p)^2|a_p|^2}{(1-|a_p|)^2}
\ll\frac{(\log p)^2}{p}.
\end{gathered}
\]
Using independence under the tilted measure and the two bounds above, we obtain
\[
 \EE|(R_j^{-1})'(t)|^2
 \ll e^{3j}\left\{
       \left(\sum_{w<p\le w_j}\frac{\log p}{p}\right)^2
       +\sum_{w<p\le w_j}\frac{(\log p)^2}{p}\right\}
 \ll e^{3j}D_j^2.
\]

Substituting these moment bounds into the Sobolev inequality, we obtain expected suprema $O(e^j)$ for $|R_j|^2$ and $O(e^{3j})$ for $|R_j|^{-2}$ on each covering interval $I$. We apply Markov's inequality on each interval and sum over the intervals and $1\le j\le m$. For the upper boundary, we obtain
\[
 \begin{aligned}
 &\PP\left(\exists\,1\le j\le m,\ t\in[1,U]:
 |R_j(t)|>\sqrt T\,e^j e^{(\log T)^{1/100}}\log\log w_j\right)\\
 &\qquad\ll\sum_{j=1}^m
 \frac{UD_je^j}
      {T e^{2j}\exp\big(2(\log T)^{1/100}\big)(b+j)^2}\\
 &\qquad\ll\frac{(\log T)^{100}}b
      \exp\big(b-2(\log T)^{1/100}\big)=o(1).
 \end{aligned}
\]
For the lower boundary, the same argument gives
\[
 \begin{aligned}
 &\PP\left(\exists\,1\le j\le m,\ t\in[1,U]:
 |R_j(t)|<(\sqrt T\,e^j)^{-B_0}\right)\\
 &\qquad\le\sum_{j=1}^m\frac{CUD_je^{3j}}{T^{B_0}e^{2B_0j}}\\
 &\qquad\ll (\log T)^{100}e^bT^{1-B_0}
       \sum_{j\ge1}e^{-(2B_0-4)j}=o(1),
 \end{aligned}
\]
since $B_0>2$ and $b=(\log T)^{1/200}+O(1)$. This proves \eqref{eq:weakgood}.
\end{proof}

It follows that restricting the Euler product integral on $[1,U]$ to $\mathcal G_T(t)$ preserves it with probability $1-o(1)$. The strong event need not have high probability. We next show that passing to its rounded version loses only a negligible weighted contribution.

\subsection{Comparison of the barrier expectations}\label{sec:barrierexpectations}
We now show that replacing the weak barrier by the strong barrier changes the weighted mean by a negligible amount. Put
\[
 J_G=\int_1^U |F_z(\tfrac12+it)|^2K_T(t)^2
          \ind_{\mathcal G_T(t)}\,dt.
\]
For the strong barrier, we use a rounded event $\mathcal H^\#(t)\subseteq\mathcal H_T(t)$, defined in Section~\ref{sec:barrierprobabilities}, and put
\[
 J_H=\int_1^U |F_z(\tfrac12+it)|^2K_T(t)^2
          \ind_{\mathcal H^\#(t)}\,dt.
\]
Our aim is to prove
\[
 \EE(J_G-J_H)=o(\EE J_H).
\]

Write $Z_v=\prod_{p\le v}(1+p^{-1})$ for $v=w,z$, and put $Z_{w,z}=\prod_{w<p\le z}(1+p^{-1})$. Thus $Z_z=Z_wZ_{w,z}$. Let $F_w(s)=\prod_{p\le w}(1+\eps_pp^{-s})$. The factorization $F_z(\tfrac12+it)=F_w(\tfrac12+it)R_m(t)$ separates the primes at most $w$ from those in $(w,z]$.

We absorb the terminal weight $|R_m(t)|^2$ by a change of measure of Girsanov type; compare \cite[Lemma~4]{Harper}. Since $\EE|R_m(t)|^2=Z_{w,z}$, the density
\[
 \frac{d\PP_t}{d\PP}
 =\frac{|R_m(t)|^2}{Z_{w,z}}
 =\prod_{w<p\le z}
   \left(1+\frac{2\eps_p\sqrt p}{p+1}\cos(t\log p)\right)
\]
defines a probability measure on the signs with $w<p\le z$. The product form shows that these signs remain independent under $\PP_t$.

For any event $\mathcal A$ determined by the primes in $(w,z]$, independence of the small prime factors gives
\[
 \begin{aligned}
 \EE\left[|F_z(\tfrac12+it)|^2\ind_{\mathcal A}\right]
 &=\EE|F_w(\tfrac12+it)|^2\,
      \EE\left[|R_m(t)|^2\ind_{\mathcal A}\right]\\
 &=Z_wZ_{w,z}\PP_t(\mathcal A)
 =Z_z\PP_t(\mathcal A).
 \end{aligned}
\]
Both barrier events depend only on these large primes. Hence
\[
 \EE(J_G-J_H)
 =Z_z\int_1^U K_T(t)^2
   \left[\PP_t(\mathcal G_T(t))
        -\PP_t(\mathcal H^\#(t))\right]dt.
\]

The change of measure reduces our task to comparing two barrier probabilities. In the Steinhaus model, multiplication by $p^{-it}$ preserves the prime laws, so the one point calculation reduces to $t=0$; see \cite[Proposition~3.5]{HSX}. For Rademacher signs the tilted moments depend on $t$. We will compare both barrier probabilities with the same Gaussian main term, uniformly for $1\le t\le U$. Since this main term can tend to zero, the comparison must have relative error $o(1)$. We first prove the comparison estimates, then apply them to the two barriers in Section~\ref{sec:barrierprobabilities}.

\subsubsection{Relative Gaussian comparison}\label{sec:relative}
For one frequency, the increments to be compared are $\log|R_j(t)|-\log|R_{j-1}(t)|$ under the one point tilt; for two frequencies, we use the vector of increments at $t,s$. The next proposition converts an absolute Gaussian approximation for each increment into a relative comparison for rounded paths. Its key requirement is that the local errors remain summable after division by the Gaussian cell probabilities. We state it for general independent vectors and verify its hypotheses for the Euler products in Proposition~\ref{prop:phasecomparison}.

We write $\|A\|_{\mathrm{op}}=\sup_{\|x\|_2=1}\|Ax\|_2$ for the operator norm, and $I_d$ for the identity matrix. For symmetric matrices, $A\le B$ means $x^{\mathsf T}Ax\le x^{\mathsf T}Bx$ for all $x$. For $h>0$ and $x\in\mathbb R^d$, define $[x]_h$ coordinatewise by $([x]_h)_\ell=h\lfloor x_\ell/h+1/2\rfloor$ for $1\le\ell\le d$. For $y\in h\mathbb Z^d$, its cell is $Q_h(y)=\prod_{\ell=1}^d[y_\ell-h/2,y_\ell+h/2)$. Rectangles below are axis-parallel with arbitrary open or closed faces.

\begin{proposition}\label{prop:pathcomparison}
Fix $d\ge1$ and constants $0<c_0<C_0$. There are constants $K$ and $\rho_0>0$, depending only on $d,c_0,C_0$, with the following property.

Let $X_1,\ldots,X_n$ be independent vectors in $\RR^d$, with means $\mu_j$ and covariance matrices $\Sigma_j$. Let $G_j\sim N(\nu_j,V_j)$ be independent Gaussian reference vectors, and write $\widetilde G_j\sim N(\mu_j,\Sigma_j)$. Take $\eta_j,\varepsilon_j\ge0$, $0<h_j\le1$ and $R_j\ge1$. Assume that:
\begin{enumerate}[label=(\arabic*),leftmargin=*,itemsep=2pt,topsep=4pt,parsep=0pt]
\item For every $j$, $c_0I_d\le V_j\le C_0I_d$, and
\[
 \|\mu_j-\nu_j\|_2+\|\Sigma_j-V_j\|_{\mathrm{op}}\le\eta_j.
\]
\item For every $j$,
\[
 \sup_{Q\ \mathrm{a\ rectangle}}
 \left|\PP(X_j\in Q)-\PP(\widetilde G_j\in Q)\right|
 \le\varepsilon_j.
\]
\item The total error satisfies
\begin{equation}\label{eq:pathrho}
 \rho:=\sum_{j=1}^n
 \left((1+R_j^2)\eta_j+
       \varepsilon_jh_j^{-d}e^{KR_j^2}\right)\le\rho_0.
\end{equation}
\end{enumerate}
Then every nonnegative function $\Phi$ on $\prod_{j=1}^n h_j\mathbb Z^d$, supported on paths $(y_j)_j$ with $\|y_j-\nu_j\|_2\le R_j$ for every $j$, satisfies
\begin{equation}\label{eq:pathcomparison}
 \mathbb E\Phi\big(([X_j]_{h_j})_{j=1}^n\big)
 =\big(1+O_{d,c_0,C_0}(\rho)\big)
   \mathbb E\Phi\big(([G_j]_{h_j})_{j=1}^n\big).
\end{equation}
The error is uniform in $n$ and in $\Phi$.
\end{proposition}
\begin{proof}
We compare the two Gaussian laws on each allowed cell first, and then use its probability to turn the absolute approximation error into a relative one. Let $y\in h_j\mathbb Z^d$ satisfy $\|y-\nu_j\|_2\le R_j$. Since $h_j\le1$ and $R_j\ge1$, each $x\in Q_{h_j}(y)$ satisfies $\|x-\nu_j\|_2\ll_dR_j$. The Gaussian density with mean $\nu_j$ and covariance $V_j$ is therefore bounded below on this cell by $c e^{-KR_j^2}$. In particular,
\[
 q_{j,y}:=\mathbb P\{G_j\in Q_{h_j}(y)\}
 \ge c h_j^d e^{-KR_j^2}.
\]
The constants here and below depend only on $d,c_0,C_0$. Choose $\rho_0$ sufficiently small. Then $\eta_j\le\rho_0$ implies that the eigenvalues of $\Sigma_j$ also have fixed positive upper and lower bounds.

Let $f_j$ and $\widetilde f_j$ be the densities of $G_j$ and $\widetilde G_j$. In the Gaussian density formula, the changes in the mean and covariance give
\[
 \left|\log\frac{\widetilde f_j(x)}{f_j(x)}\right|
 \ll \eta_j\big(1+\|x-\nu_j\|_2^2\big)
 \ll (1+R_j^2)\eta_j
 \qquad (x\in Q_{h_j}(y)).
\]
This bound is small by \eqref{eq:pathrho}, so integration over the cell gives the same relative bound for its two Gaussian probabilities. The approximation of $X_j$ by $\widetilde G_j$ has the additional absolute error $\varepsilon_j$. Dividing this error by $q_{j,y}$, we obtain
\[
 \begin{aligned}
 \left|
 \frac{\mathbb P\{X_j\in Q_{h_j}(y)\}}{q_{j,y}}-1
 \right|
 &\ll (1+R_j^2)\eta_j+\frac{\varepsilon_j}{q_{j,y}}\\
 &\ll (1+R_j^2)\eta_j+
       \varepsilon_jh_j^{-d}e^{KR_j^2}
 =:r_j.
 \end{aligned}
\]
Increasing $K$ and decreasing $\rho_0$ if necessary, we have $\sum_jr_j\ll\rho$ and $r_j\le1/2$. These bounds are uniform over the allowed cells. By independence, the logarithm of the probability ratio for each path is $O(\sum_jr_j)$. Hence
\[
 \prod_{j=1}^n\mathbb P\{X_j\in Q_{h_j}(y_j)\}
 =\big(1+O(\rho)\big)
   \prod_{j=1}^n\mathbb P\{G_j\in Q_{h_j}(y_j)\}.
\]
Since $\Phi$ is nonnegative, the same relative bound holds after multiplying by $\Phi((y_j)_j)$ and summing over paths. This proves \eqref{eq:pathcomparison}.
\end{proof}

In particular, the proposition applies to indicators of any set of allowed rounded paths, even when their reference probability tends to zero. Removing the rounding requires a separate estimate for the Gaussian probability lost when the barriers are moved.

We next verify the comparison conditions for logarithmic Euler product increments. Under the tilted measures, their means and covariances are controlled by the prime sums from Section~\ref{sec:local}. Berry-Esseen estimates give the Gaussian approximation for each block. The following proposition combines these estimates with the path comparison.

We allow fixed complex phases in the Euler factors, so that the result also applies to the character twist in Section~\ref{sec:twists}.

Let $\mathcal P_1,\ldots,\mathcal P_n$ be disjoint finite sets of primes, and put $\mathcal P=\bigcup_j\mathcal P_j$. Fix deterministic complex numbers $\omega_p$ with $|\omega_p|=1$, and write
\[
 \mathcal R_\omega(t)
 =\prod_{p\in\mathcal P}
       (1+\varepsilon_p\omega_pp^{-1/2-it}).
\]
For $d=1$ or $2$ and $\mathbf t=(t_1,\ldots,t_d)\in\mathbb R^d$, the vector of logarithmic increments in the $j$th block is
\[
 \mathbf Y_j^\omega(\mathbf t)
 =\left(
   \sum_{p\in\mathcal P_j}
   \log|1+\varepsilon_p\omega_pp^{-1/2-it_\ell}|
   \right)_{\ell=1}^d.
\]
We write $Y_j^\omega(t_\ell)$ for the $\ell$th coordinate of this vector. Put $u_j=\sum_{p\in\mathcal P_j}p^{-2}$. The terms responsible for the frequency dependence are
\[
 \Gamma_j(\xi)
 =\sum_{p\in\mathcal P_j}
       \frac{\Re(\omega_p^2p^{-i\xi})}{p}.
\]
For two points, the ordinary difference frequency also contributes through
\[
 D_j(\xi)
 =\sum_{p\in\mathcal P_j}\frac{\cos(\xi\log p)}p.
\]

\begin{proposition}\label{prop:phasecomparison}
Put $\alpha_j=\sum_{p\in\mathcal P_j}p^{-1}$ and suppose $c_0\le\alpha_j\le C_0$ for fixed positive constants $c_0,C_0$. Let $d=1$ or $2$, and define
\begin{equation}\label{eq:phasebudget}
 \zeta_j(\mathbf t)=
 \begin{cases}
 u_j+|\Gamma_j(2t)|,&d=1,\quad\mathbf t=(t),\\
 u_j+|\Gamma_j(2t)|+|\Gamma_j(2s)|
       +|D_j(t-s)|+|\Gamma_j(t+s)|,
       &d=2,\quad\mathbf t=(t,s).
 \end{cases}
\end{equation}
Let $G_j\sim N(\alpha_j\mathbf1_d,(\alpha_j/2)I_d)$ be independent. Here $\mathbf1_d=(1,\ldots,1)^{\mathsf T}$ is a column vector. Take $0<h_j\le1$ and $R_j\ge1$, and put
\begin{equation}\label{eq:phaserho}
 \rho_\omega(\mathbf t)
 =\sum_{j=1}^n
 \left((1+R_j^2)\zeta_j(\mathbf t)
       +\tau_jh_j^{-d}e^{KR_j^2}\right),
\end{equation}
where $\tau_j=\sum_{p\in\mathcal P_j}p^{-3/2}$ and $K$ is a sufficiently large constant depending only on $d,c_0,C_0$. If $\rho_\omega(\mathbf t)$ is sufficiently small, then every nonnegative function $\Phi$ on $\prod_{j=1}^n h_j\mathbb Z^d$, supported on ${\|y_j-\alpha_j\mathbf1_d\|_2}\le R_j$ for every $j$, satisfies
\begin{equation}\label{eq:phasecomparison}
 \begin{aligned}
 &\mathbb E\left[
    \prod_{\ell=1}^d|\mathcal R_\omega(t_\ell)|^2\,
    \Phi\big(([\mathbf Y_j^\omega(\mathbf t)]_{h_j})_{j=1}^n\big)
    \right]\\
 &\qquad=Z_{\mathcal P}^{d}
       \big(1+O_{d,c_0,C_0}(\rho_\omega(\mathbf t))\big)
       \mathbb E\Phi\big(([G_j]_{h_j})_{j=1}^n\big),
 \end{aligned}
\end{equation}
where $Z_{\mathcal P}=\prod_{p\in\mathcal P}(1+p^{-1})$. The estimates are uniform in the phases, the frequencies, and the number of blocks subject to these hypotheses.
\end{proposition}
\begin{proof}
We use the probability measure with density proportional to $\prod_{\ell=1}^d|\mathcal R_\omega(t_\ell)|^2$. The signs remain independent under this measure, as do the block vectors. Write $\mathbb E_t$ and $\mathbb E_{t,s}$ for the one point and two point expectations.

For one frequency, the normalizing factor is $Z_{\mathcal P}$. At a prime $p$, the tilted density is $1+\varepsilon_pc_{p,\omega}(t)$, where $c_{p,\omega}(t)=2\Re(\omega_pp^{-it})/(\sqrt p(1+1/p))$. Applying the calculation of Lemma~\ref{lem:tilt} with $\omega_pp^{-it}$ in place of $p^{-it}$, we obtain
\[
 \begin{aligned}
 \mathbb E_tY_j^\omega(t)
 &=\alpha_j+\tfrac12\Gamma_j(2t)+O(u_j),\\
 \operatorname{Var}_tY_j^\omega(t)
 &=\tfrac12\alpha_j+\tfrac12\Gamma_j(2t)+O(u_j).
 \end{aligned}
\]
For two frequencies, use
\[
 2\Re(\omega_pp^{-it})\Re(\omega_pp^{-is})
 =\cos((t-s)\log p)+\Re(\omega_p^2p^{-i(t+s)})
\]
in the two point calculation of Lemma~\ref{lem:tilt}. This gives
\[
 \begin{aligned}
 \mathbb E_{t,s}Y_j^\omega(t)
 &=\alpha_j+\tfrac12\Gamma_j(2t)
       +D_j(t-s)+\Gamma_j(t+s)+O(u_j),\\
 \operatorname{Cov}_{t,s}(Y_j^\omega(t),Y_j^\omega(s))
 &=\tfrac12\big(D_j(t-s)+\Gamma_j(t+s)\big)+O(u_j).
 \end{aligned}
\]
The other mean follows by exchanging $t$ and $s$, and the diagonal variances have the one point main terms above. For the two point normalizing factor, independence gives
\[
 \frac{\mathbb E|\mathcal R_\omega(t)\mathcal R_\omega(s)|^2}
      {Z_{\mathcal P}^2}
 =\prod_{p\in\mathcal P}
       (1+c_{p,\omega}(t)c_{p,\omega}(s)).
\]
Taking logarithms and using the same phase identity, we find
\[
 \log\frac{\mathbb E|\mathcal R_\omega(t)\mathcal R_\omega(s)|^2}
               {Z_{\mathcal P}^2}
 =2\sum_j\big(D_j(t-s)+\Gamma_j(t+s)\big)
       +O\left(\sum_ju_j\right)
 =O(\rho_\omega).
\]

Under either tilt, each centered prime term has norm $O_d(p^{-1/2})$ and third absolute moment $O_d(p^{-3/2})$. For sufficiently small $\rho_\omega$, the block covariance eigenvalues are bounded above and away from zero by the moment formulas. Lemma~\ref{lem:BE} and Proposition~\ref{prop:pathcomparison} therefore apply with $\eta_j\ll\zeta_j$ and $\varepsilon_j\ll\tau_j$, after adjusting $K$. The relative error is $O(\rho_\omega)$, and multiplication by the normalizing factor proves \eqref{eq:phasecomparison}.
\end{proof}

\begin{remark}
The block conditions allow either sign of the frequencies but exclude the untwisted correlations at $t=s$ and $t=-s$, where the paths are identical in modulus. Nearby frequencies are treated in Section~\ref{sec:secondmoment}. Uniformity concerns deterministic frequencies, not a simultaneous continuum coupling. For other phases, a short interval central limit theorem still requires arithmetic transfer and fourth moment estimates in addition to this Euler product comparison.
\end{remark}

\subsubsection{Expected loss from strengthening the barrier}\label{sec:barrierprobabilities}
We return to the comparison of $\EE J_G$ and $\EE J_H$. We first compute the strong barrier mean by applying the preceding result to a rounded version of $\mathcal H_T(t)$.

For the logarithmic paths, write $L=\log T$, $a=L/2$ and $A=L^{1/100}$. Thus $b=L^{1/200}+O(1)$ and $b\ll A=o(a)$. Put $S_j(t)=\log|R_j(t)|$ and $X_j(t)=S_j(t)-S_{j-1}(t)$. Each increment is a sum over the primes in $(w_{j-1},w_j]$.

We use a rounded version of the strong event $\mathcal H_T(t)$. With $\Delta_j=e^{-2(b+j-1)}$, put ${[X]_j=\Delta_j\lfloor X/\Delta_j+1/2\rfloor}$. Any rounded partial sum differs from the original one by at most $\epsilon_g=\frac12\sum_{j=1}^m\Delta_j\ll e^{-2b}$. Define $\mathcal H^\#(t)$ by
\begin{equation}\label{eq:strongbarrier}
 -B_0(a+j)+\epsilon_g
 \le\sum_{i\le j}[X_i(t)]_i
 \le a+j-A-5\log(b+j)-\epsilon_g,
 \qquad 1\le j\le m.
\end{equation}
Since $\log\log w_j=b+j$, the inward shifts by $\epsilon_g$ imply
\[
 \mathcal H^\#(t)\subseteq\mathcal H_T(t)\subseteq\mathcal G_T(t).
\]
Conversely, if both logarithmic boundaries of $\mathcal H_T(t)$ are narrowed by $2\epsilon_g$, the resulting event implies $\mathcal H^\#(t)$. This is the only modification of the strong product bounds needed for the grid comparison.

Put $\alpha_j=\sum_{w_{j-1}<p\le w_j}p^{-1}$. The prime number theorem and partial summation give $\alpha_j=1+O(e^{-c\sqrt{\log w_{j-1}}})$, uniformly in $j$. Let $G_j\sim N(\alpha_j,\alpha_j/2)$ be independent. Write $p^\#$ for the probability that their rounded partial sums satisfy \eqref{eq:strongbarrier}. The following estimates follow from Proposition~\ref{prop:phasecomparison} with $\omega_p=1$. The one point estimate is needed here; we also record the two point estimate for Section~\ref{sec:barrierintegration}.

Here the abstract vectors are the increments $X_j(t)$ under the one point tilt, or $(X_j(t),X_j(s))$ under the two point tilt. Their Gaussian reference means and variances are $\alpha_j$ and $\alpha_j/2$. The grid size in Proposition~\ref{prop:phasecomparison} is $h_j=\Delta_j$, and its allowed increment radius is $C_{B_0}(a+j+1)$. The lower and upper path bounds imply this radius condition; the proof below checks the remaining error sum.

\begin{lemma}\label{lem:relative}
Put $r_0=e^{-b/4}$. Uniformly for $1\le t\le U$,
\begin{equation}\label{eq:onepoint}
 \mathbb E[|R_m(t)|^2\ind_{\cH^\#(t)}]
 =Z_{w,z}p^\#(1+O_{B_0}(e^{-b/2})),
\end{equation}
where $Z_{w,z}=\prod_{w<p\le z}(1+p^{-1})$. Uniformly for $1\le t,s\le U$ with $|t-s|\ge r_0$,
\begin{equation}\label{eq:twopoint}
 \mathbb E[|R_m(t)|^2|R_m(s)|^2
       \ind_{\cH^\#(t)}\ind_{\cH^\#(s)}]
 =Z_{w,z}^2(p^\#)^2(1+O_{B_0}(e^{-b/2})).
\end{equation}
\end{lemma}
\begin{proof}
For one frequency take $d=1$ and $d_*=1$; for two frequencies take $d=2$ and $d_*=\min\{1,|t-s|\}$. Since the frequencies are positive, the sums at $2t$, $2s$, and $t+s$ are all away from zero. Lemma~\ref{lem:primeosc} therefore bounds the block mean and covariance perturbations by
\begin{equation}\label{eq:etaj}
 \eta_j\ll\frac{e^{-b-j}}{d_*}
       +(1+U)e^{-ce^{(b+j-1)/2}}.
\end{equation}
The terms $u_j$ are absorbed by the second term.

A path satisfying \eqref{eq:strongbarrier} has each increment bounded by $C_{B_0}(a+j)$. In Proposition~\ref{prop:phasecomparison}, we therefore take grid size $h_j=\Delta_j$ and increment radius $C_{B_0}(a+j+1)$, after increasing the fixed constant. Since $\tau_j\ll w_{j-1}^{-1/2}$, the rectangular approximation part of \eqref{eq:phaserho} is bounded by
\[
 \begin{aligned}
 \sum_{j=1}^m
       \tau_j\Delta_j^{-d}e^{C_{B_0}(a+j)^2}
 &\ll\sum_{j\ge1}
       \exp\{-c e^{b+j-1}
              +C_{B_0}(a+j)^2+C_{B_0}(b+j)\}\\
 &\ll e^{-b}.
 \end{aligned}
\]
Here $e^b$ dominates every fixed power of $L$, and the exponential in $j$ dominates the polynomial terms uniformly in the number of blocks. The mean and covariance part satisfies
\[
 \begin{aligned}
 \sum_{j=1}^m(1+(a+j)^2)\eta_j
 &\ll_{B_0}\frac{(a+1)^2e^{-b}}{d_*}
       +(a+1)^2(1+U)e^{-c'e^{b/2}}\\
 &\ll_{B_0}e^{-b/2}.
 \end{aligned}
\]
The last inequality uses $d_*\ge e^{-b/4}$, $a=L/2$, $b=L^{1/200}+O(1)$, and $\log U\ll L$. Consequently $\rho_\omega=O_{B_0}(e^{-b/2})$.

Apply \eqref{eq:phasecomparison} to the indicator of the one or two strong barrier events. Under the reference law the $d$ coordinate paths are independent, so their probability is $(p^\#)^d$. This proves \eqref{eq:onepoint} and \eqref{eq:twopoint}.
\end{proof}

In particular, the estimates remain relative when $p^\#$ tends to zero.

The one point estimate and independence of the small prime factors give
\[
 \EE J_H
 =Z_zp^\#\left(1+O_{B_0}(e^{-b/2})\right)
       \int_1^U K_T(t)^2\,dt.
\]
Thus the main term is determined by $p^\#$, the survival probability of the Gaussian random walk.

Put
\[
 p_{z,T}=\frac1{\sqrt{2\pi}}
       \int_{-L/\sqrt{2M}}^{L/\sqrt{2M}}e^{-u^2/2}\,du.
\]
This is the same Gaussian probability as $p_x$, up to a relative error tending to zero:
\begin{equation}\label{eq:scalecomparison}
 p_{z,T}/p_x\to1.
\end{equation}
Indeed, $L=\log(x/y)+O(1)$, $\log(x/y)\to\infty$, and $M/\log\log x\to1$, so the two integration limits have ratio tending to one. A multiplicative change $1+o(1)$ in the limit of this Gaussian integral changes the integral by the same relative order, also when that limit tends to zero. If $s=\min\{1,L/\sqrt M\}$, we consequently have $p_{z,T}\asymp s$, $\vR\asymp ys$ and $\vR\gg y/\sqrt M$.

We now evaluate $p^\#$ and compare it with the Gaussian weak barrier probability. After subtracting the means, the Gaussian boundary estimates of Harper, Soundararajan and Xu reduce the calculation to a centered walk staying below height $L/2$ for $m\sim M$ steps. These estimates control both the change from the weak to the strong boundary and the shifts introduced by rounding.

\begin{lemma}\label{lem:ballot}
Let $p^\#$ be the probability that the partial sums of the rounded Gaussian increments satisfy the bounds in \eqref{eq:strongbarrier}. Then
\begin{equation}\label{eq:ballotasymp}
 p^\#\sim p_{z,T}.
\end{equation}
The same relative asymptotic holds for the Gaussian weak barrier corresponding to \eqref{eq:weakbarrier}, and for either barrier after moving each logarithmic boundary by at most $2\epsilon_g$.
\end{lemma}
\begin{proof}
Write $W_j=\sum_{i\le j}(G_i-\alpha_i)$ and $V_j=\operatorname{Var}(W_j)$. The prime sum estimate above gives $\sum_{i\le j}\alpha_i=j+O(e^{-ce^{b/2}})$ uniformly for $j\le m$. Thus $V_j=j/2+o(1)$, and the increment variances are bounded above and away from zero. Let $\mathcal U_0$ be the event that $W_j\le a$ for all $j\le m$, and put $q_m(a)=\mathbb P(\mathcal U_0)$. By \cite[Lemma~3.4]{HSX},
\[
 q_m(a)=\left(1+O\left(\frac{\log a}{a}\right)\right)
       \frac1{\sqrt{2\pi}}
       \int_{-a/\sqrt{V_m}}^{a/\sqrt{V_m}}e^{-u^2/2}\,du.
\]
Since $a=L/2$, $V_m\sim M/2$ and $m\sim M$, the comparison in \eqref{eq:scalecomparison} gives $q_m(a)\sim p_{z,T}$ and $q_m(a)\asymp\min\{1,L/\sqrt m\}$. It remains to compare the barrier events with $\mathcal U_0$.

After centering, the strong and weak upper boundaries are respectively
\[
 a-A-5\log(b+j)+O(1),\qquad
 a+A+\log(b+j)+O(1).
\]
Let $\mathcal U_H$ and $\mathcal U_G$ be the corresponding upper events, without the lower boundary. For large $L$, we have $\mathcal U_H\subseteq\mathcal U_0\subseteq\mathcal U_G$, so we first prove
\[
 \mathbb P(\mathcal U_G\setminus\mathcal U_H)=o(q_m(a)).
\]
On a path in $\mathcal U_G\setminus\mathcal U_H$, let $k$ be the first time the strong upper bound fails. Put $h_k=a+A+\log(b+k)+C$ and $D_k=C_1(A+\log(b+k))$, with $C,C_1$ sufficiently large. Then the first $k$ steps lie below $h_k$ and $W_k\in[h_k-D_k,h_k]$. To bound this event, split the first $k$ increments into three consecutive groups of comparable length, as in the estimate underlying \cite[(3.19)]{HSX}. The first group stays below $h_k$. In the last group, read the increments in reverse order with their signs reversed. The resulting partial sums are $W_{k-r}-W_k\le D_k$. The two groups are independent, so by \cite[Lemma~3.4]{HSX} the probability of these two constraints is $O(\min\{1,h_k/\sqrt k\}\min\{1,D_k/\sqrt k\})$.
Given both groups, the middle sum has a Gaussian density bounded by $O(k^{-1/2})$. Requiring $W_k$ to lie in the strip costs a further factor $O(\min\{1,D_k/\sqrt k\})$. Since $h_k\ll L+\log(b+k)$, the resulting bound is
\[
 \begin{aligned}
 &\mathbb P\bigl(\max_{j\le k}W_j\le h_k,\ W_k\ge h_k-D_k\bigr)\\
 &\qquad\ll\min\left\{1,\frac{L+\log(b+k)}{\sqrt k}\right\}
       \frac{(A+\log(b+k))^2}{k}.
 \end{aligned}
\]
For bounded $k$, the same bound follows from the trivial probability bound since $A\to\infty$.
Given the first $k$ steps, the remaining weak upper constraints imply
\[
 W_{k+r}-W_k\le2A+6\log(b+k)+\log(1+r)+O(1),
 \qquad 1\le r\le m-k.
\]
The remaining increments are independent centered Gaussians, with variances $\alpha_j/2$ in $[1/20,20]$ for large $x$. Since $\log(1+r)\le\log 2+\log r$, the boundary is at most $H_k+\log r$ for some $H_k\asymp A+\log(b+k)$. Thus \cite[Probability Result~1]{Harper} bounds the conditional probability by $O(\min\{1,H_k/\sqrt{m-k+1}\})$. For bounded $m-k$ this follows from the trivial bound. Using $H_k\asymp A+\log(b+k)$ and summing over the first crossing time gives
\[
 \mathbb P(\mathcal U_G\setminus\mathcal U_H)
 \ll\sum_{k=1}^m
   \min\left\{1,\frac{L+\log(b+k)}{\sqrt k}\right\}
   \frac{(A+\log(b+k))^3}{k\sqrt{m-k+1}}.
\]
For $k\le\min\{m,L^{10}\}$, $A+\log(b+k)\ll A$; splitting at $m/2$ bounds this part by $O(A^3\log L/\sqrt m)$. If $L^{10}<m$, the remaining part is at most
\[
 L^{1+3/100}\sum_{L^{10}<k\le m}
       \frac{(\log k)^4}{k^{3/2}\sqrt{m-k+1}}
 \ll\frac{L^{-4+3/100}(\log L)^4}{\sqrt m}.
\]
Thus $\mathbb P(\mathcal U_G\setminus\mathcal U_H)\ll L^{1/25}/\sqrt m$. The near full range gives $L\le(\tfrac12+o(1))M$, hence $L\ll m$. The ratio of this bound to $q_m(a)$ is $O(L^{-24/25})$ when $L\le\sqrt m$, and $O(m^{-23/50})$ otherwise. Hence
\[
 \mathbb P(\mathcal U_H)\sim q_m(a)
 \sim\mathbb P(\mathcal U_G).
\]

We now restore the common lower boundary $-B_0a-(B_0+1)j+O(1)$. Its violation at time $k$ implies $W_k<-(a+k)$ for large $a$. To retain the factor $q_m(a)$ in the error, we keep the condition on the remaining path. Suppose first that $a\le\sqrt m$. For $k\le m/2$ and $u<0$, independence of the remaining increments and the flat boundary estimate give
\[
 \mathbb P(\mathcal U_0\mid W_k=u)
 \ll\min\left\{1,\frac{a+|u|+1}{\sqrt m}\right\}.
\]
Using the Gaussian density of $W_k$, we obtain
\[
 \begin{aligned}
 \mathbb P(\mathcal U_0,\ W_k<-(a+k))
 &\ll\frac1{\sqrt{mV_k}}
       \int_{-\infty}^{-(a+k)}(a+|u|+1)e^{-u^2/(2V_k)}\,du\\
 &\ll\frac a{\sqrt m}e^{-c(a+k)}.
 \end{aligned}
\]
The last bound uses $V_k\asymp k$ and $(a+k)^2/k\ge a+k$, with the polynomial factors absorbed into the exponential. Summing over $k\le m/2$ gives $O(e^{-ca}q_m(a))$. For $k>m/2$, the unconditional Gaussian tails sum to $O(e^{-c(a+m)})$, which is smaller than the same bound. If $a>\sqrt m$, then $q_m(a)\asymp1$ and the unconditional tails suffice for all $k$. We have therefore proved
\[
 \mathbb P\bigl(\mathcal U_0,\ 
    W_k<-B_0a-(B_0+1)k+O(1)\text{ for some }k\le m\bigr)
 \ll e^{-ca}q_m(a).
\]
Since $\mathcal U_H\subseteq\mathcal U_0$, this bounds the loss from the strong upper event. For the weak upper event, the loss is at most this quantity plus $\mathbb P(\mathcal U_G\setminus\mathcal U_0)=o(q_m(a))$. Both full barrier probabilities are consequently asymptotic to $q_m(a)$.

All these estimates allow bounded shifts of the barriers. Each rounded partial sum differs from the unrounded one by at most $\epsilon_g$, so the rounded events lie between unrounded events with boundaries moved by at most $2\epsilon_g$. These have the same asymptotic $q_m(a)\sim p_{z,T}$, proving the lemma.
\end{proof}

In particular, $p^\#\asymp\min\{1,L/\sqrt M\}$. The error in this Gaussian ballot asymptotic is used only as an $o(1)$ relative error. It is separate from the quantitative relative error in \eqref{eq:onepoint}; no uniform $M^{-c}$ rate is asserted when $T\to\infty$ arbitrarily slowly.

It follows that
\[
 \EE J_H\sim Z_zp_{z,T}\int_1^U K_T(t)^2\,dt.
\]
To compare the weak barrier mean with this quantity, we bound the original weak event between two rounded events with slightly shifted boundaries. The relative comparison applies to both bounding events, and their Gaussian probabilities have the same main term $p_{z,T}$. This gives the required uniform estimate for the difference.

\begin{proposition}\label{lem:barriergap}
We have
\[
 \EE(J_G-J_H)=o(\EE J_H).
\]
More precisely, uniformly for $1\le t\le U$,

\begin{equation}\label{eq:barrierdelete}
 \EE\left[
  |F_z(\tfrac12+it)|^2
  (\ind_{\mathcal G_T(t)}-\ind_{\mathcal H^\#(t)})
 \right]
 =o(Z_zp_{z,T}).
\end{equation}
\end{proposition}
\begin{proof}
Let $\mathcal G_T^{\#,-}(t)$ and $\mathcal G_T^{\#,+}(t)$ be the events obtained by imposing the logarithmic weak bounds from \eqref{eq:weakbarrier} on the rounded partial sums, with both boundaries moved inward and outward by $\epsilon_g$, respectively. The deterministic rounding bound gives
\[
 \mathcal G_T^{\#,-}(t)
 \subseteq\mathcal G_T(t)
 \subseteq\mathcal G_T^{\#,+}(t).
\]
Both rounded events restrict every increment to $O_{B_0}(a+j)$. We compare each event with its Gaussian counterpart by Proposition~\ref{prop:phasecomparison}, as in \cref{lem:relative}. By \cref{lem:ballot}, their Gaussian probabilities are both $(1+o(1))p_{z,T}$. The inclusions therefore imply, under the terminal one point tilted probability $\PP_t$,
\[
 \PP_t(\mathcal G_T(t))=(1+o(1))p_{z,T}.
\]
Also, \eqref{eq:onepoint} and \cref{lem:ballot} give $\PP_t(\mathcal H^\#(t))=(1+o(1))p_{z,T}$. Since the strong event is contained in the weak event, their tilted probabilities differ by $o(p_{z,T})$.

The change of measure identity at the start of Section~\ref{sec:barrierexpectations} identifies the weighted difference in \eqref{eq:barrierdelete} with
\[
 Z_z\left(\PP_t(\mathcal G_T(t))
             -\PP_t(\mathcal H^\#(t))\right)
 =o(Z_zp_{z,T}),
\]
uniformly for $1\le t\le U$.

After multiplying \eqref{eq:barrierdelete} by $K_T(t)^2$ and integrating over $1\le t\le U$, we obtain
\[
 \EE(J_G-J_H)
 =o\left(Z_zp_{z,T}\int_1^U K_T(t)^2\,dt\right)
 =o(\EE J_H).\qedhere
\]
\end{proof}

Since $J_G-J_H\ge0$, we apply Markov's inequality to this difference and use $\EE(J_G-J_H)=o(\EE J_H)$ to obtain $J_G-J_H=o_{\mathbb P}(\EE J_H)$. We may therefore prove concentration under the strong barrier in Section~\ref{sec:barrierintegration} and transfer it to the weakly restricted integral.

\begin{remark}\label{rem:steinhausbarrier}
The Steinhaus analogue has normalizing factor $\prod_{p\le v}(1-p^{-1})^{-1}$, tilted block means $\alpha_j+O(u_j)$, variances $\alpha_j/2+O(u_j)$ and the same third moment bound, where $u_j=\sum_{w_{j-1}<p\le w_j}p^{-2}$. The same argument gives the expected loss estimate. Translation invariance permits evaluation at zero frequency there; the relative comparison also handles the Rademacher frequency dependence.
\end{remark}

\subsection{Concentration under the strong barrier and completion}\label{sec:barrierintegration}
Section~\ref{sec:barrierexpectations} shows that passing from $J_G$ to $J_H$ loses only a negligible weighted mean. We now prove concentration of $J_H$. Write $\mu^\#=Z_zp^\#$. Since $\int_1^U K_T(t)^2\,dt\asymp T$, it is enough to show that
\[
 \EE\left|J_H-\mu^\#\int_1^U K_T(t)^2\,dt\right|^2
 =o\bigl((T\mu^\#)^2\bigr).
\]
We first estimate this second moment, and then remove the barriers and the frequency restrictions to complete the conditional variance estimate.

\subsubsection{The second moment under the strong barrier}\label{sec:secondmoment}
The second moment is
\[
 \begin{aligned}
 \EE J_H^2
 &=\int_1^U\int_1^U K_T(t)^2K_T(s)^2\\
 &\qquad\times
 \EE\left[
 |F_z(\tfrac12+it)F_z(\tfrac12+is)|^2
 \ind_{\mathcal H^\#(t)}\ind_{\mathcal H^\#(s)}
 \right]\,dt\,ds.
 \end{aligned}
\]
The barrier events only involve primes in $(w,z]$. Since $F_z(\tfrac12+it)=F_w(\tfrac12+it)R_m(t)$, independence separates the integrand into two factors:
\begin{equation}\label{eq:secondmomentfactorization}
 \begin{aligned}
 &\EE\left[|F_z(\tfrac12+it)F_z(\tfrac12+is)|^2
       \ind_{\mathcal H^\#(t)}\ind_{\mathcal H^\#(s)}\right]\\
 &\quad=\EE|F_w(\tfrac12+it)F_w(\tfrac12+is)|^2\\
 &\qquad\quad\times
 \EE\left[|R_m(t)R_m(s)|^2
       \ind_{\mathcal H^\#(t)}\ind_{\mathcal H^\#(s)}\right].
 \end{aligned}
\end{equation}
We split the pairs at $r_0=e^{-b/4}$. For $|t-s|\ge r_0$, Lemma~\ref{lem:relative} replaces the second factor by $Z_{w,z}^2(p^\#)^2$ with a small relative error. The first factor need not separate pointwise; we control its integral by a Walsh expansion. For $|t-s|<r_0$, the initial large prime blocks are also correlated, and the path comparison must start later. We prepare the small prime estimate first, then the conditional comparison for these nearby frequencies.

These two estimates lead to Proposition~\ref{prop:concentration}. Its relative error terms $e^{-b/2}$, $e^{5b}/T$ and $e^{-A}$ come from the large prime comparison, the fluctuation of the small prime integral and the strong barrier bound for nearby frequencies, respectively.

\paragraph{The small prime integral.}\label{sec:smallvariance}
For separated frequencies, \eqref{eq:secondmomentfactorization} leaves the second moment of the small prime integral. The following bound controls its fluctuation after integration.

\begin{lemma}\label{lem:smallvariance}
Suppose $g$ is integrable on $\RR$ and of bounded variation, with $\|g\|_1\ll T$ and $\|g\|_\infty+\operatorname{TV}(g)\ll1$. Then
\begin{equation}\label{eq:smallvariance}
 \EE\left|\int g(t)|F_w(\tfrac12+it)|^2\dd t-Z_w\int g(t)\dd t\right|^2
 \ll T(\log(2w))^7.
\end{equation}
\end{lemma}
\begin{proof}
Let $\mathcal S_w$ be the square-free positive integers supported on primes at most $w$, and put $\eps_d=\prod_{p\mid d}\eps_p$ for $d\in\mathcal S_w$. From \eqref{eq:localdensity} we obtain
\begin{equation}\label{eq:walshsmall}
 \frac{|F_w(\tfrac12+it)|^2}{Z_w}
 =\sum_{d\in\mathcal S_w}\eps_d\frac{\sqrt d}{\sigma(d)}
       \sum_{rs=d}e^{-it\log(r/s)},
\end{equation}
where $\sigma(d)=\sum_{r\mid d}r$. The term $d=1$ is constant. Put $\widehat g(\xi)=\int_{\RR}g(t)e^{-it\xi}\dd t$. The assumptions imply $|\widehat g(\xi)|\ll\min\{T,|\xi|^{-1}\}$ for $\xi\ne0$.

We first use orthogonality of the variables $\eps_d$ and then apply Cauchy's inequality to each inner sum:
\[
 \begin{aligned}
 &Z_w^{-2}\EE\left|
   \int g(t)\bigl(|F_w(\tfrac12+it)|^2-Z_w\bigr)\dd t
   \right|^2\\
 &=\sum_{\substack{d\in\mathcal S_w\\d>1}}
     \frac{d}{\sigma(d)^2}
     \left|\sum_{rs=d}\widehat g(\log(r/s))\right|^2\\
 &\le\sum_{\substack{d\in\mathcal S_w\\d>1}}
     \frac{\tau(d)}d
     \sum_{rs=d}|\widehat g(\log(r/s))|^2\\
 &\ll T\sum_{\substack{r,s\in\mathcal S_w\\(r,s)=1\\r\ne s}}
       \frac{\tau(r)\tau(s)}{rs|\log(r/s)|}.
 \end{aligned}
\]
We used $d/\sigma(d)^2\le1/d$, $|\widehat g(\xi)|^2\ll T/|\xi|$, and $\tau(d)=\tau(r)\tau(s)$ when $rs=d$ is square-free.

Dropping coprimality increases the last sum. The part with $r/s\notin(1/2,2)$ is $O((\log(2w))^4)$ by the Euler product for $\sum_{r\in\mathcal S_w}\tau(r)/r$. For comparable $r,s$, use $2\tau(r)\tau(s)\le\tau(r)^2+\tau(s)^2$ and symmetry. Since $|\log(r/s)|\asymp|r-s|/r$,
\[
 \sum_{\substack{r/2<s<2r\\s\ne r}}
       \frac1{rs|\log(r/s)|}
 \ll\frac{\log(2r)}r.
\]
This part is at most
\[
 \sum_{r\in\mathcal S_w}\frac{\tau(r)^2\log(2r)}r
 \ll(\log(2w))^5,
\]
as follows by expanding $\log r=\sum_{p\mid r}\log p$ in the Euler product $\prod_{p\le w}(1+4/p)$. Finally multiply by $Z_w^2\ll(\log(2w))^2$.
\end{proof}

For $g(t)=K_T(t)^2\ind_{[1,U]}(t)$, the mean is $Z_w\int g\asymp TZ_w$. Since $Z_w\asymp\log w=e^b$, the variance bound is small relative to the square of this mean:
\[
 \frac{T(\log(2w))^7}{(TZ_w)^2}
 \asymp\frac{e^{5b}}T=o(1).
\]

\paragraph{Restarting the comparison for nearby frequencies.}
Let $d=|t-s|$ with $0<d<r_0$. The ordinary two point terminal moment is controlled once $d\log w_k\gg1$. For the relative path comparison, \eqref{eq:restarterror} below gives
\[
 \sum_{r=1}^{m-k}\bigl(1+(a+k+r)^2\bigr)\eta_{k+r}
 \ll_{B_0}\frac{(a+k+1)^2}{d\log w_k}.
\]
At the natural separation scale $\log w_k\asymp d^{-1}$, this bound need not tend to zero. We therefore keep that scale for the terminal moment but restart the path comparison only when $\log w_k\gg d^{-2}$. The error then becomes $O((a+k+1)^2d)=o(1)$. Take a sufficiently large fixed constant $C$ and set
\begin{align}\label{eq:cutpoints}
 k_0&=\min\{m,\max(0,\lceil\log(C/d)-b\rceil)\},\notag\\
 k_1&=\min\{m,\max(0,\lceil\log(C/d^2)-b\rceil)\}.
\end{align}
At the first cutoff, $d\log w_{k_0}\gg1$ unless $k_0=m$. At the second, $\log w_{k_1}\gg d^{-2}$ unless $k_1=m$. These cutoffs give the decomposition
\[
 R_m(t)=R_{k_0}(t)\,
        \frac{R_{k_1}(t)}{R_{k_0}(t)}\,
        \frac{R_m(t)}{R_{k_1}(t)}.
\]
On the first segment, we use the strong upper bound for one path and keep the prefix barrier probability for the other. On the middle segment, we drop the path constraints and estimate only the terminal two point moment. On the final segment, we condition on the signs already exposed and restart the relative comparison, which supplies a Gaussian barrier probability for each of the two remaining paths. The next lemma makes this last step precise.

\begin{lemma}\label{lem:restart}
Let $1\le t,s\le U$ with $0<d=|t-s|<r_0$, and let $k=k_1(d)$ be as in \eqref{eq:cutpoints}. Put $\cF_{w_k}=\sigma(\eps_p:p\le w_k)$. Under the original Rademacher law, the following bound holds uniformly over all initial sign choices satisfying the first $k$ inequalities in \eqref{eq:strongbarrier} at both $t$ and $s$.
\begin{equation}\label{eq:restart}
 \begin{aligned}
 &\EE\left[
 \left|\frac{R_m(t)R_m(s)}{R_k(t)R_k(s)}\right|^2
 \ind_{\mathcal H^\#(t)}\ind_{\mathcal H^\#(s)}
 \,\middle|\,\cF_{w_k}
 \right]\\
 &\qquad\ll_{B_0}
 e^{2(m-k)}
 \min\left\{1,\frac{a+k+1}{\sqrt{m-k+1}}\right\}^{2}.
 \end{aligned}
\end{equation}
\end{lemma}
\begin{proof}
When $k=m$, the tail products are empty and the ratios equal one, so \eqref{eq:restart} is immediate. We may therefore assume $k<m$. Condition on the prime signs up to $w_k$. On initial paths admitted by the strong barriers, both logarithmic paths have heights of absolute value $O_{B_0}(a+k+1)$. After subtracting these heights and the reference Gaussian means, the future upper barriers are at most $C_{B_0}(a+k+1)$. On paths satisfying the barrier constraints, the increment in the $r$th remaining block has size at most $C_{B_0}(a+k+r+1)$.

The future blocks remain independent under the terminal two point tilt. Using \eqref{eq:etaj} in Proposition~\ref{prop:phasecomparison}, the mean and covariance part of the error is bounded by
\begin{equation}\label{eq:restarterror}
 \begin{aligned}
 &\sum_{r=1}^{m-k}
   \bigl(1+(a+k+r)^2\bigr)\eta_{k+r}\\
 &\quad\ll_{B_0}
   \frac1{d\log w_k}
   \sum_{r\ge1}(a+k+r+1)^2e^{-r}\\
 &\quad\ll_{B_0}
   \frac{(a+k+1)^2}{d\log w_k}
   \ll_{B_0}(a+k+1)^2d=o(1).
 \end{aligned}
\end{equation}
The prime number theorem remainders are absorbed here as in the proof of Lemma~\ref{lem:relative}. The last bound is uniform: $k\le2\log(1/d)+O(1)$, $d<e^{-b/4}$, and $a$ is a fixed power of $b$. The rectangular approximation errors are also summable by that proof, since the first remaining prime block starts still higher.

Apply the relative comparison to the event with both future barrier constraints. Under the Gaussian reference law the two future paths are independent. We may now drop their lower constraints and enlarge the upper constraints to the constant height above. The bound for the survival probability of the Gaussian random walk in \cite[Lemma~3.4]{HSX} gives one factor $\min\{1,(a+k+1)/\sqrt{m-k+1}\}$ for each path.

By Lemma~\ref{lem:primeosc} and $d\log w_k\gg1$, the terminal two point normalizing factor for the remaining blocks is $O(e^{2(m-k)})$. Multiplying this factor by the two Gaussian probability bounds proves \eqref{eq:restart}.
\end{proof}

Both probability factors in \eqref{eq:restart} are needed for comparison with the square $(p^\#)^2$ of the one point probability.

We also need a one point estimate for an initial segment. Let $\mathcal H^\#_{\le k}(u)$ denote the event obtained by keeping only the first $k$ inequalities in \eqref{eq:strongbarrier}; for $k=0$ it is unrestricted. Uniformly for $1\le u\le U$ and $0\le k\le m$,
\[
 \EE\left[|R_k(u)|^2\ind_{\mathcal H^\#_{\le k}(u)}\right]
 \ll_{B_0}e^k
       \min\left\{1,\frac{a+1}{\sqrt{k+1}}\right\}.
\]
For $k\ge1$, apply Proposition~\ref{prop:phasecomparison} to the first $k$ blocks. Their error sum is bounded by the full one point sum in the proof of Lemma~\ref{lem:relative}, and their normalizing factor is $\asymp e^k$. After centering, the Gaussian prefix event implies $\sum_{i\le j}(G_i-\alpha_i)\le a+O(1)$ for $j\le k$. Dropping the lower boundary and applying \cite[Lemma~3.4]{HSX} bounds the survival probability of the Gaussian random walk by the minimum above. The case $k=0$ is immediate.

We can now combine the estimates for the two frequency ranges.

\begin{proposition}\label{prop:concentration}
With the notation above,
\begin{equation}\label{eq:concentration}
 \EE\left|J_H-\mu^\#\int_1^U K_T(t)^2\,dt\right|^2
 \ll_{B_0}
 \left(e^{-b/2}+\frac{e^{5b}}T+e^{-A}\right)(T\mu^\#)^2.
\end{equation}
\end{proposition}
\begin{proof}
The one point comparison, together with independence across the cutoff $w$, implies
\[
 \EE J_H
 =\mu^\#(1+O(e^{-b/2}))\int_1^U K_T(t)^2\,dt.
\]

First consider $|t-s|\ge r_0$. Put $I_w=\int_1^U K_T(t)^2|F_w(\tfrac12+it)|^2\,dt$. Apply \eqref{eq:secondmomentfactorization} and \eqref{eq:twopoint} to the integrand, then enlarge the domain of the small prime integral using nonnegativity. The contribution of these pairs to $\EE J_H^2$ is at most
\[
 (1+Ce^{-b/2})Z_{w,z}^2(p^\#)^2\,\EE I_w^2.
\]
Lemma~\ref{lem:kernel} allows us to apply Lemma~\ref{lem:smallvariance} with $g(t)=K_T(t)^2\ind_{[1,U]}(t)$, so
\[
 \EE I_w^2
 =Z_w^2\left(\int_1^U K_T(t)^2\,dt\right)^2
  +O\bigl(T(\log w)^7\bigr).
\]
The relative small prime error is $O(e^{5b}/T)$, as computed after Lemma~\ref{lem:smallvariance}.

We next treat $d=|t-s|<r_0$. Let $k_0,k_1$ be as in \eqref{eq:cutpoints}, and write $\pi_j=\min\{1,(a+1)/\sqrt{j+1}\}$ and $q_j=\min\{1,(a+j+1)/\sqrt{m-j+1}\}$. At time $k_0$, the strong upper boundary implies
\begin{equation}\label{eq:earlybarrier}
 |R_{k_0}(t)|^2
 \le\frac{Te^{2k_0}e^{-2A}}{(b+k_0)^{10}}.
\end{equation}
This also holds for $k_0=0$, since $R_0=1$ and $Te^{-2A}/b^{10}\to\infty$.

If $k_0<k_1$, then $d\log w_{k_0}\gg1$. Expanding the two point normalizing factor and applying Lemma~\ref{lem:primeosc} at $t-s$ and $t+s$, we find
\[
 \EE\left|
 \frac{R_{k_1}(t)R_{k_1}(s)}
      {R_{k_0}(t)R_{k_0}(s)}
 \right|^2
 \ll e^{2(k_1-k_0)}.
\]
The same bound holds for an empty segment.

Condition on the signs up to $w_{k_1}$. For every admitted pair of initial paths, Lemma~\ref{lem:restart} bounds the remaining contribution by $C_{B_0}e^{2(m-k_1)}q_{k_1}^2$. We may then discard the barrier constraints between $k_0$ and $k_1$. On the earlier blocks, use \eqref{eq:earlybarrier} for the $t$ path and retain the prefix event for the $s$ path. The prefix estimate above bounds that part by $C_{B_0}e^{k_0}\pi_{k_0}$. The primes at most $w$ contribute at most $C(\log w)^6$ to the fourth moment. Independence of these prime segments yields
\begin{equation}\label{eq:nearpoint}
 \begin{aligned}
 &\EE\left[
 |F_z(\tfrac12+it)F_z(\tfrac12+is)|^2
 \ind_{\mathcal H^\#(t)}\ind_{\mathcal H^\#(s)}
 \right]\\
 &\quad\ll_{B_0}
 (\log w)^6
 \frac{Te^{2k_0}e^{-2A}}{(b+k_0)^{10}}\,
 e^{k_0}\pi_{k_0}\,
 e^{2(k_1-k_0)}e^{2(m-k_1)}q_{k_1}^2\\
 &\quad=
 \frac{Te^{-2A}(\log w)^6e^{2m+k_0}}{(b+k_0)^{10}}\,
 \pi_{k_0}q_{k_1}^2.
 \end{aligned}
\end{equation}
When $k_1=m$, the final segment is empty and $q_m=1$. When $k_0=m$, both the intermediate and final segments are empty.

It remains to integrate in $d$. We claim that
\begin{equation}\label{eq:nearseries}
 \int_0^{r_0}
 \frac{e^{k_0(d)}}{(b+k_0(d))^{10}}\,
 \pi_{k_0(d)}q_{k_1(d)}^2\,dd
 \ll\pi_m^2.
\end{equation}
We first treat $k_0=0$ and $k_0=m$. If $k_0=0$, then $d\ge Ce^{-b}$ and $k_1\le b+O(1)$. Since $b=o(a)$ and $b=o(m)$, we have $q_{k_1}^2\ll\pi_m^2$. This set has length at most $r_0$, so its contribution is $O(r_0b^{-10}\pi_m^2)$. If $k_0=m$, then $k_1=m$ and the set has length $O(e^{-b-m})$. Integrating over this set gives
\[
 \begin{aligned}
 \int_{\substack{0<d<r_0\\k_0(d)=m}}
       \frac{e^m\pi_m}{(b+m)^{10}}\,dd
 &\ll e^{-b}(b+m)^{-10}\pi_m\\
 &\ll\pi_m^2,
 \end{aligned}
\]
where we used $q_m=1$ and $\pi_m\ge(m+1)^{-1/2}$.

For $1\le k_0=k<m$, the corresponding set of distances has length $O(e^{-b-k})$, and $k_1\le\min\{m,b+2k+O(1)\}$. The remaining contribution is therefore at most
\[
 Ce^{-b}\sum_{1\le k<m}(b+k)^{-10}\pi_kq_{k_1}^2.
\]
If $m\le16(a+b+1)^2$, then $\pi_m\asymp1$, and the bound follows from $\pi_k,q_{k_1}\le1$. Otherwise split at $k=m/8$. For $k\le m/8$, the remaining length $m-k_1+1$ is comparable to $m$, so
\[
 q_{k_1}^2
 \ll\pi_m^2\left(1+\frac{b+k}{a+1}\right)^2.
\]
For the terms with $k\le m/8$, we use $\pi_k\le1$ and obtain
\[
 Ce^{-b}\pi_m^2
 \sum_{k\ge1}
 \left((b+k)^{-10}
       +\frac{(b+k)^{-8}}{(a+1)^2}\right)
 \ll\pi_m^2.
\]
For $k>m/8$, we insert $\pi_k\ll\pi_m$ and $q_{k_1}\le1$ into the sum to obtain
\[
 Ce^{-b}\pi_m\sum_{k>m/8}(b+k)^{-10}
 \ll e^{-b}\pi_m m^{-9}
 \ll\pi_m^2.
\]
This proves \eqref{eq:nearseries}.

Finally, $\int_{\RR}K_T(t)^2K_T(t+d)^2\,dt\ll T$. Integrate \eqref{eq:nearpoint} first in $t$ and then in $d$, and apply \eqref{eq:nearseries}. The near frequency contribution is therefore at most
\[
 C_{B_0}T^2e^{-2A}(\log w)^6e^{2m}\pi_m^2
 \ll_{B_0}e^{4b-2A}(T\mu^\#)^2
 \ll_{B_0}e^{-A}(T\mu^\#)^2,
\]
using $Z_z^2\asymp e^{2b+2m}$ and $p^\#\asymp\pi_m$. Inserting the two second moment bounds and the first moment estimate into the centered square proves \eqref{eq:concentration}.
\end{proof}

\subsubsection{Completion of the concentration argument}\label{sec:CVcompletion}
We first prove the tail estimate \eqref{eq:tail-negligible}, which was deferred in Section~\ref{sec:CVtransfer}. We then remove the barriers on $[1,U]$ and use \eqref{eq:frequency-reduction} to recover the full conditional variance.

The exact second moment alone gives a relative tail bound $O(T/(Up_{z,T}))$, which need not vanish. We therefore impose barriers on dyadic tail intervals to keep the Gaussian probability factor.

Put $R_j=2^jU$ for $j\ge0$ and $s_z(R)=\min\{1,\log R/\sqrt M\}$. For each $R=R_j$ with $\log R<\sqrt M$, keep $z$ fixed and repeat the prime partition and barrier definitions with $R$ in place of $T$. Thus the new block count is $m_R=\lceil M-(\log R)^{1/200}\rceil$, the new initial time is $b_R=M-m_R$, and the prime endpoints are $\exp(e^{b_R+\ell})$ for $0\le\ell\le m_R$. Let $\mathcal G_R(t)$ denote the resulting weak barrier. These choices concern the barrier estimates; the tail integral still has the original kernel $K_T$.

We check why the earlier estimates apply uniformly on $[R,2R]$. First, $m_R=M-o(M)$, since $(\log R)^{1/200}\le M^{1/400}$. The proof of Proposition~\ref{lem:weakbarrier}, with an interval of length $R$, gives failure probability $\ll\exp(-c(\log R)^{1/100})$. In the one point comparison of Lemma~\ref{lem:relative}, replace $L,b,a$ by $\log R,b_R,(\log R)/2$. Its two error sums are still $O(e^{-b_R/2})$: the logarithm of the lower prime endpoint is $e^{b_R}$, which dominates every fixed power of $\log R$, and the frequencies are at most $2R$. Lemma~\ref{lem:ballot} then bounds the Gaussian weak barrier probability by a constant times $s_z(R)$. Applying the same change of measure and rounding argument as in Proposition~\ref{lem:barriergap} proves the required weighted mean bound. All these estimates are uniform as $R\ge U\to\infty$. For $\log R\ge\sqrt M$, take $\mathcal G_R$ unrestricted and use \eqref{eq:Fzsecond}. Thus
\begin{align}\label{eq:tailbarrier}
 \PP(\mathcal G_R\text{ fails somewhere on }[R,2R])
 &\ll e^{-c(\log R)^{1/100}},\notag\\
 \EE\left[|F_z(\tfrac12+it)|^2\ind_{\mathcal G_R(t)}\right]
 &\ll Z_zs_z(R)
 \qquad(R\le t\le2R).
\end{align}
The sum of the failure probabilities over $R_j$ tends to zero. On their common event, the positive tail equals the sum of the restricted integrals. Integrate the weighted mean bound in \eqref{eq:tailbarrier} over each dyadic interval and use $K_T(t)^2\ll T^2/t^2$. Since $s_z(R_j)/s_z(T)\le\log R_j/\log T$, summing the resulting bounds geometrically yields
\[
 \frac{T^2Z_z\sum_{j\ge0}s_z(R_j)/R_j}
      {TZ_zp_{z,T}}
 \ll\frac TU\frac{\log U+1}{\log T}
 \ll L^{-100}.
\]
We apply Markov's inequality to the sum of the restricted tail integrals and add the probability that a tail barrier fails. The preceding estimate and the vanishing sum of failure probabilities show that the positive tail is $o_{\PP}(TZ_zp_{z,T})$. By evenness, the negative tail satisfies the same bound. Since $p_{z,T}/p_x\to1$ by \eqref{eq:scalecomparison}, this proves \eqref{eq:tail-negligible} and completes the reduction \eqref{eq:frequency-reduction}.

\begin{theorem}\label{thm:CV}
We have
\begin{equation}\label{eq:CVintegral}
 \frac{\displaystyle\int_{\RR}|F_z(\tfrac12+it)|^2K_T(t)^2\,dt}
      {\displaystyle Z_zp_{z,T}\int_{\RR}K_T(t)^2\,dt}
 \xrightarrow{\PP}1.
\end{equation}
Consequently
\begin{equation}\label{eq:CVgoal}
 \frac{\Vf}{\vR}\xrightarrow{\PP}1.
\end{equation}
\end{theorem}
\begin{proof}
We apply Chebyshev's inequality to $J_H-\mu^\#\int_1^U K_T(t)^2\,dt$ and use Proposition~\ref{prop:concentration}. Since $\int_1^U K_T(t)^2\,dt\asymp T$, we obtain
\[
 J_H=(1+o_{\PP}(1))Z_zp^\#\int_1^U K_T(t)^2\,dt.
\]
By Lemma~\ref{lem:ballot}, we may replace $p^\#$ by $p_{z,T}$.

Let $\Omega_G$ be the event that $\mathcal G_T(t)$ holds for every $t\in[1,U]$, and put $I_+=\int_1^U|F_z(\tfrac12+it)|^2K_T(t)^2\,dt$. On $\Omega_G$, we have $I_+=J_G$. Apply Markov's inequality to $J_G-J_H$ and add $\PP(\Omega_G^c)$. For every fixed $\varepsilon>0$,
\[
 \begin{aligned}
 \PP\{I_+-J_H>\varepsilon TZ_zp_{z,T}\}
 &\le\PP(\Omega_G^c)
     +\frac{\EE(J_G-J_H)}
            {\varepsilon TZ_zp_{z,T}}\\
 &=o(1),
 \end{aligned}
\]
by Propositions~\ref{lem:weakbarrier} and~\ref{lem:barriergap}. Hence $I_+=J_H+o_{\PP}(TZ_zp_{z,T})$.

The reduction \eqref{eq:frequency-reduction}, together with $p_{z,T}\sim p_x$, now gives
\[
 \int_{\mathbb R}|F_z(\tfrac12+it)|^2K_T(t)^2\,dt
 =2J_H+o_{\PP}(TZ_zp_{z,T}).
\]
Since $2\int_1^U K_T^2\sim\int_{\mathbb R}K_T^2$ by Lemma~\ref{lem:kernel}, this proves \eqref{eq:CVintegral}.

Finally, Proposition~\ref{prop:transfer} and the exact kernel normalization in Lemma~\ref{lem:kernel} yield
\[
 \Vf
 =(1+o_{\PP}(1))d_zZ_z\,yp_{z,T}
 =(1+o_{\PP}(1))\frac6{\pi^2}yp_x.
\]
Here $d_zZ_z\to6/\pi^2$, $p_{z,T}/p_x\to1$, and the error $O(yM^{-50})$ in the transfer is $o(\vR)$ because $\vR\gg y/\sqrt M$. This proves \eqref{eq:CVgoal}.
\end{proof}

Together with the estimates in Section~\ref{sec:elementary}, this establishes \eqref{eq:smooth-target}-\eqref{eq:variance-target}. It remains to prove \eqref{eq:fluctuation-target}, which we do in Section~\ref{sec:conditionalCLT}.

\section{Conditional fourth moment}\label{sec:conditionalCLT}
We now prove \eqref{eq:fluctuation-target}, the last of the four estimates stated at the beginning of Section~\ref{sec:conditional}. Put $\cQ_z=\sum_{p>z}A_p^2=\sum_{p>z}\Delta_p^2$. Its conditional mean is $\Vf=\EE[\cQ_z\mid\cF_z]$. Define
\[
 \cE_z=\EE[(\cQ_z-\Vf)^2\mid\cF_z].
\]
Our goal is
\[
 \frac{\cE_z}{\vR^2}\xrightarrow{\PP}0.
\]
The proof proceeds as follows. Section~\ref{sec:sixedge} expresses the smoothed conditional error as a four frequency integral and gives two bounds for its arithmetic kernel. Section~\ref{sec:localization} uses these bounds to leave only a narrow two point integral. Section~\ref{sec:endpoint} controls that integral by endpoint barriers, uniformly including zero frequency. Section~\ref{sec:smoothing} removes the smoothing and the frequency truncation, completing \eqref{eq:fluctuation-target}. Finally, Section~\ref{sec:completion} combines this estimate with the results of Section~\ref{sec:conditional} to prove Theorem~\ref{thm:main}.

\subsection{Arithmetic reduction and the four frequency integral}\label{sec:sixedge}
We first expand $\cE_z$ using conditional orthogonality. The resulting square product relation is described by six edge variables and a common factor. After smoothing and Mellin inversion, it gives a four frequency integral. The two complementary kernel bounds at the end of this subsection will be used for localization in Section~\ref{sec:localization}.

Using \eqref{eq:roughdecomp}, a largest prime decomposition of the rough variables gives
\[
 \cQ_z=\sum_{\substack{r_1,r_2\in\mathcal A_z\\P(r_1)=P(r_2)}}
 a(r_1)a(r_2) f\left(\frac{r_1r_2}{(r_1,r_2)^2}\right).
\]
The diagonal is $\Vf$. Conditional orthogonality on primes greater than $z$ therefore gives
\begin{equation}\label{eq:Efour}
 \cE_z=
 \sum_{\substack{r_1,r_2,r_3,r_4\in\mathcal A_z\\
 r_1r_2r_3r_4=\square\\
 P(r_1)=P(r_2),\ P(r_3)=P(r_4)\\
 r_1\ne r_2,\ r_3\ne r_4}}
 a(r_1)a(r_2)a(r_3)a(r_4).
\end{equation}
Here $\square$ denotes a perfect square. The expression is nonnegative as a whole; individual products of coefficients need not be nonnegative. Absolute values are not inserted until after the Mellin representation.

\begin{lemma}\label{lem:sixedge}
Four square-free integers satisfy $r_1r_2r_3r_4=\square$ if and only if they have a unique representation
\begin{equation}\label{eq:edges}
 r_1=\rho gab,\quad r_2=\rho gcd,\quad
 r_3=\rho hac,\quad r_4=\rho hbd,
\end{equation}
where $\rho=\gcd(r_1,r_2,r_3,r_4)$ and
\[
 (g,h,a,b,c,d)=\frac1\rho\bigl(\gcd(r_1,r_2),\gcd(r_3,r_4),
 \gcd(r_1,r_3),\gcd(r_1,r_4),\gcd(r_2,r_3),\gcd(r_2,r_4)\bigr).
\]
The seven variables are square-free and pairwise coprime. If all prime factors of the $r_i$ exceed $z$, the same is true of the seven variables. When $\rho=1$, the maximum prime and off diagonal constraints in \eqref{eq:Efour} are equivalent to $g,h>1$, $abcd>1$, $P(abcd)<\min\{P(g),P(h)\}$.
\end{lemma}
\begin{proof}
Each prime divides either two or four of the $r_i$. Those dividing all four form $\rho$; each remaining prime occurs in exactly one of the six gcd quotients. This gives \eqref{eq:edges}, pairwise coprimality and uniqueness. Conversely, \eqref{eq:edges} gives $r_1r_2r_3r_4=\rho^4(ghabcd)^2$.

If $\rho=1$, the common maximum prime of $r_1,r_2$ must lie in $g$, and similarly the common maximum prime of $r_3,r_4$ must lie in $h$. This is equivalent to the last inequality and $g,h>1$. Under pairwise coprimality, $r_1=r_2$ occurs exactly when $a=b=c=d=1$, which is also equivalent to $r_3=r_4$. Thus the two off diagonal conditions become $abcd>1$.
\end{proof}

\paragraph{Smoothing and finite Mellin expansion.}
Put $\delta=\log(1+y/x)$, $\eta=(\log x)^{-1000}$ and $W_0=(\log x)^{2000}$. Following \cite[Section~4, (4.5)-(4.9)]{HSX}, smooth the interval multiplier to $\psi_\eta(u)=\frac1\eta\int_0^\eta \ind_{e^{-v}<u\le e^{\delta-v}}\dd v$. It is bounded by one and differs from the indicator only on two endpoint intervals of logarithmic length $\eta$. Mellin inversion gives
\begin{equation}\label{eq:smoothingkernel}
 \psi_\eta(u)=\frac1{2\pi}\int_{\RR}u^{-1/2-it}W(t)\dd t,
 \quad
 W(t)=\frac{e^{\delta s}-1}{s}\frac{1-e^{-\eta s}}{\eta s},
\end{equation}
where $s=\tfrac12+it$, and
\begin{equation}\label{eq:Wbound}
 |W(t)|\ll\min\{\delta,(1+|t|)^{-1},
                  \eta^{-1}(1+|t|)^{-2}\}.
\end{equation}
Thus the smoothed coefficient is
\[
 a_\eta(r)=\frac1{2\pi}\int_{\RR}
 F_z(\tfrac12+it)(x/r)^{1/2+it}W(t)\dd t.
\]
Let $\cE_{z,\eta}$ denote \eqref{eq:Efour} with $a(r)$ replaced by $a_\eta(r)$.

Before Mellin inversion, extend the rough variables to square-free integers supported on $(z,3x]$; the added terms have $a_\eta(r)=0$. This is a finite sum, so exchanging it with the four integrals gives
\begin{equation}\label{eq:Efourier}
 \cE_{z,\eta}\le Cx^2\int_{\RR^4}\prod_{j=1}^4|F_z(\tfrac12+it_j)W(t_j)|
 |G_R(\mathbf t)|\dd\mathbf t,
\end{equation}
where $G_R(\mathbf t)=\sum^{*}_{r_1,r_2,r_3,r_4} \prod_{j=1}^4r_j^{-1/2-it_j}$. The star imposes the square product, maximum prime and off diagonal constraints of \eqref{eq:Efour}, and restricts prime factors to $(z,3x]$.

\paragraph{Complementary Euler product bounds.}
Put $u=\frac{\log(3x)}{\log z}\asymp\log M$. For $v\ne0$, put $\gamma(v)=\min\{1,(|v|\log(3x))^{-1}\}$, and set $\gamma(0)=1$.

\begin{lemma}\label{lem:roughkernel}
For $|t_j|\le W_0$,
\begin{align}\label{eq:roughkernel}
 |G_R(\mathbf t)|\ll u^6\min\big\{&
 \gamma(t_1+t_2)\gamma(t_3+t_4),\notag\\[-2pt]
 &\gamma(t_1+t_3)+\gamma(t_1+t_4)
 +\gamma(t_2+t_3)+\gamma(t_2+t_4)\big\}+O(u^6/z).
\end{align}
Also $|G_R(\mathbf t)|\ll u^6$ for arbitrary real frequencies.
\end{lemma}
\begin{proof}
The parametrization \eqref{eq:edges} supplies the harmonic weight
$\rho^{-2}(ghabcd)^{-1}$. The contribution of $\rho>1$ is bounded by
\[
 \left(\prod_{z<p\le3x}(1+p^{-2})-1\right)
 \prod_{z<p\le3x}(1+p^{-1})^6\ll u^6/z.
\]
In the $\rho=1$ sum, remove pairwise coprimality of the six edge variables. The error is bounded by the total weight of an overlap prime, again $O(u^6\sum_{p>z}p^{-2})$.

For $z\le v\le3x$, set $E_v(\xi)=\prod_{z<p\le v}(1+p^{-1-i\xi})$. The partial summation argument in \cite[(4.13) and the proof of Lemma~4.4]{HSX} gives, uniformly for $|\xi|\le2W_0$,
\begin{equation}\label{eq:Eulerosc}
 |E_v(\xi)|\ll u,\qquad
 |E_v(\xi)-1|\ll u\gamma(\xi).
\end{equation}
Indeed, the square-free and reciprocal Euler products satisfy
\[
 E_v(\xi)=\bigl(1+O(z^{-1})\bigr)
          \prod_{z<p\le v}(1-p^{-1-i\xi})^{-1}.
\]
Here $\prod_{z<p\le v}(1-p^{-2-2i\xi})=1+O(z^{-1})$. Since $W_0$ is a fixed power of $\log x$, the prime number theorem error remains negligible, and the additional $O(u/z)$ is absorbed by $u\gamma(\xi)$.

Fix $a,b,c,d$ with $abcd>1$, and put $D=\max\{z,P(abcd)\}$. The maximum prime constraints in Lemma~\ref{lem:sixedge} require $P(g),P(h)>D$. After removing coprimality as above, the $g$ sum is
\[
 \sum_{\substack{g\in\mathcal A_z\\D<P(g)\le3x}}
       \frac1{g^{1+i(t_1+t_2)}}
 =E_{3x}(t_1+t_2)-E_D(t_1+t_2).
\]
The $h$ sum has the same form with $t_3+t_4$ in place of $t_1+t_2$. Subtracting $1$ from both Euler products and using \eqref{eq:Eulerosc} bounds these two differences by $O(u\gamma(t_1+t_2))$ and $O(u\gamma(t_3+t_4))$, respectively. The four remaining harmonic sums contribute $O(u^4)$. This proves the first kernel bound.

Alternatively fix $g,h$ and put $P_*=\min\{P(g),P(h)\}$. The strict maximum prime inequality in Lemma~\ref{lem:sixedge} means that every cross edge is supported on primes $z<p<P_*$. Accordingly set $E_{P_*^-}(\xi)=\prod_{z<p<P_*}(1+p^{-1-i\xi})$. The cross edge sum is
\[
 \begin{aligned}
 &E_{P_*^-}(t_1+t_3)E_{P_*^-}(t_1+t_4)\\
 &\qquad\times E_{P_*^-}(t_2+t_3)E_{P_*^-}(t_2+t_4)-1.
 \end{aligned}
\]
The subtraction of $1$ removes $a=b=c=d=1$. Write this difference as a telescoping sum and apply \eqref{eq:Eulerosc}; the resulting bound is
\[
 Cu^4\bigl(\gamma(t_1+t_3)+\gamma(t_1+t_4)
          +\gamma(t_2+t_3)+\gamma(t_2+t_4)\bigr).
\]
Summing $g,h$ absolutely contributes $u^2$. Combining the two estimates proves \eqref{eq:roughkernel}.
\end{proof}

\subsection{Reduction to a narrow two point integral}\label{sec:localization}
We use the complementary kernel bounds to reduce \eqref{eq:Efourier} to the narrow two point integral in \eqref{eq:Icore}. We first need a half moment estimate for the discarded frequency configurations. The cutoffs below are independent of the scale $T$ used for conditional variance concentration.

\paragraph{A uniform two point half moment.}
For $t,s\ge0$ write
$S_z(v)=\sum_{p\le z}p^{-1}\cos(v\log p)$. Expanding a single Euler factor and taking its expectation gives, uniformly in $t,s$,
\begin{align}\label{eq:halfmomentexpansion}
 \log\EE|F_z(\tfrac12+it)F_z(\tfrac12+is)|
 ={}&\tfrac12S_z(0)+\tfrac12S_z(t-s)+\tfrac12S_z(t+s)\notag\\
 &-\tfrac14S_z(2t)-\tfrac14S_z(2s)+O(1).
\end{align}
The error is summable over primes. This formula retains the negative reflection corrections; discarding them near the origin is not permissible.

\begin{lemma}\label{lem:halfmoment}
For $0\le t,s\le W_0$, one has
\begin{equation}\label{eq:halfmoment}
 \EE|F_z(\tfrac12+it)F_z(\tfrac12+is)|
 \ll \Lambda^{1/2}M^{20}
 \left(1+\min\{\Lambda,|t-s|^{-1}\}^{1/2}\right),
\end{equation}
with the minimum interpreted as $\Lambda$ at equality.
\end{lemma}
\begin{proof}
The comparison $S_z(v)=\log|\zeta(1+\Lambda^{-1}+iv)|+O(1)$ follows by bounding the difference between the sharp and exponentially damped prime sums and the contribution of higher prime powers. For $|v|\le1$, the pole and the absence of zeros on the adjacent compact region give $S_z(v)=\log\min\{\Lambda,1/|v|\}+O(1)$. For $1/2\le|v|\le2W_0$, split the prime sum at $v_0=\exp(CM^2)$ with $C$ sufficiently large. The primes at most $v_0$ contribute at most $2\log M+O(1)$ in absolute value by Mertens' theorem. For $v_0<p\le z$, partial summation as in Lemma~\ref{lem:primeosc} gives an oscillatory integral of size $O((|v|\log v_0)^{-1})$ and an error $O((1+|v|)e^{-c\sqrt{\log v_0}})$. Both are $O(1)$ in this range.

Put $D(v)=\max\{|v|,\Lambda^{-1}\}$. When $t+s\le1/2$, we have $D(2t),D(2s)\le2D(t+s)$, and hence
\[
 S_z(t+s)-\tfrac12S_z(2t)-\tfrac12S_z(2s)
 =\log\frac{\sqrt{D(2t)D(2s)}}{D(t+s)}+O(1)\le C.
\]
When $t+s>1/2$, the first term is at most $2\log M+O(1)$ and each of $S_z(2t),S_z(2s)$ is at least $-2\log M-O(1)$, using the small frequency estimate when necessary. Thus the same expression is bounded above by $4\log M+O(1)$. Inserting these upper bounds and the estimates for $S_z(0)$ and $S_z(t-s)$ in \eqref{eq:halfmomentexpansion} proves the lemma.
\end{proof}

\paragraph{Folding and localization.}
Because $|F_z(\tfrac12-it)|=|F_z(\tfrac12+it)|$ and $|W(-t)|=|W(t)|$, $\gamma(t_i+t_j)\le\gamma(|t_i|-|t_j|)$. Thus the four dimensional integral in \eqref{eq:Efourier} may be bounded by a constant times the corresponding positive quadrant integral with every edge frequency a difference. This is a deterministic inequality. No probabilistic independence of the positive and negative axes is involved.

Put $\delta_0=y/x$ and $w(t)=\min\{\delta_0,(1+t)^{-1}\}$ for $t\ge0$. Choose $r_1=M^{200}/\Lambda$ and $r_2=M^{1000}/\Lambda$, so that $r_1=o(r_2)$ and $r_2\to0$. The first kernel bound removes configurations with a large distance inside either pair $(t_1,t_2)$ or $(t_3,t_4)$. When both pairs are close internally, the second bound removes configurations in which the pairs are far apart. All four remaining frequencies are therefore close, and integrating out two of them gives the following bound.

\begin{proposition}\label{prop:localization}
There is a fixed constant $c_1>0$ such that, outside an event of probability $O(M^{-c_1})$, the smoothed conditional error is bounded by
\begin{equation}\label{eq:localization}
 \cE_{z,\eta}\ll \vR^2M^{-10}+\mathcal I,
\end{equation}
up to frequency tail and common prime errors whose normalized expectations tend to zero faster than every fixed power of $M^{-1}$, where
\begin{equation}\label{eq:Icore}
 \mathcal I=\frac{x^2u^6}{\Lambda^2}
 \iint_{\substack{0\le t,s\le W_0\\|t-s|\le Cr_2}}
 |F_z(\tfrac12+it)F_z(\tfrac12+is)|^2w(t)^2w(s)^2\dd t\dd s.
\end{equation}
\end{proposition}
\begin{proof}
We divide the frequency configurations as in the proof of \cite[Lemma~4.5]{HSX}, using the Rademacher half moment bound above. All double integrals below are over $[0,W_0]^2$.

First suppose that either $|t_1-t_2|$ or $|t_3-t_4|$ exceeds $r_1$. The first member of the minimum in \eqref{eq:roughkernel} bounds this contribution by $Cx^2u^6J_0J_1$, where
\begin{align*}
 J_0&=\iint |F_z(\tfrac12+it)F_z(\tfrac12+is)|w(t)w(s)\gamma(t-s)\dd t\dd s,\\
 J_1&=\iint_{|t-s|>r_1}|F_z(\tfrac12+it)F_z(\tfrac12+is)|w(t)w(s)\gamma(t-s)\dd t\dd s.
\end{align*}
We use \cref{lem:halfmoment} to estimate the expected Euler product in $J_0$ and $J_1$, and apply $\int_0^\infty w(t)w(t+h)\dd t\ll\delta_0$ to integrate over one frequency. The resulting bounds are
\begin{equation}\label{eq:Jexpectations}
 \EE J_0\ll\delta_0M^{21},\qquad
 \EE J_1\ll\delta_0M^{-79}.
\end{equation}
For $J_0$, split the difference $h$ at $\Lambda^{-1}$ and $1$. For $J_1$, the singular part over $r_1<h<1$ is bounded by
\[
 \frac{\delta_0M^{20}}{u\sqrt\Lambda}
 \int_{r_1}^{1}h^{-3/2}\dd h
 \ll\delta_0M^{-80}.
\]
The part with $h>1$ is smaller: its logarithmic length is $O(M)$ and it has the factor $\Lambda^{-1/2}$. We apply Markov's inequality separately to $J_0$ and $J_1$, using the expectations in \eqref{eq:Jexpectations}. Outside an event of probability $O(M^{-8})$, we obtain
\begin{equation}\label{eq:Jgood}
 J_0\le\delta_0M^{30},\qquad
 J_1\le\delta_0M^{-60}.
\end{equation}
This part of the four frequency integral is therefore $\ll y^2u^6M^{-30}$. Since $p_x\gg M^{-1/2}$, we have $\vR^2\gg y^2/M$. Together with $u\asymp\log M$, this gives a contribution $\ll\vR^2M^{-10}$.

Next suppose that both paired distances are at most $r_1$, but the two pairs are separated by more than a constant multiple of $r_2$. All cross distances then exceed a constant multiple of $r_2$, so the second member of the minimum in \eqref{eq:roughkernel} bounds the kernel by $Cu^5M^{-1000}$. Define
\[
 J_2=\iint_{|t-s|\le r_1}|F_z(\tfrac12+it)F_z(\tfrac12+is)|w(t)w(s)\dd t\dd s.
\]
For $J_2$, we use \cref{lem:halfmoment} and split the integral over $|t-s|$ at $\Lambda^{-1}$. We obtain
\[
 \EE J_2\ll\delta_0M^{120}.
\]
Applying Markov's inequality to $J_2$ at level $\delta_0M^{140}$, we obtain $J_2\le\delta_0M^{140}$ outside an event of probability $O(M^{-20})$. Outside this exceptional event, the contribution of these configurations is at most
\[
 Cx^2u^5M^{-1000}J_2^2\ll y^2u^5M^{-720},
\]
which is again $\ll\vR^2M^{-10}$.

Finally, all four folded frequencies lie within $Cr_2$ of each other. Write $F_j=F_z(\tfrac12+it_j)$ and use the first member of the minimum in \eqref{eq:roughkernel}, together with
\[
 |F_1F_2F_3F_4|\le\tfrac12\bigl(|F_1F_3|^2+|F_2F_4|^2\bigr).
\]
The two resulting terms have the same bound after relabeling. Consider the term with $|F_1F_3|^2$. Since $|t_1-t_2|,|t_3-t_4|\le Cr_2=o(1)$, we have $w(t_2)\ll w(t_1)$ and $w(t_4)\ll w(t_3)$. We keep $t_1,t_3$ and integrate out $t_2,t_4$. Enlarging the inner domains, this contribution is at most
\[
 \begin{aligned}
 &Cx^2u^6\iint_{\substack{0\le t_1,t_3\le W_0\\|t_1-t_3|\le Cr_2}}
 |F_1F_3|^2w(t_1)^2w(t_3)^2\\
 &\qquad\times
 \left(\int_{|t_2-t_1|\le Cr_2}\gamma(t_1-t_2)\,dt_2\right)
 \left(\int_{|t_4-t_3|\le Cr_2}\gamma(t_3-t_4)\,dt_4\right)
 \,dt_1\,dt_3.
 \end{aligned}
\]
For either inner integral, the change of variables to the difference gives
\[
 \int_{|h|\le Cr_2}\gamma(h)\,dh
 \ll\frac{1+\log(r_2\log(3x))}{\log(3x)}
 \ll\frac1\Lambda.
\]
Here $\log(3x)=u\Lambda$, $r_2\log(3x)=uM^{1000}$ and $u\asymp\log M$. The two integrations therefore supply the factor $\Lambda^{-2}$. Relabeling $t_1,t_3$ as $t,s$ gives the integral in \eqref{eq:Icore}, up to an absolute constant.

For the kernel error $O(u^6/z)$, apply H\"older's inequality to the product of the four Euler factors, using $\EE|F_z(\tfrac12+it)|^4\ll\Lambda^6$ for each frequency, and then integrate the four weights $|W(t_j)|$. The expected total contribution, including shared primes and the removal of edge coprimality, is $\ll x^2u^6z^{-1}\Lambda^6(\log(1/\eta))^4$. Division by $\vR^2\gg x^2/(M\log x)$ makes this negligible. The frequency tails are treated in \cref{sec:smoothing}.
\end{proof}

\subsection{Endpoint barriers and the narrow integral}\label{sec:endpoint}
We now estimate the narrow integral $\mathcal I$ left by Proposition~\ref{prop:localization}. We shall prove that, for some absolute constant $c>0$,
\begin{equation}\label{eq:narrow-goal}
 \PP(\mathcal I>\vR^2M^{-c})\ll M^{-c}.
\end{equation}
We use a weak path restriction that holds throughout the frequency range, together with a stronger restriction at the last time. The contribution lost at the last time will be bounded by a weighted strip estimate. On the restricted event, the terminal upper bound controls one Euler product directly.

\paragraph{The path and terminal restrictions.}
For this estimate we use a different prime partition from Section~\ref{sec:barriers}. Put $J=\lfloor M-10\log M\rfloor$, $b_0=M-J$ and $w_j=\exp(e^{b_0+j})$ for $0\le j\le J$. Thus $\log w_0\asymp M^{10}$, $w_J=z$ and $J\sim M$. Set $R_j(t)=\prod_{w_0<p\le w_j}(1+\eps_pp^{-1/2-it})$ and $S_j(t)=\log|R_j(t)|$. The primes at most $w_0$ form an independent initial factor. They account for only $O(\log M)$ of the $M$ logarithmic time units, while the starting point makes the block Berry-Esseen errors smaller than $e^{-J}$, uniformly over the required frequencies. We verify this below. The two prime partitions need not be nested; their high probability events will be intersected at the end.

For $a\ge1$, define the weak event by
\begin{equation}\label{eq:endbarriers}
 \mathcal G_t(a)=\{S_j(t)\le j+a+30\log M\text{ for }1\le j\le J\}.
\end{equation}
The stronger event is $\mathcal H_t(a)=\mathcal G_t(a)\cap\{S_J(t)\le J+a-D_0\log M\}$, where $D_0=10000$. Thus $\mathcal G_t(a)\setminus\mathcal H_t(a)$ consists of weakly restricted paths whose final value lies in a strip of width $(D_0+30)\log M$. We need the following two estimates: one for the weighted mean under the weak restriction, and one for the weighted loss in this terminal strip.

\begingroup
\setlength{\abovedisplayskip}{8pt plus 2pt minus 2pt}
\setlength{\belowdisplayskip}{8pt plus 2pt minus 2pt}
\begin{lemma}\label{lem:endpoint}
For each fixed $C>0$, uniformly for $0\le t\le e^{CM}$ and $a\ge1$, we have
\begin{equation}\label{eq:endpointmean}
 \EE[|F_z(\tfrac12+it)|^2\ind_{\mathcal G_t(a)}]
 \ll_C Z_z\min\{1,(a+\log M)/\sqrt M\}.
\end{equation}
The contribution from the endpoint strip satisfies
\begin{equation}\label{eq:endpointstrip}
 \EE[|F_z(\tfrac12+it)|^2\ind_{\mathcal G_t(a)\setminus\mathcal H_t(a)}]
 \ll_C Z_z\min\{1,a/\sqrt M\}M^{-1/5}.
\end{equation}
\end{lemma}
\endgroup

The proof of this lemma is given after Proposition~\ref{prop:transitionstrip}.

\paragraph{A probability estimate for the tilted walk.}
To prove these estimates, we absorb the terminal weight $|R_J(t)|^2$ into the probability measure. The tilted Rademacher increments have an additional drift near zero frequency, and their variances change at the same scale. For $t>0$, put $k_t=\min\{J,\max(0,\log(1/(2t\log w_0)))\}$, and set $k_t=J$ at $t=0$. This is the time at which $2t\log w_j$ reaches one, up to the truncation at $0$ and $J$. Before this time the cosine terms in the tilted moments have not begun to oscillate; after it their cumulative contribution is bounded.

As verified in the proof of Lemma~\ref{lem:endpoint}, the cumulative tilted mean is $j+\frac12\min(j,k_t)+O(1)$ and the cumulative variance is $\frac12(j+\min(j,k_t))+O(1)$. The next proposition allows this change of drift and variance. Its first conclusion bounds the weak event, and its second bounds the terminal strip. It also applies to the character twist in Section~\ref{sec:twists}, where the additional drift is absent and $k=0$.

\begin{proposition}\label{prop:transitionstrip}
Fix positive constants $c_*,C_*,K,C_0,D_0$. Let $Y_1,\ldots,Y_J$ be independent real random variables, where $J\ge3$, and write $S_j=\sum_{i\le j}Y_i$, $\mu_j=\EE S_j$ and $V_j=\Var(S_j)$. Suppose $c_*\le\Var(Y_j)\le C_*$ and, for some $k\in[0,J]$, we have
\[
 \left|\mu_j-j-\tfrac12\min(j,k)\right|\le K
 \qquad(1\le j\le J).
\]
Suppose also that
\[
 \left|V_j-\tfrac12\bigl(j+\min(j,k)\bigr)\right|\le K
 \qquad(1\le j\le J),
\]
and that the standardized increments satisfy
\[
 \max_{j\le J}\sup_{u\in\RR}
 \left|\PP\left(\frac{Y_j-\EE Y_j}{\sqrt{\Var(Y_j)}}\le u\right)-\Phi(u)\right|
 \le e^{-J},
\]
where $\Phi$ is the standard Gaussian distribution function. For $a\ge1$, define $\mathcal A(a)=\{S_j\le j+a+C_0\log J\text{ for }1\le j\le J\}$. Then, uniformly in $k$ and $a$,
\[
 \PP(\mathcal A(a))
 \ll\min\{1,(a+\log J)/\sqrt J\}.
\]
The probability of ending in a strip next to the upper boundary satisfies
\[
 \PP\bigl(\mathcal A(a),\ S_J>J+a-D_0\log J\bigr)
 \ll\min\{1,a/\sqrt J\}J^{-1/5}.
\]
The implied constants depend only on $c_*,C_*,K,C_0,D_0$.
\end{proposition}
\begin{proof}
We first compare the discrete path with Brownian motion. For the standard Brownian motion facts used below, see \cite[Chapter~2]{LeGallBM}. Put $\sigma_j^2=\Var(Y_j)$ and couple $Y_j$ with $\EE Y_j+\sigma_jZ_j$ by independent uniform quantiles, where the $Z_j$ are independent standard Gaussians. For fixed $D>0$, all $|Z_j|\le C_D\sqrt{\log J}$ outside an event of probability $O_D(J^{-D})$. The Gaussian density there is bounded below by a negative power of $J$, so the Kolmogorov bound gives
\[
 \sum_{j\le J}|Y_j-\EE Y_j-\sigma_jZ_j|
 \ll_D J^{C_D}e^{-J}=o(1).
\]

Extend these Gaussian increments to a Brownian motion $B$ so that $B_{V_j}=\sum_{i\le j}\sigma_iZ_i$. Set $\nu_j=(j+\min(j,k))/2$. Since $|V_j-\nu_j|\le K$, Gaussian tails and a union bound control $B_{V_j}-B_{\nu_j}$ by $C_D\sqrt{\log J}$ for all $j$, outside another event of probability $O_D(J^{-D})$. Including the bounded error in the means gives
\begin{equation}\label{eq:pathcoupling}
 \PP\left(\max_{j\le J}
 \left|S_j-j-B_{\nu_j}-\tfrac12\min(j,k)\right|
 >C_D\sqrt{\log J}\right)\ll_D J^{-D}.
\end{equation}

To pass from the discrete checks to continuous survival, consider $X(r)=B_r+\tfrac12\min(r,k)$ for $0\le r\le R=(J+k)/2$. The identity $\min(\nu_j,k)=\min(j,k)$ shows that $X(\nu_j)$ is the reference path in \eqref{eq:pathcoupling}. The time spacings $\nu_j-\nu_{j-1}$ lie in $[1/2,1]$. Brownian bridge bounds on these intervals, followed by a union bound, show that every bridge differs from its linear interpolation by at most $C_D\sqrt{\log J}$, with failure probability $O_D(J^{-D})$. The deterministic drift is Lipschitz, so its interpolation error is bounded. Consequently $\mathcal A(a)$ implies continuous survival below
\[
 U=a+C_0\log J+C_D\sqrt{\log J},
\]
apart from these exceptional events. On the strip event, $0\le U-X(R)\le D_J$, where $D_J=(C_0+D_0)\log J+2C_D\sqrt{\log J}$.

The drift of $X$ is nonnegative. Thus survival of $X$ below $U$ implies survival of $B$ below $U$. Since $J/2\le R\le J$, the reflection principle gives $\PP(\mathcal A(a))\ll\min\{1,(a+\log J)/\sqrt J\}$ after absorbing the exceptional probability.

For the strip bound, first suppose $a\ge J^{1/4}$. The terminal variable $X(R)$ is Gaussian with variance $R\asymp J$, so its probability of lying in an interval of length $D_J$ is $O(\log J/\sqrt J)$. This is $O(\min\{1,a/\sqrt J\}J^{-1/5})$.

It remains to consider $1\le a<J^{1/4}$. If $k\ge J/4$, survival at time $J/4$ forces $B_{J/4}\le U-J/8$. Since $U=o(J)$ in this range of $a$, this has probability $O(e^{-cJ})$. If $k<J/4$, put $Z(r)=U-X(r)$ and kill this process on reaching zero. Up to time $k$ it is a nonnegative killed supermartingale, and hence
\[
 \EE[Z(k);\text{survival up to }k]\le U.
\]
After time $k$ the drift vanishes, and the remaining time is $R-k=(J-k)/2\ge3J/8$. By the reflection principle, the killed Brownian transition density is
\begin{equation}\label{eq:killed}
 p_r^+(u,v)=\frac1{\sqrt{2\pi r}}
 \left(e^{-(u-v)^2/(2r)}-e^{-(u+v)^2/(2r)}\right)
 \ll\frac{uv}{r^{3/2}}.
\end{equation}
We integrate the killed transition density over $0<v<D_J$ and then average over $Z(k)$ on the survival event. The bound on $\EE[Z(k);\text{survival up to }k]$ gives
\[
 \begin{aligned}
 &\PP\bigl(X(r)\le U\ (0\le r\le R),\ 0\le U-X(R)\le D_J\bigr)\\
 &\qquad\le C\frac{(a+\log J)(\log J)^2}{J^{3/2}}
 \ll\min\{1,a/\sqrt J\}J^{-1/5}.
 \end{aligned}
\]
Choosing $D$ large absorbs all coupling and bridge errors. The bounded range of smaller $J$ is covered by increasing the implied constants.
\end{proof}

\paragraph{Verification for the Euler product.}
\begin{proof}[Proof of Lemma~\ref{lem:endpoint}]
Tilt the primes above $w_0$ by $|R_J(t)|^2$, divided by its expectation. Write $\alpha_j=\sum_{w_{j-1}<p\le w_j}p^{-1}$. The frequency dependent term is
\[
 \beta_j(t)=\sum_{w_{j-1}<p\le w_j}\frac{\cos(2t\log p)}p.
\]
By Lemma~\ref{lem:tilt}, the mean is $\alpha_j+\beta_j/2+O(w_{j-1}^{-1})$ and the variance is $(\alpha_j+\beta_j)/2+O(w_{j-1}^{-1})$ for the $j$th logarithmic increment.

Writing $B_j=\log w_{j-1}$, the variance is $\int_1^e\cos^2(tB_jv)\,dv/v+o(1)$. This continuous function of $tB_j$ is positive at every finite argument, equals one at zero, and tends to $1/2$ at infinity. Its infimum is therefore positive. The block variances lie in a fixed positive compact interval, uniformly in $t$. The prime number theorem error is uniform on the stated range because $e^{CM}e^{-cM^5}=o(1)$.

For the transition time $k_t$ defined above, put $\tau_j(k)=\min(j,k)-\min(j-1,k)$. We claim that
\begin{equation}\label{eq:betasum}
 \sum_{j=1}^J|\beta_j(t)-\tau_j(k_t)|\ll_C1.
\end{equation}
Before $2t\log w_j$ reaches one, expand the cosine and sum the squared frequency errors geometrically. After $2t\log w_{j-1}$ exceeds one, integration by parts gives a summable error $O((2t\log w_{j-1})^{-1})$. Only a bounded number of blocks meet the transition. At $t=0$ use $\beta_j=\alpha_j$. Since $\sum_{i\le j}\alpha_i=j+o(1)$ uniformly in $j$, the cumulative tilted mean satisfies
\begin{equation}\label{eq:reflectiontime}
 \mu_j=j+\tfrac12\min(j,k_t)+O_C(1).
\end{equation}
The corresponding cumulative variance satisfies
\[
 V_j=\tfrac12\bigl(j+\min(j,k_t)\bigr)+O_C(1).
\]

The prime contributions remain independent under this product tilt. Their third centered absolute moments are $O(p^{-3/2})$. Lemma~\ref{lem:BE} with $d=1$, applied using the variance lower bound, bounds the standardized block Kolmogorov errors by
\[
 \delta_j(t)\ll\sum_{p>w_{j-1}}p^{-3/2}
 \ll w_{j-1}^{-1/2}
 \ll\exp(-cM^{10}e^{j-1}).
\]
In particular these errors are at most $e^{-J}$ for sufficiently large $M$. Proposition~\ref{prop:transitionstrip} applies with $k=k_t$. Since $J\sim M$, replacing the fixed multiples of $\log M$ in the two boundaries by slightly larger fixed multiples of $\log J$ only enlarges the events. The two conclusions of that proposition are precisely the required tilted probability bounds. The high prime tilt normalization and the exact mean of the independent small prime product multiply to $Z_z$, proving \eqref{eq:endpointmean} and \eqref{eq:endpointstrip}.
\end{proof}

\paragraph{A simultaneous upper barrier event.}
For applications put $T_0=x/y$, $a_t=\tfrac12\log(T_0+t)$, $\mathcal G^*(t)=\mathcal G_t(a_t)$, $\mathcal H^*(t)=\mathcal H_t(a_t)$. For large $x$, $a_t\ge1$. The derivative/Sobolev argument in the proof of \cref{lem:weakbarrier}, now on $[0,W_0]$, gives
\begin{equation}\label{eq:globalendbarrier}
 \PP(\exists t\in[0,W_0]:\mathcal G^*(t)\text{ fails})\ll M^{-20}.
\end{equation}
To check the powers, the expected supremum of $|R_j|^2$ on a mesh interval of length $(\log w_j)^{-1}$ is $O(e^j)$, whereas the threshold is
$e^{2j}(T_0+t)M^{60}$. The weighted number of intervals is bounded by
$\log w_j\int_0^{W_0}(T_0+t)^{-1}\dd t+O(1)$.
Since $\log w_j\asymp M^{10}e^j$, $J\asymp M$, and $\log(1+W_0/T_0)\ll M$, the total probability is $O(M^{-48})$, which implies the stated weaker bound. This uniform event uses no frequencywise Gaussian coupling.

\paragraph{Controlling the narrow integral.}
On the event \eqref{eq:globalendbarrier}, every frequency in the domain of $\mathcal I$ satisfies the weak barrier. We split the integral according to whether both points satisfy the endpoint restriction. The part where at least one point fails it will be bounded by a product of two one dimensional integrals. The remaining part will use the terminal upper bound.

Let $s_t=\min\{1,\log(T_0+t)/\sqrt M\}$, $s_0=\min\{1,\log T_0/\sqrt M\}\asymp p_x$. Set
\begin{align*}
 I_1&=\Lambda^{-1}\int_0^{W_0}|F_z(\tfrac12+it)|^2
                   \ind_{\mathcal G^*(t)}w(t)^2\dd t,\\
 I_2&=\Lambda^{-1}\int_0^{W_0}|F_z(\tfrac12+it)|^2
          (\ind_{\mathcal G^*(t)}-\ind_{\mathcal H^*(t)})w(t)^2\dd t.
\end{align*}
Here $I_1$ is the weakly restricted integral, and $I_2$ is the contribution removed by the endpoint restriction. To estimate the two expectations, integrate \eqref{eq:endpointmean} and \eqref{eq:endpointstrip} against $\Lambda^{-1}w(t)^2\dd t$. Use the elementary estimate
\[
 \int_0^{W_0}w(t)^2s_t\dd t\ll\delta_0s_0,
\]
and its version with $\log(T_0+t)+\log M$ in place of $\log(T_0+t)$ in the definition of $s_t$. The latter introduces a factor of at most $\log M$, so
\[
 \EE I_1\ll\delta_0s_0\log M,\qquad
 \EE I_2\ll\delta_0s_0M^{-1/5}.
\]
Both estimates retain the factor $s_0\asymp p_x$, which is needed for the scale $y^2s_0^2\asymp\vR^2$. We apply Markov's inequality to $I_1$ and $I_2$ at levels $\delta_0s_0M^{1/40}$ and $\delta_0s_0M^{-1/10}$, respectively. The sum of the two exceptional probabilities is $O(M^{-1/40}\log M+M^{-1/10})=O(M^{-1/50})$. Hence, outside this exceptional event,
\begin{equation}\label{eq:I1I2}
 I_1\le\delta_0s_0M^{1/40},\qquad
 I_2\le\delta_0s_0M^{-1/10}.
\end{equation}
On the event \eqref{eq:globalendbarrier}, both points satisfy the weak barrier. The contribution to \eqref{eq:Icore} where at least one point fails the endpoint restriction is therefore
\begin{equation}\label{eq:stripcontribution}
 \begin{aligned}
 &\frac{x^2u^6}{\Lambda^2}
 \iint_{\substack{0\le t,s\le W_0\\|t-s|\le Cr_2}}
 |F_z(\tfrac12+it)F_z(\tfrac12+is)|^2\\
 &\qquad\times
 \bigl(\ind_{\mathcal G^*(t)\cap\mathcal G^*(s)}
       -\ind_{\mathcal H^*(t)\cap\mathcal H^*(s)}\bigr)
 w(t)^2w(s)^2\dd t\dd s\\
 &\quad\ll x^2u^6I_1I_2
 \ll y^2s_0^2u^6M^{-3/40}\ll\vR^2M^{-1/20}.
 \end{aligned}
\end{equation}

For the remaining part, factor $F_z=F_{w_0}R_J$ and apply the terminal constraint to $R_J(s)$ on $\mathcal H^*(s)$. For the mixed moment below, independence separates the primes above $w_0$ from the two small prime factors. Thus
\[
 \begin{aligned}
 |F_z(\tfrac12+is)|^2\ind_{\mathcal H^*(s)}
 &\le |F_{w_0}(\tfrac12+is)|^2(T_0+s)e^{2J}M^{-2D_0},\\
 \EE[|F_z(\tfrac12+it)|^2|F_{w_0}(\tfrac12+is)|^2]
 &\ll\frac{\Lambda}{\log w_0}(\log w_0)^6
 \ll\Lambda M^{50}.
 \end{aligned}
\]
The second bound uses the small prime fourth moment. Since $e^{2J}\asymp\Lambda^2/M^{20}$, these inequalities imply
\begin{equation}\label{eq:endstrongmean}
 \EE[|F_z(\tfrac12+it)F_z(\tfrac12+is)|^2\ind_{\mathcal H^*(s)}]
 \ll (T_0+s)\Lambda^3M^{30-2D_0}.
\end{equation}
For the part where both points satisfy the endpoint restriction, apply \eqref{eq:endstrongmean}, integrate over $|t-s|\le Cr_2$, and then use $\int_0^\infty(T_0+t)w(t)^4\dd t\ll\delta_0^2$. Since $r_2=M^{1000}/\Lambda$ and $D_0=10000$, this gives
\begin{equation}\label{eq:strongcontribution}
 \begin{aligned}
 &\EE\Biggl[
 \frac{x^2u^6}{\Lambda^2}
 \iint_{\substack{0\le t,s\le W_0\\|t-s|\le Cr_2}}
 |F_z(\tfrac12+it)F_z(\tfrac12+is)|^2\\
 &\qquad\times\ind_{\mathcal H^*(t)\cap\mathcal H^*(s)}
 w(t)^2w(s)^2\dd t\dd s\Biggr]\\
 &\quad\ll y^2u^6M^{1030-2D_0}
 \ll\vR^2M^{-1000}.
 \end{aligned}
\end{equation}
Apply Markov's inequality to this nonnegative integral at level $\vR^2M^{-100}$. By \eqref{eq:strongcontribution}, the exceptional probability is $O(M^{100}M^{-1000})=O(M^{-900})$.

The strip estimate \eqref{eq:stripcontribution} and the terminal bound \eqref{eq:strongcontribution} prove \eqref{eq:narrow-goal}, absorbing $u^6\asymp(\log M)^6$ into fixed powers of $M$. Together with Proposition~\ref{prop:localization}, this gives a fixed power saving for the smoothed conditional error.

\subsection{Removal of smoothing and completion of the conditional estimate}\label{sec:smoothing}
We now return to the sharp interval and restore the discarded frequencies. We first prove an $L^2$ bound for the smoothing error, then control the Fourier tail and combine the estimates in Proposition~\ref{prop:QV}.

Let $A_{p,\eta}$ be obtained from $A_p$ by the multiplier $\psi_\eta$, and define
$\cQ_{z,\eta}=\sum_{p>z}A_{p,\eta}^2$,
$V_{z,\eta}=\EE[\cQ_{z,\eta}\mid\cF_z]$.
For each $p$, let $E_p$ count the integers where the two coefficient functions differ. Then $E_p$ counts two endpoint intervals of total length $O(\eta x/p)$, with at most a bounded number of endpoint terms added.

\begin{lemma}\label{lem:smoothing}
\begin{equation}\label{eq:smoothingerror}
 \|\cQ_z-\cQ_{z,\eta}\|_2
 \ll (x\eta^{1/4}+\eta^{-2})(\log(2x))^5.
\end{equation}
\end{lemma}
\begin{proof}
Write $\cQ_z-\cQ_{z,\eta}$ as $\sum_{p>z}(A_p-A_{p,\eta})(A_p+A_{p,\eta})$. We apply the triangle inequality in $L^2$ to this sum and H\"older's inequality to each product:
\[
 \|\cQ_z-\cQ_{z,\eta}\|_2
 \le\sum_{p>z}\|A_p-A_{p,\eta}\|_4\|A_p+A_{p,\eta}\|_4.
\]
Apply \cref{lem:fourth} to the two fourth norms in each summand. This bounds the smoothing error by
\[
 \|\cQ_z-\cQ_{z,\eta}\|_2
 \ll(\log(2x))^4\sum_{p>z}(x/p)^{3/4}E_p^{1/4}.
\]
For $p\le\eta x$, $E_p\ll\eta x/p$, and this contributes
$O(x\eta^{1/4}\log(2x))$ before the outside logarithm.
For $p>\eta x$, the difference can occur only when, for some integer $n\ll\eta^{-1}$, $p$ lies within $O(\eta x/n)$ of $x/n$ or $(x+y)/n$. Discarding primality bounds the contribution by
\[
 \sum_{n\ll\eta^{-1}}n^{3/4}(\eta x/n+1)
 \ll x\eta^{1/4}+\eta^{-7/4}
 \ll x\eta^{1/4}+\eta^{-2}.
\]
This proves \eqref{eq:smoothingerror}.

\end{proof}

The Fourier tail is similarly harmless but must be recorded. By \eqref{eq:Wbound},
\[
 \int_{|t|>W_0}|W(t)|\dd t\ll \eta^{-1}/W_0=(\log x)^{-1000},
 \qquad \int_{\RR}|W(t)|\dd t\ll\log\log x.
\]
The crude bound $|G_R|\ll u^6$ and
$\sup_t\EE|F_z(\tfrac12+it)|^4\ll\Lambda^6$ bound the expected part of \eqref{eq:Efourier} with at least one large frequency by
\[
 Cx^2u^6\Lambda^6(\log\log x)^3(\log x)^{-1000}=o(\vR^2M^{-B})
\]
for every fixed $B$. Hence truncating and restoring frequencies is justified in probability.

We can now state the conditional fluctuation estimate for the original interval.

\begin{proposition}\label{prop:QV}
There exists an absolute constant $c>0$ such that
\begin{equation}\label{eq:QVgoal}
 \mathbb P(\cE_z>\vR^2M^{-c})\ll M^{-c}.
\end{equation}
In particular $\cE_z/\vR^2\to0$ in probability.
\end{proposition}
\begin{proof}
Let $\cE_{z,\eta}$ denote the conditional quadratic variation error formed with the smoothed interval multiplier of Section~\ref{sec:sixedge}. By \cref{prop:localization}, outside an event of probability $O(M^{-c_1})$, $\cE_{z,\eta}\ll \vR^2M^{-10}+\mathcal I$, where $\mathcal I$ is the narrow two point integral in \eqref{eq:Icore}. On the event \eqref{eq:globalendbarrier}, split $\mathcal I$ according to the weak and endpoint barriers. Equations \eqref{eq:stripcontribution} and \eqref{eq:strongcontribution} give $\mathcal I\ll \vR^2M^{-c_2}$ for some fixed $c_2>0$, outside an event of probability $O(M^{-c_2})$. The common prime, overlap prime and frequency tail errors are smaller than $\vR^2M^{-B}$ in expectation for every fixed $B$.

It remains to return from the smoothed coefficients to the original interval. Put $D=(\cQ_z-\cQ_{z,\eta})- \mathbb E[\cQ_z-\cQ_{z,\eta}\mid\cF_z]$. Conditional Minkowski gives
\[
 \sqrt{\cE_z}
 \le \sqrt{\cE_{z,\eta}}
     +\sqrt{\mathbb E[D^2\mid\cF_z]}.
\]
By \cref{lem:smoothing},
$\|\cQ_z-\cQ_{z,\eta}\|_2=o(\vR M^{-B})$ for every fixed $B$ after increasing the numerical smoothing exponent if necessary. Since conditional expectation is an $L^2$ contraction,
$\mathbb E D^2=o(\vR^2M^{-B})$. For any fixed $A>0$, apply Markov's inequality to the nonnegative variable $\mathbb E[D^2\mid\cF_z]$ at level $\vR^2M^{-2A}$. Choosing $B>3A$, we obtain
\[
 \PP\left(\sqrt{\mathbb E[D^2\mid\cF_z]}>\vR M^{-A}\right)
 \le\frac{M^{2A}\mathbb E D^2}{\vR^2}=o(M^{-A}).
\]
Thus the second term in the conditional Minkowski bound is negligible at any fixed power of $M^{-1}$. Combining the estimates and decreasing $c$ if necessary proves \eqref{eq:QVgoal}.
\end{proof}

This proves \eqref{eq:fluctuation-target}, the last estimate required by the martingale criterion.

\subsection{Proof of Theorem~\ref{thm:main}}\label{sec:completion}
\begin{proof}
The smooth part estimate and \cref{cor:increments} verify
\eqref{eq:smooth-target} and \eqref{eq:increments-target}. By \cref{thm:CV}, $V_f/v_R\to1$
in probability, and \cref{prop:QV} gives $\cE_z/v_R^2\to0$ in probability.
The martingale assembly criterion, \cref{thm:assembly}, proves
\eqref{eq:targetclt}. The arguments for \cref{cor:fullpicture,rem:transition}
were given in the introduction.
\end{proof}

\begin{remark}
The theorem is qualitative. If $T$ tends to infinity very slowly, the barrier comparison may have an error of size $\exp(-c(\log T)^{1/200})$, which is smaller than one but need not be bounded by a fixed power of $M^{-1}$. We therefore do not state a uniform Berry-Esseen rate.
\end{remark}

\section{Application to a fixed nonreal Dirichlet twist}\label{sec:twists}
Let $\chi$ be a fixed primitive nonreal Dirichlet character modulo $q$, and put $\psi=\chi^2$. Thus $\psi$ is nonprincipal. We first consider the range $x/\sqrt{\log x}\le y=o(x)$, with the parameters of Section~\ref{sec:conditional} and $x$ large enough that $z>q$. The proof uses the martingale decomposition and the weak and strong barrier arguments developed in Sections~\ref{sec:conditional} and~\ref{sec:conditionalCLT}. The smooth part estimate, the fourth moment bounds for individual increments, and the sieve and Mellin-Parseval reduction carry over with the coefficients $\chi(n)$.

The character twist simplifies the frequency dependence of the tilted moments. On the large prime blocks used in the Gaussian comparison, the additional terms at frequencies $2t$ and $t+s$ contain the nonprincipal character $\psi$. The prime number theorem makes these terms negligible, uniformly including zero frequency. We can therefore apply the relative Gaussian comparison on both half axes without a separate treatment near zero. The correlation arising from $t-s$ remains, so the second moment argument still requires the restart estimate.

The complex Gaussian limit requires an additional quadratic estimate. To state the two estimates needed, for $p>z$ put $A_{p,\chi}=\sum_{x/p<n\le(x+y)/p,\,P(n)<p}f(n)\chi(n)$ and $\Delta_p^\chi=\varepsilon_p\chi(p)A_{p,\chi}$. The largest prime decomposition gives $S_{f,\chi}=S_{\chi,\le z}+S_{\chi,>z}$, where $S_{\chi,\le z}$ is the sum over $P(n)\le z$ and $S_{\chi,>z}=\sum_{p>z}\Delta_p^\chi$. Write $Q_\chi=\sum_{p>z}|A_{p,\chi}|^2$ and $R_\chi=\sum_{p>z}\chi(p)^2A_{p,\chi}^2$. For every real $\theta$,
\begin{equation}\label{eq:twistprojection}
 \sum_{p>z}\bigl(\Re(e^{-i\theta}\Delta_p^\chi)\bigr)^2
 =\frac12\bigl(Q_\chi+\Re(e^{-2i\theta}R_\chi)\bigr).
\end{equation}
We shall prove $Q_\chi/v_\chi\to1$ and $R_\chi/v_\chi\to0$ in probability. The second estimate ensures that every real projection has the same limiting variance. When $p_x\to0$, both estimates must hold on the smaller scale $v_\chi=c_qyp_x$. There is also a change in the frequency localization: the absolute value of the twisted Euler product is no longer even in the frequency, so we use conjugation and changes of variables in place of the reflection used in Section~\ref{sec:conditionalCLT}.

Section~\ref{sec:twistvariance} proves concentration of the conditional mean $V_\chi=\mathbb E[Q_\chi\mid\mathcal F_z]=\mathbb E[|S_{\chi,>z}|^2\mid\mathcal F_z]$. Section~\ref{sec:twistquadratic} then establishes $Q_\chi-V_\chi=o_{\mathbb P}(v_\chi)$ by adapting the localization and endpoint barrier estimates of Section~\ref{sec:conditionalCLT}. For $R_\chi$, the six edge parametrization places the nonprincipal character weights on the two variables carrying the largest prime conditions, and their cancellation gives the required bound. These estimates prove the complex Gaussian limit in the near full range. Section~\ref{sec:twistshort} treats the shorter intervals by a weighted martingale argument and the counting estimates used in the classical range, completing the proof of Theorem~\ref{thm:twisted}.

\subsection{Concentration of the conditional variance}\label{sec:twistvariance}
Put $F_{z,\chi}(s)=\prod_{p\le z,\,p\nmid q}(1+\varepsilon_p\chi(p)p^{-s})$ and $Z_{z,\chi}=\prod_{p\le z,\,p\nmid q}(1+p^{-1})$. Orthogonality gives $\mathbb E|F_{z,\chi}(\tfrac12+it)|^2=Z_{z,\chi}$ for every real $t$, and $d_zZ_{z,\chi}\to c_q$. The following estimate controls the additional frequency terms in the tilted moments.

\begin{lemma}\label{lem:charprime}
Let $\psi$ be a fixed nonprincipal Dirichlet character. There is $c=c(\psi)>0$ such that, for $v$ sufficiently large and real $\xi$,
\begin{equation}\label{eq:charblock}
 \sum_{v<p\le v^e}\frac{\psi(p)p^{-i\xi}}p
 \ll_\psi (1+|\xi|)e^{-c\sqrt{\log v}}.
\end{equation}
Moreover, if $z\le w\le3x$ and $|\xi|\le2W_0$, then for every fixed $A>0$,
\begin{equation}\label{eq:charEuler}
 E_w^\psi(\xi):=
 \prod_{z<p\le w}\left(1+\frac{\psi(p)}{p^{1+i\xi}}\right)
 =1+O_{\psi,A}((\log x)^{-A}).
\end{equation}
\end{lemma}
\begin{proof}
The prime number theorem for a fixed nonprincipal character gives $\sum_{p\le u}\psi(p)\log p\ll_\psi u e^{-c\sqrt{\log u}}$. Partial summation proves \eqref{eq:charblock}; the factor $1+|\xi|$ comes from differentiating $u^{-1-i\xi}$. For \eqref{eq:charEuler}, apply the same argument starting at $z=x^{1/\log\log\log x}$. Since $W_0=(\log x)^{2000}$, its contribution is smaller than every fixed negative power of $\log x$. The higher terms in the logarithm of the Euler product are $O(1/z)$. See \cite[Section~11.3]{MV}.
\end{proof}

Proposition~\ref{prop:transfer} applies with the bounded coefficients $\chi(n)$: its sieve estimates, square-free correction and Mellin-Parseval argument give, with probability $1-O(M^{-10})$,
\begin{equation}\label{eq:twistedtransfer}
 V_\chi=d_z\frac{xe^{1/T}}{2\pi T^2}
 \int_{\mathbb R}|F_{z,\chi}(\tfrac12+it)|^2K_T(t)^2\,dt
 +O(yM^{-50}).
\end{equation}
The error is $o(v_\chi)$, since $v_\chi\gg_q y/\sqrt M$. We next apply the phase comparison to prove concentration of the integral.

\begin{proposition}\label{prop:twistedV}
Suppose $x/\sqrt{\log x}\le y=o(x)$. Then
\begin{equation}\label{eq:twistedV}
 \frac{V_\chi}{v_\chi(x,y)}\xrightarrow{\mathbb P}1.
\end{equation}
\end{proposition}
\begin{proof}
For large $x$ we have $w>q$. Let $R_{j,\chi}(t)$ be the product $R_j(t)$ with $\varepsilon_p$ replaced by $\varepsilon_p\chi(p)$, and use the weak and rounded strong barriers for these products. Write $J_{G,\chi}$ and $J_{H,\chi}$ for the corresponding restricted integrals over $[-U,U]$. The four moment bounds for $R_j,R_j^{-1}$ and their derivatives in the proof of Proposition~\ref{lem:weakbarrier} hold for $R_{j,\chi}$ as well. Indeed, the orthogonality bound and the inverse factor bounds there apply with $a_p=\chi(p)p^{-1/2-it}$, using $|a_p|=p^{-1/2}$ and $|a'_p|=(\log p)p^{-1/2}$. The covering of $[-U,U]$ still uses $O(U\log w_j)$ intervals in block $j$. Hence the weak barrier holds on all of $[-U,U]$ with probability $1-o(1)$.

Take $\omega_p=\chi(p)$ in Proposition~\ref{prop:phasecomparison}. By Lemma~\ref{lem:charprime}, each $\Gamma_j(\xi)$ at $\xi=2t,2s,t+s$ is $O_q((1+U)e^{-ce^{(b+j-1)/2}})$ for $|t|,|s|\le U$, including when $\xi=0$. Lemma~\ref{lem:primeosc} applies to $t-s$ of either sign. Thus \eqref{eq:etaj} holds for one frequency on $[-U,U]$ with $d_*=1$ and for two distinct frequencies with $d_*=\min\{1,|t-s|\}$. The rectangular approximation error is unchanged. In particular, \eqref{eq:onepoint} and \eqref{eq:twopoint} hold on the signed interval for $R_{m,\chi}$, with relative error $O_q(e^{-b/2})$. The reference Gaussian paths are the same as before. The one point comparison for the weak barrier and Lemma~\ref{lem:ballot} give
\[
 \mathbb E(J_{G,\chi}-J_{H,\chi})=o(TZ_{z,\chi}p_{z,T}).
\]

For separated frequencies we use $g(t)=K_T(t)^2\ind_{[-U,U]}(t)$. The derivative bound in the proof of Lemma~\ref{lem:kernel} and the two endpoint terms show that $g$ satisfies the hypotheses of Lemma~\ref{lem:smallvariance}. In \eqref{eq:walshsmall} for $F_{w,\chi}$ the normalizing factor is $Z_{w,\chi}$, and each term with $rs=d$ carries the phase $\chi(r)\overline{\chi(s)}$, with $r,s$ supported on primes not dividing $q$. The Cauchy-Schwarz inequality bounds the squared Walsh coefficients as before, so the small prime integral has relative variance $O_q(e^{5b}/T)$.

For nearby pairs with $0<d=|t-s|<r_0$, keep the two cutoffs in \eqref{eq:cutpoints}. On a nonempty middle segment, divide the two point normalizing factor by $\prod_{w_{k_0}<p\le w_{k_1}}(1+1/p)^2$. By the calculation in Proposition~\ref{prop:phasecomparison}, the logarithm of this ratio is $O_q(1)$. Indeed, the sum of $D_j(t-s)$ is $O(1)$ since $d\log w_{k_0}\gg1$, and the sums of $|\Gamma_j(t+s)|$ and $u_j$ are bounded. Hence the middle segment contributes $O_q(e^{2(k_1-k_0)})$. On the last segment, for $k=k_1<m$, the weighted phase error in \eqref{eq:phaserho} is
\[
 \ll_{q,B_0}\frac{(a+k+1)^2}{d\log w_k}+O_q(e^{-b})=o(1)
\]
by the same summation as in \eqref{eq:restarterror}. The two Gaussian paths after conditioning on the prefix are independent, so \eqref{eq:restart} holds with both factors on this signed interval. If $k=m$ the segment is empty. The prefix bound follows from the one point comparison for either sign. Finally, the Cauchy-Schwarz inequality and the single factor fourth moment bound $1+6/p+O(p^{-2})$ give $\mathbb E|F_{w,\chi}(\tfrac12+it)F_{w,\chi}(\tfrac12+is)|^2\ll_q(\log w)^6$. We obtain \eqref{eq:nearpoint} and can integrate in $d$ using \eqref{eq:nearseries}, since $\int_{\mathbb R}K_T(t)^2K_T(t\pm d)^2\,dt\ll T$. The near pairs contribute $O_q(e^{-A}(TZ_{z,\chi}p^\#)^2)$.

The separated and nearby bounds give concentration of $J_{H,\chi}$ about $Z_{z,\chi}p^\#\int g$ as in Proposition~\ref{prop:concentration}. The mean loss from the weak to the strong barrier is negligible, and the weak barrier holds simultaneously on $[-U,U]$. Applying Markov's inequality to $J_{G,\chi}-J_{H,\chi}$ and using $p^\#\sim p_{z,T}$ gives the same concentration for the unrestricted integral on $[-U,U]$.

For the tails, set $R=2^jU$ as in Section~\ref{sec:CVcompletion}. When $\log R<\sqrt M$, run the weak barrier at scale $R$ on each of $[R,2R]$ and $[-2R,-R]$. The covering bounds are the same on each interval and the tilted moment estimates just used are uniform for $|t|\le2R$. Thus \eqref{eq:tailbarrier} holds on each interval with $Z_z$ replaced by $Z_{z,\chi}$. When $\log R\ge\sqrt M$, use $\mathbb E|F_{z,\chi}(\tfrac12+it)|^2=Z_{z,\chi}$ on both intervals. Summing the failure probabilities and the weighted means as in Section~\ref{sec:CVcompletion} gives $o_{\mathbb P}(TZ_{z,\chi}p_{z,T})$ for the two tails together. Since $\int_{-U}^U K_T(t)^2\,dt\sim\int_{\mathbb R}K_T(t)^2\,dt$, we obtain
\[
 \frac{\displaystyle\int_{\mathbb R}|F_{z,\chi}(\tfrac12+it)|^2K_T(t)^2\,dt}
 {\displaystyle Z_{z,\chi}p_{z,T}\int_{\mathbb R}K_T(t)^2\,dt}
 \xrightarrow{\mathbb P}1.
\]
Now use \eqref{eq:twistedtransfer}, the normalization above and \eqref{eq:scalecomparison}.
\end{proof}

\subsection{The two quadratic terms and the near full range}\label{sec:twistquadratic}
We first show that $Q_\chi-V_\chi=o_{\mathbb P}(v_\chi)$. Define $a_\chi(r)=\sum_{x/r<n\le(x+y)/r,\,P(n)\le z}f(n)\chi(n)$ for $r\in\mathcal A_z$. The rough expansion is
\begin{equation}\label{eq:Qchiexpand}
 Q_\chi=
 \sum_{\substack{r_1,r_2\in\mathcal A_z\\P(r_1)=P(r_2)}}
 a_\chi(r_1)\overline{a_\chi(r_2)}
 \chi(r_1)\overline{\chi(r_2)}
 f\!\left(\frac{r_1r_2}{(r_1,r_2)^2}\right).
\end{equation}
The diagonal is $V_\chi$. Expanding $\mathcal E_\chi=\mathbb E[|Q_\chi-V_\chi|^2\mid\mathcal F_z]$ gives the same square product condition as in Section~\ref{sec:sixedge}. Under the parametrization \eqref{eq:edges}, the character factor becomes
\begin{equation}\label{eq:Qcharchar}
 \chi(r_1)\overline{\chi(r_2)}\,
 \overline{\chi(r_3)}\chi(r_4)
 =\psi(b)\overline{\psi(c)}.
\end{equation}
Thus two of the cross factors carry nonprincipal characters. This changes the second kernel estimate, while the first remains valid in absolute value.

\begingroup
\setlength{\abovedisplayskip}{8pt plus 2pt minus 2pt}
\setlength{\belowdisplayskip}{8pt plus 2pt minus 2pt}
\begin{lemma}\label{lem:Qchikernel}
Let $G_{Q,\chi}(\mathbf t)$ be the rough kernel obtained by the smoothing and Mellin inversion of Section~\ref{sec:sixedge}, applied to $\mathcal E_\chi$. For $|t_j|\le W_0$,
\begin{equation}\label{eq:Qkernel1}
 |G_{Q,\chi}(\mathbf t)|
 \ll_q u^6\gamma(t_1+t_2)\gamma(t_3+t_4)+O_q(u^6/z).
\end{equation}
For every fixed $A>0$, we also have
\begin{equation}\label{eq:Qkernel2}
 \begin{split}
 |G_{Q,\chi}(\mathbf t)|
 &\ll_q u^4\{\gamma(t_1+t_3)+\gamma(t_2+t_4)\}\\
 &\qquad+O_{q,A}(u^4(\log x)^{-A})+O_q(u^6/z).
 \end{split}
\end{equation}
\end{lemma}
\begin{proof}
Fixing $a,b,c,d$ and summing $g,h$ gives \eqref{eq:Qkernel1}, exactly as in Lemma~\ref{lem:roughkernel}, since the character factors have modulus one. For the second bound, fix $g,h$ and put $P_*=\min\{P(g),P(h)\}$. The cross factor sum is
\[
 E_{P_*^-}(t_1+t_3)
 E_{P_*^-}^{\psi}(t_1+t_4)
 E_{P_*^-}^{\overline\psi}(t_2+t_3)
 E_{P_*^-}(t_2+t_4)-1.
\]
By \eqref{eq:charEuler}, the two character factors are $1+O_{q,A}((\log x)^{-A})$. Telescope the remaining two factors and use \eqref{eq:Eulerosc}. Summing $g,h$ absolutely costs $O(u^2)$, giving \eqref{eq:Qkernel2}. Shared prime factors and the removal of coprimality contribute $O_q(u^6/z)$ as in Lemma~\ref{lem:roughkernel}.
\end{proof}

\endgroup

\paragraph{Frequency localization.}
We explain the change needed in Section~\ref{sec:localization}, since $|F_{z,\chi}|$ need not be even. The four Mellin factors have character types $(\chi,\overline\chi,\overline\chi,\chi)$. Set $(u_1,u_2,u_3,u_4)=(t_1,-t_2,-t_3,t_4)$. The identity $F_{z,\overline\chi}(\tfrac12-iu)=\overline{F_{z,\chi}(\tfrac12+iu)}$, together with $|W(-u)|=|W(u)|$, puts all four absolute values in terms of $F_{z,\chi}(\tfrac12+iu_j)$. The arguments in \eqref{eq:Qkernel1} become $u_1-u_2$ and $u_4-u_3$; those in \eqref{eq:Qkernel2} become $u_1-u_3$ and $u_4-u_2$. We can therefore localize by differences directly on $[-W_0,W_0]^4$.

The required half moment bound holds uniformly for $|t|,|s|\le W_0$:
\begin{equation}\label{eq:twistedhalfmoment}
 \mathbb E|F_{z,\chi}(\tfrac12+it)F_{z,\chi}(\tfrac12+is)|
 \ll_q\Lambda^{1/2}M^{20}
 \left(1+\min\{\Lambda,|t-s|^{-1}\}^{1/2}\right).
\end{equation}
To verify this, put $D_{z,q}(v)=\sum_{p\le z,\,p\nmid q}\cos(v\log p)/p$ and $\Gamma_{z,\psi}(v)=\Re\sum_{p\le z}\psi(p)p^{-1-iv}$. The local expansion in \eqref{eq:halfmomentexpansion} now gives
\[
 \begin{aligned}
 \log\mathbb E|F_{z,\chi}(\tfrac12+it)F_{z,\chi}(\tfrac12+is)|
 ={}&\tfrac12D_{z,q}(0)+\tfrac12D_{z,q}(t-s)
       +\tfrac12\Gamma_{z,\psi}(t+s)\\
    &-\tfrac14\Gamma_{z,\psi}(2t)
       -\tfrac14\Gamma_{z,\psi}(2s)+O_q(1).
 \end{aligned}
\]
Take $v_0=\exp(C_\chi M^2)$ with $C_\chi$ large. Partial summation as in Lemma~\ref{lem:charprime} bounds the sum over $v_0<p\le z$ by $O_q((1+|\xi|)e^{-c\sqrt{\log v_0}})=O_q(1)$ for $|\xi|\le2W_0$. The remaining primes contribute at most $2\log M+O_q(1)$. Thus $|\Gamma_{z,\psi}(\xi)|\le2\log M+O_q(1)$, and the bounds for $D_{z,q}$ from Lemma~\ref{lem:halfmoment} prove \eqref{eq:twistedhalfmoment}.

With $w(t)=\min\{y/x,(1+|t|)^{-1}\}$ on $\mathbb R$, the bound $\int_{\mathbb R}w(t)w(t+h)\,dt\ll y/x$ gives the same $J_0,J_1,J_2$ estimates as in Section~\ref{sec:localization}. The two kernel bounds and Proposition~\ref{prop:localization} leave
\begin{equation}\label{eq:twistednarrow}
 \mathcal I_\chi=\frac{x^2u^6}{\Lambda^2}
 \iint_{\substack{|t|,|s|\le W_0\\|t-s|\le Cr_2}}
 |F_{z,\chi}(\tfrac12+it)F_{z,\chi}(\tfrac12+is)|^2
 w(t)^2w(s)^2\,dt\,ds.
\end{equation}
Here $r_2=M^{1000}/\Lambda$ is unchanged. The error $u^4(\log x)^{-A}$ in \eqref{eq:Qkernel2} is negligible after choosing $A$ large, by the crude fourth moment bound and $v_\chi^2\gg_q y^2/M$.

\paragraph{The endpoint bound.}
Use the partition and events in \eqref{eq:endbarriers} with $F_{z,\chi}$ in place of $F_z$.

\begin{lemma}\label{lem:twistedendpoint}
The weak and endpoint strip estimates \eqref{eq:endpointmean} and \eqref{eq:endpointstrip} hold with $F_z,Z_z$ replaced by $F_{z,\chi},Z_{z,\chi}$, uniformly for $|t|\le W_0$ and $a\ge1$, with implied constants depending on $q$.
\end{lemma}
\begin{proof}
On this partition $\log w_0\asymp M^{10}$. The moment calculation in the proof of Proposition~\ref{prop:phasecomparison}, together with Lemma~\ref{lem:charprime}, gives cumulative tilted mean $\mu_j=j+O_q(1)$ and variance $V_j=j/2+O_q(1)$, uniformly for $|t|\le W_0$. The block variances are bounded above and away from zero. The centered third prime moments are $O(p^{-3/2})$, giving standardized Berry-Esseen errors $\ll_q\exp(-cM^{10}e^{j-1})$. Proposition~\ref{prop:transitionstrip} therefore applies with $k=0$. Since $J\sim M$, it gives both stated probability bounds. Restoring the terminal tilt and the mean of the independent small prime product multiplies them by $Z_{z,\chi}$.
\end{proof}

\begin{proposition}\label{prop:twistedQ}
There exists $c>0$ such that
\[
 \mathbb P(\mathcal E_\chi>v_\chi(x,y)^2M^{-c})\ll_q M^{-c}.
\]
Consequently,
\begin{equation}\label{eq:QminusVchi}
 \frac{Q_\chi-V_\chi}{v_\chi(x,y)}\xrightarrow{\mathbb P}0.
\end{equation}
\end{proposition}
\begin{proof}
The preceding localization leaves \eqref{eq:twistednarrow}. Put $T_0=x/y$ and $a_t=\tfrac12\log(T_0+|t|)$. The derivative and Sobolev argument for \eqref{eq:globalendbarrier} applies on both half axes with $F_{z,\chi}$, so the weak endpoint barrier holds throughout $[-W_0,W_0]$ outside an event of probability $O_q(M^{-20})$.

Apply Lemma~\ref{lem:twistedendpoint} with $a=a_t$. Since the deterministic weights depend on $|t|$, the strip and terminal estimates \eqref{eq:stripcontribution} and \eqref{eq:strongcontribution} extend to the full axis with only fixed factors. They retain $s_0\asymp p_x$ and give a fixed power saving for $\mathcal I_\chi$ on the scale $v_\chi^2\asymp_q y^2p_x^2$. The common prime errors, smoothing and frequency tails satisfy the bounds of Section~\ref{sec:smoothing}, since $|\chi(n)|\le1$. This proves the asserted estimate for $\mathcal E_\chi$. For fixed $\epsilon>0$, condition on $\mathcal F_z$ and apply Chebyshev's inequality to $Q_\chi-V_\chi$. Averaging and separating the exceptional event in the proposition, we obtain
\[
 \mathbb P(|Q_\chi-V_\chi|>\epsilon v_\chi)
 \le\mathbb P(\mathcal E_\chi>v_\chi^2M^{-c})
       +\epsilon^{-2}M^{-c}
 \ll_{q,\epsilon}M^{-c}.
\]
This proves \eqref{eq:QminusVchi}.
\end{proof}

It remains to prove that $R_\chi$ is negligible on the same scale. Its rough expansion is
\[
 R_\chi=
 \sum_{\substack{r_1,r_2\in\mathcal A_z\\P(r_1)=P(r_2)}}
 a_\chi(r_1)a_\chi(r_2)\chi(r_1r_2)
 f\!\left(\frac{r_1r_2}{(r_1,r_2)^2}\right).
\]
In the expansion of $\mathbb E|R_\chi|^2$, the same square product parametrization gives
\begin{equation}\label{eq:Rcharchar}
 \chi(r_1r_2)\overline{\chi(r_3r_4)}
 =\psi(g)\overline{\psi(h)}.
\end{equation}
The nonprincipal weights now lie on $g,h$, the factors carrying the largest prime conditions. Their cancellation bounds the whole rough kernel.

\begin{proposition}\label{prop:Rchismall}
Suppose $x/\sqrt{\log x}\le y=o(x)$. Then
\begin{equation}\label{eq:Rsmall}
 \frac{R_\chi}{v_\chi(x,y)}\xrightarrow{\mathbb P}0.
\end{equation}
\end{proposition}
\begin{proof}
With the multiplier $\psi_\eta$ from Section~\ref{sec:sixedge}, define
\[
 A_{p,\chi,\eta}
 =\sum_{\substack{n\ge1\\P(n)<p}}
       f(n)\chi(n)\psi_\eta(pn/x),
\]
and put $R_{\chi,\eta}=\sum_{p>z}\chi(p)^2A_{p,\chi,\eta}^2$.
We apply the Mellin expansion of Section~\ref{sec:sixedge} to $\mathbb E|R_{\chi,\eta}|^2$, and denote the resulting rough kernel by $G_{R,\chi}(\mathbf t)$. As in Lemma~\ref{lem:roughkernel}, the terms with $\rho>1$ and the removal of edge coprimality contribute $O_q(u^6/z)$. In the remaining sum, fix $a,b,c,d$ and put $D=\max\{z,P(abcd)\}$. By \eqref{eq:Rcharchar}, the $g$ sum is
\[
 \sum_{\substack{g\in\mathcal A_z\\D<P(g)\le3x}}
       \frac{\psi(g)}{g^{1+i(t_1+t_2)}}
 =E_{3x}^{\psi}(t_1+t_2)-E_D^{\psi}(t_1+t_2).
\]
The $h$ sum is $E_{3x}^{\overline\psi}(t_3+t_4)-E_D^{\overline\psi}(t_3+t_4)$. For $|t_j|\le W_0$, both frequencies have absolute value at most $2W_0$. Estimate \eqref{eq:charEuler} holds uniformly for every cutoff $z\le D\le3x$, so each difference is $O_{q,A}((\log x)^{-A})$ for every fixed $A>0$. Thus the maximum prime restriction preserves the character cancellation. Summing $a,b,c,d$ absolutely costs $O(u^4)$, and, since $A$ is arbitrary, we obtain
\[
 |G_{R,\chi}(\mathbf t)|
 \ll_{q,A}u^4(\log x)^{-A}+u^6/z.
\]
To estimate $\mathbb E|R_{\chi,\eta}|^2$, apply H\"older's inequality to the four Mellin factors. Using $\mathbb E|F_{z,\chi}(\tfrac12+it)|^4\ll(\log z)^6$ and $\int_{\mathbb R}|W(t)|\,dt\ll\log\log x$ with the preceding kernel bound, we obtain
\[
 \mathbb E|R_{\chi,\eta}|^2=o(v_\chi^2)
\]
for $A$ sufficiently large. Here we use $v_\chi^2\gg_q y^2/M$ and $x/y\le\sqrt{\log x}$. Apply Markov's inequality to $|R_{\chi,\eta}|^2$, then remove the frequency tails and smoothing as in Section~\ref{sec:smoothing}, with the inserted coefficients bounded by one. This proves \eqref{eq:Rsmall}.
\end{proof}

For the near full range, Lemmas~\ref{lem:smooth} and~\ref{lem:fourth} make the smooth part negligible and give $\sum_{p>z}\mathbb E|\Delta_p^\chi|^4=o(v_\chi^2)$. For each fixed $\theta$, put $X_p=\sqrt{2/v_\chi}\,\Re(e^{-i\theta}\Delta_p^\chi)$. Proposition~\ref{prop:twistedV}, \eqref{eq:QminusVchi} and \eqref{eq:Rsmall}, together with \eqref{eq:twistprojection}, imply $\sum_pX_p^2\to1$ in probability. As in the proof of Proposition~\ref{thm:assembly}, we obtain $\mathbb E\max_p|X_p|\to0$ by H\"older's inequality. The real martingale central limit theorem and Cram\'er-Wold therefore yield
\[
 \frac{S_{f,\chi}(x,y)}{\sqrt{v_\chi(x,y)}}
 \xrightarrow{d}\mathcal N_{\mathbb C}(0,1).
\]
This proves Theorem~\ref{thm:twisted} for $x/\sqrt{\log x}\le y=o(x)$.

\subsection{Shorter intervals and completion of the proof}\label{sec:twistshort}
It remains to treat $x^{1/5}\log x\ll y\le x/\sqrt{\log x}$. In this range $\kappa_x\to\infty$ and $p_x\to1$, so it suffices to prove the complex Gaussian limit with normalization $\sqrt{c_qy}$. We use the counting estimates of Chatterjee and Soundararajan and of Soundararajan and Xu, together with the following weighted form of the real martingale argument.

For a finite set $\mathcal A$ of square-free integers greater than $1$, let $E_{\mathrm{off}}(\mathcal A)$ count the ordered quadruples in $\mathcal A^4$ whose product is a square and whose entries cannot be paired into equal pairs.

\begin{lemma}\label{lem:bounded-real-weights}
Let $\mathcal A$ and deterministic real coefficients $b_n$, $n\in\mathcal A$, vary along a sequence, with $|b_n|\le1$. Put $V=\sum_{n\in\mathcal A}b_n^2$ and $B_p=\sum_{n\in\mathcal A,\,P(n)=p}b_n^2$. Suppose that $V\to\infty$, $E_{\mathrm{off}}(\mathcal A)=o(V^2)$ and $\max_pB_p=o(V)$. Then
\[
 \frac1{\sqrt V}\sum_{n\in\mathcal A}b_nf(n)
 \xrightarrow{d}N(0,1).
\]
\end{lemma}
\begin{proof}
Reveal the prime signs in increasing order. The variables
\[
 D_p=\frac1{\sqrt V}\sum_{\substack{n\in\mathcal A\\P(n)=p}}b_nf(n)
\]
are martingale differences, and $Q=\sum_pD_p^2$ has mean $1$. Rademacher orthogonality leaves only square products. The three diagonal pairings in $\sum_p\mathbb E D_p^4$ cost at most $3\sum_pB_p^2/V^2$. In $\mathbb E Q^2$, the pairing $n_1=n_2$, $n_3=n_4$ contributes $1$; either other pairing forces a common largest prime factor and costs at most $\sum_pB_p^2/V^2$. Bounding non-diagonal coefficients by $1$, we obtain
\begin{align*}
 \sum_p\mathbb E D_p^4
 &\le\frac{3\sum_pB_p^2+E_{\mathrm{off}}(\mathcal A)}{V^2}=o(1),\\
 \mathbb E(Q-1)^2
 &\le\frac{2\sum_pB_p^2+E_{\mathrm{off}}(\mathcal A)}{V^2}=o(1),
\end{align*}
since $\sum_pB_p^2\le(\max_pB_p)V=o(V^2)$. Apply the real martingale central limit theorem, \cite[Theorem~2.3]{SX}.
\end{proof}

We next compute the variance of each real projection. The required square-free counts are as follows.

\begin{lemma}\label{lem:squarefree-fixed-q}
Let $\chi$ be a fixed nonreal Dirichlet character modulo $q$. If $x^{1/5}\log x\ll y\le x$, then
\begin{equation}\label{eq:twisted-density}
 \sum_{x<n\le x+y}\mu^2(n)|\chi(n)|^2=c_qy+o_q(y),
\end{equation}
and
\begin{equation}\label{eq:twisted-pseudovariance}
 \sum_{x<n\le x+y}\mu^2(n)\chi(n)^2=o_q(y),
\end{equation}
where $c_q$ is defined before Theorem~\ref{thm:twisted}.
\end{lemma}
\begin{proof}
We obtain the square-free count in each reduced residue class from the absolute large square estimate in \cite[Section~2.1]{PandeySQ}. For some sufficiently small fixed $0<\eta<1/100$, uniformly for $h\asymp X^{1/5-\eta}$, put $D=hX^{-\eta}$. To estimate the contribution of $d>D$, apply equations~(2) and~(6) there on each dyadic range and sum the estimates. We obtain
\[
 R(X,h;D):=\sum_{D<d\le\sqrt{X+h}}
       \#\{m\ge1:X<md^2\le X+h\}
 \ll hX^{-\eta}\log X=o(h).
\]
Here $D\gg\sqrt h$, so each $d$ contributes at most one $m$. This unsigned estimate remains valid after imposing a congruence condition. The identity $\mu^2(n)=\sum_{d^2\mid n}\mu(d)$ and the Chinese remainder theorem give, for $(a,q)=1$,
\begin{align*}
 \sum_{\substack{X<n\le X+h\\n\equiv a\pmod q}}\mu^2(n)
 &=\sum_{\substack{d\le D\\(d,q)=1}}\mu(d)
       \left(\frac{h}{qd^2}+O(1)\right)+O(R(X,h;D))\\
 &=\frac hq\prod_{p\nmid q}(1-p^{-2})
       +O_q\bigl(D+h/D+R(X,h;D)\bigr)\\
 &=\frac hq\prod_{p\nmid q}(1-p^{-2})+o_q(h).
\end{align*}
For general $y$ in the stated range, divide $(x,x+y]$ into $J=\lfloor y/x^{1/5-\eta}\rfloor$ equal intervals of length $h=y/J$. Then $x^{1/5-\eta}\le h\le2x^{1/5-\eta}$ for large $x$, and all left endpoints are comparable to $x$. Summing the uniform estimate gives the residue class asymptotic for $(x,x+y]$. Sum over the reduced classes with weights $1$ and $\chi(a)^2$. The latter main term vanishes because $\chi^2$ is nonprincipal, while the former has coefficient
\[
 \frac{\varphi(q)}q\prod_{p\nmid q}(1-p^{-2})
 =\frac6{\pi^2}\prod_{p\mid q}\frac p{p+1}=c_q.
\]
This proves both assertions.
\end{proof}

\begin{proof}[Completion of the proof of Theorem~\ref{thm:twisted}]
Let $\mathcal T$ be the square-free integers in $(x,x+y]$. We choose $\mathcal S\subseteq\mathcal T$ with $|\mathcal T\setminus\mathcal S|=o(y)$ and $E_{\mathrm{off}}(\mathcal S)=o(y^2)$.

If $y\le x/(\log x)^2$, take $\mathcal S=\mathcal T$. Writing $\delta_0=y/x$, Proposition~3.1 of \cite{CS} gives
\[
 E_{\mathrm{off}}(\mathcal T)
 \ll y^2\delta_0\log x(1+\delta_0\log x)=o(y^2).
\]
For $x/(\log x)^2<y\le x/\sqrt{\log x}$, take $\mathcal S=\{n\in\mathcal T:\Omega(n)\le(1+\epsilon)\log\log x\}$, where $\epsilon>0$ is fixed and sufficiently small. The proof of \cite[Corollary~9.2]{SX} gives the two required counting estimates throughout this range, with $\epsilon$ chosen so that $2\log2-1+3\epsilon<1/2$. In both ranges, \cite[Lemma~3.3]{SX} gives
\[
 \max_p\#\{x<n\le x+y:P(n)=p\}=o(y).
\]

For fixed $\theta\in\mathbb R$, put $b_n=\Re(e^{-i\theta}\chi(n))$. By Lemma~\ref{lem:squarefree-fixed-q} and the bound for $|\mathcal T\setminus\mathcal S|$,
\begin{align*}
 V_\theta:=\sum_{n\in\mathcal S}b_n^2
 &=\frac12\sum_{x<n\le x+y}\mu^2(n)|\chi(n)|^2
   +\frac12\Re\!\left(e^{-2i\theta}
       \sum_{x<n\le x+y}\mu^2(n)\chi(n)^2\right)+o(y)\\
 &=\frac{c_qy}{2}+o_q(y).
\end{align*}
The counting estimates verify the hypotheses of Lemma~\ref{lem:bounded-real-weights} for these coefficients on $\mathcal S$. The omitted terms are negligible in $L^2$, since
\[
 \mathbb E\left|\sum_{n\in\mathcal T\setminus\mathcal S}f(n)\chi(n)\right|^2
 =\sum_{n\in\mathcal T\setminus\mathcal S}|\chi(n)|^2=o(y).
\]
Thus $\sqrt{2/(c_qy)}\,\Re(e^{-i\theta}S_{f,\chi}(x,y))$ converges to $N(0,1)$ for every fixed $\theta$. The Cram\'er-Wold theorem yields
\begin{equation}\label{eq:twisted-classical-limit}
 \frac{S_{f,\chi}(x,y)}{\sqrt{c_qy}}
 \xrightarrow{d}\mathcal N_{\mathbb C}(0,1).
\end{equation}
Since $p_x\to1$ in this range, Slutsky's theorem replaces $c_qy$ by $v_\chi(x,y)=c_qyp_x$. Together with the near full range proved above, this completes Theorem~\ref{thm:twisted} throughout $x^{1/5}\log x\ll y=o(x)$.
\end{proof}

\appendix
\section{Probability details}\label{app:probability}
We record how the four estimates combine through the real martingale central limit theorem.

\begin{proposition}\label{thm:assembly}
Write $S=S_f(x,y)$. Suppose that, in the range of Theorem~\ref{thm:main}, the four conditions \eqref{eq:smooth-target}-\eqref{eq:fluctuation-target} hold.
Then \eqref{eq:targetclt} holds. More precisely, for each fixed $t\in\RR$,
\begin{align}\label{eq:mastererror}
 \left|\EE e^{itS/\sqrt{\vR}}-e^{-t^2/2}\right|
 \ll_t{}&
 \left(\frac{\EE S_{\le z}^2}{\vR}\right)^{1/2}
 +\mathcal L_x^{1/4}
 +\EE\min\left\{1,\left|\frac{\Vf}{\vR}-1\right|\right\}\notag\\
 &+\EE\min\left\{1,\frac{\sqrt{\cE_z}}{\vR}\right\}.
\end{align}
Here $\cQ_z=\sum_{p>z}\Delta_p^2$,
$\cE_z=\EE[(\cQ_z-\Vf)^2\mid\cF_z]$, and
$\mathcal L_x=\sum_{p>z}\EE\Delta_p^4/\vR^2$.
\end{proposition}
\begin{proof}
We condition on $\cF_z$ and apply the Cauchy-Schwarz inequality to $|\cQ_z-\Vf|$. Since the clipped variable is also at most $1$, we obtain
\[
 \EE\left[\min\left\{1,\frac{|\cQ_z-\Vf|}{\vR}\right\}\mid\cF_z\right]
 \le\min\left\{1,\frac{\sqrt{\cE_z}}{\vR}\right\}.
\]
To bound $\EE\min\{1,|\cQ_z/\vR-1|\}$, split $\cQ_z-\vR=(\cQ_z-\Vf)+(\Vf-\vR)$, apply $\min(1,a+b)\le\min(1,a)+\min(1,b)$ and then take expectations. Using the preceding conditional estimate, we bound this expectation by the last two terms in \eqref{eq:mastererror}. For the maximum increment, apply H\"older's inequality to obtain
\[
 \EE\max_{p>z}\frac{|\Delta_p|}{\sqrt{\vR}}
 \le\left(\frac{\sum_{p>z}\EE\Delta_p^4}{\vR^2}\right)^{1/4}.
\]
Proposition~2.2 of \cite{SX}, applied to $X_p=\Delta_p/\sqrt{\vR}$, now proves \eqref{eq:mastererror} for $S_{>z}$. Restoring $S_{\le z}$ costs at most $|t|(\EE S_{\le z}^2/\vR)^{1/2}$. Finally, for fixed $0<\epsilon<1$, split according to whether $|Y_x|\le\epsilon$ to obtain
\[
 \EE\min\{1,|Y_x|\}\le\epsilon+\PP(|Y_x|>\epsilon).
\]
If $Y_x\to0$ in probability, first let $x\to\infty$ and then $\epsilon\downarrow0$. Thus the bounded error expectations vanish. Under \eqref{eq:smooth-target}-\eqref{eq:fluctuation-target}, all terms in \eqref{eq:mastererror} tend to zero, proving the conclusion.
\end{proof}

\begin{remark}\label{rem:noUI}
Neither $\EE\cE_z=o(\vR^2)$ nor $\EE|\Vf/\vR-1|\to0$ is required. Replacing the clipped errors by untruncated expectations would impose stronger conditions that can be inappropriate in the near full regime.
\end{remark}

\section*{AI declaration}

The strategy for adapting the proof of Harper, Soundararajan and Xu from the Steinhaus model to the Rademacher model is due to the authors. We worked out the proofs by using relative Gaussian comparison and a parametrization by six edge variables. ChatGPT was used to help us read the literature, check the manuscript for errors, and polish the English presentation. We take full responsibility for the mathematical content of the paper.

\end{document}